\documentclass[11pt,longbibliography]{article}

\usepackage[utf8]{inputenc}

\usepackage[T1]{fontenc}

\usepackage[a4paper, margin=2.4cm]{geometry}

\usepackage{amsmath, amsthm, amsfonts, amssymb}
\usepackage{mathrsfs}           
\usepackage{enumitem}
\usepackage[colorlinks=true,linkcolor=blue,citecolor=blue,urlcolor=blue,breaklinks]{hyperref}

\usepackage{bbm} 

\numberwithin{equation}{section}

\usepackage{enumitem}
\setenumerate[1]{label=\thesection.\arabic*.}
\setenumerate[2]{label*=\arabic*.}

\usepackage{breakurl}
\usepackage{url}
\usepackage{doi}
\def\N{\mathbb{N}}
\def\R{\mathbb{R}}

\def\E{\mathbb{E}}
\def\P{\mathbb{P}}
\def\F{\mathcal{F}}
\def\mcC{\mathcal{C}}

\DeclareMathOperator{\supp}{supp}

\DeclareMathOperator{\ee}{e}

\def\la{\langle}
\def\ra{\rangle}

\newcommand{\LM}[1]{\textcolor{black}{#1}}

\newtheorem{thm}{Theorem}[section]
\newtheorem{corollary}[thm]{Corollary}
\newtheorem{lemma}[thm]{Lemma}
\newtheorem{proposition}[thm]{Proposition}
\newtheorem{conj}[thm]{Conjecture}

\theoremstyle{definition}
\newtheorem{dfn}[thm]{Definition}
\theoremstyle{remark}
\newtheorem{rem}[thm]{Remark}
\theoremstyle{example}

\title{Weak Existence and Uniqueness for Super-Brownian Motion with Irregular Drift }
\author{ Leonid Mytnik \and Johanna Weinberger}

\AtEndDocument{\bigskip{\footnotesize
		\noindent
		\textsc{Leonid Mytnik. The Faculty of Data and Decision Sciences, Technion - Israel Institute of Technology.}
		\textit{E-mail address}: \texttt{\href{mailto: leonidm@technion.ac.il}{ leonidm@technion.ac.il}}
		\par
		\addvspace{\medskipamount}
		\noindent
		\textsc{Johanna Weinberger. The Faculty of Data and Decision Sciences, Technion - Israel Institute of Technology.}
		\textit{E-mail address}: \texttt{\href{mailto:wjohanna@campus.technion.ac.il}{wjohanna@campus.technion.ac.il}}
}}

\begin{document}
	
	\maketitle
	
	\begin{abstract}
	\noindent We establish weak existence and uniqueness for random field solutions of the one-dimensional SPDE
		\begin{equation*}
			d_tX_t = \frac{1}{2}\Delta X_t +h(X_t)+ \sqrt{X_t}\dot{W}, \quad t\geq 0,
		\end{equation*}
		where $\dot{W}$ denotes space–time white noise and $h$ is a bounded drift with $h(0)\geq 0$. The proof relies on an extension of the duality relation for super-Brownian motion, which allows us to treat a broad class of admissible drifts, including functions that are non-Lipschitz or discontinuous at zero. In particular, well-posedness is derived for certain drifts that are H\"older continuous at zero with exponent $\alpha\in(0,1)$. We also allow discontinuous drifts of the form $h(x) = b_0\mathbbm{1}_{x = 0} + b_1\mathbbm{1}_{x>0},$
where $b_0 \geq  0$, $b_1 \in \R$. Additionally, if $h(0)=0$ and the initial condition is continuous and compactly supported, we show that the Lebesgue measure of the non-zero set of $X$ is finite.  The proofs are based on duality. We use a log-Laplace equation  which is perturbed by jump noise as the equation for the dual process and the jumps of the dual process are allowed to take infinite values. We believe that the results for the dual process are also of independent interest.  \par
\noindent Under suitable assumptions on $h$ we also prove survival of $X$ with positive probability, using rescaling and comparison to the KPP-equation with branching noise.

	\end{abstract}
	
	\section{Introduction}
	\label{sec:intro}
	We consider one-dimensional SPDEs of the form  
	\begin{equation}\label{eq:SPDE1}
		d_tX_t(x) = \frac{1}{2}\Delta X_t(x) +h(X_t(x))+ \sqrt{X_t(x)}\dot{W}(t,x), \quad X_0 = f, x\in \R, t\geq 0,
	\end{equation}
	where $f$ is a non-negative, continuous, bounded function,  $\dot{W}$ is space-time white noise and
	\begin{equation*}
		h(x) = h_1(x)+ h_{\infty}(x), \quad x\in \R_+.
	\end{equation*}
	Hereby, 
	\begin{equation}
		\label{eq:h1}h_{1}(x) = \int_0^\infty (\ee^{-\lambda x})(\nu^2-\nu^1)(d\lambda)  \quad \forall x \in \R_+,
	\end{equation}
	and
	\begin{equation}\label{eq:hinfty}
		h_{\infty}(x) = b_0\mathbbm{1}_{x = 0} + b_1\mathbbm{1}_{x>0}, \quad x\in \R_+,
	\end{equation}
for $b_0 \geq  \la \nu^1-\nu^2,1\ra$, $b_1 \in \R $ and finite measures $\nu^1, \nu^2$ on $(0,\infty)$. This in particular implies
\begin{equation}\label{eq:h0}
		h(0) \geq 0.
	\end{equation}
The above assumptions on $h$ include a wide range of examples that are irregular at zero (see Section~\ref{sec:examples}). \\
In this article,  using duality methods, we establish weak well-posedness for solutions of \eqref{eq:SPDE1} for drifts satisfying conditions (\ref{eq:h1} -- \ref{eq:h0}).
A solution to~\eqref{eq:SPDE1} with the drift $h(x)=cx$, $c\in \R$, describes the density of the classical super-Brownian motion (SBM) with drift in one dimension, which has been extensively studied in the literature (see, for instance, \cite{Watanabe68, Konno88, Dawson93, Dynkin94, LeGall99, Perkins02}), and for which weak well-posedness is well known. Our goal is to establish weak well-posedness for the much broader class of drifts satisfying conditions (\ref{eq:h1} -- \ref{eq:h0}), beyond  the classical case of $h(x)=cx$. We will comment further on this below.

We consider solutions in the space  $\mathcal{C}(\R_+,\mathcal{C}_{tem}^+)$, which is the space of continuous functions taking values in 
\begin{equation*}
	\mathcal{C}_{tem}^+ = \{ f \in \mathcal{C}(\R) \colon f\geq 0,  |f|_{(\lambda)}<\infty \quad \forall \lambda <0\},
\end{equation*}
where
\begin{equation*}
		|f|_{(\lambda)} = \sup_{x\in \R}|\ee^{\lambda|x|}f(x)|, \quad  f\in \mathcal{C}(\R).
\end{equation*}
We also allow random initial conditions distributed according to a measure in $\mathcal{P}_p(\mathcal{C}_{tem}^+)$ (see Section~\ref{sec:prelim}). 
We now state the main results of the paper. 
	\begin{thm}\label{th:uniqueness}
	Let  $\eta \in \mathcal{P}_p(\mathcal{C}_{tem}^+)$ with $p >2$ and assume $X_0 \sim \eta$ is such that
	\begin{equation}\label{eq:cond_x_0}
		\sup_{x\in \R}\E[X_0(x)] <\infty.
	\end{equation}
	Then there exists a unique weak solution to \eqref{eq:SPDE1} in $\mathcal{C}(\R_+,\mathcal{C}_{tem}^+)$ in the sense of Definition~\ref{def:weak_sol}.
\end{thm}
A natural question concerning the properties of the solution  is the  following: How large  can the non-zero set of the solution be in the case where there is no upward drift at the zero points of the solution? That is, in the case when 
$h(0)= 0$. In this setting, we prove that  the non-zero set of the solution to \eqref{eq:SPDE1} has finite Lebesgue measure.
\begin{thm}\label{thm:scp}
Suppose that almost surely, $X_0 = f \in \mathcal{C}_c^+(\R)$. Let $(X_t)_{t\geq 0}$ be the solution to \eqref{eq:SPDE1} with $h$ given by 
\begin{equation}\label{eq:h_special_supp}
	h(x) = \int_0^\infty(1-\ee^{-\lambda x})\nu^1(d\lambda) + b_1\mathbbm{1}_{x>0} , \quad x\in \R_+,
\end{equation}
 with $b_1 \in \R$ such that $h(0)\geq 0$. Then there exists $T_0>0$ such that for all $t \in [0,T_0]$ we have
\begin{equation}\label{eq:th_supp}
	\mathbb{P}(\mathrm{Leb}(\{ x\in \R \colon X_t(x) > 0\}) < \infty) = 1.
\end{equation}
\end{thm}

Further fundamental questions concern the size of the closure of the non-zero set, in particular whether the compact support property holds under the assumptions of Theorem~\ref{thm:scp}. We conjecture the following.
\begin{conj}
	Suppose that almost surely, $X_0 = f \in \mathcal{C}_c^+(\R)$. Let $(X_t)_{t\geq 0}$ be the solution to \eqref{eq:SPDE1} with $h$ given by \eqref{eq:h_special_supp}. Then $X_t$ is compactly supported with probability one for all $t\geq 0$.
\end{conj} 
Besides the non-zero set and the support of the solution it is natural to investigate the longtime behavior, in particular the survival probability of the process. 
\begin{thm}\label{thm:survival}
 	Suppose that $h\neq 0$ is given by
	\begin{equation*}
		h(x) = \int_0^\infty(1-\ee^{-\lambda x})\nu^1(d\lambda) + b_1\mathbbm{1}_{x>0} , \quad x\in \R_+,
	\end{equation*}
	 with $b_1\geq 0$ and let $(X_t)_{t\geq 0}$ be the solution to \eqref{eq:SPDE1} with $X_0 =f\in \mathcal{C}_c^+(\R)$ and $f\neq 0$. Then 
	 \begin{equation*}
	 	\mathbb{P}(\la X_t,1\ra >0, \quad \forall t > 0) >0. 
	 \end{equation*}
\end{thm}
The main tool for the proof of Theorems~\ref{th:uniqueness} and \LM{\ref{thm:scp}} is  the duality relation derived in~Proposition~\ref{prop:duality}.
 \\
\subsection{Admissible drift terms and comparison with scalar SDEs}
Our assumptions about the structure of the drift function $ h$ may seem peculiar, but they actually allow for significant irregularities of $h$ at zero, precisely the point where the noise coefficient $\sqrt{x}$ becomes non-Lipschitz and degenerate.
For example, by choosing
\begin{equation}
\label{eq:choice_disc}
\nu^1 = \nu^2 \equiv 0, \quad b_0 \ne b_1, \quad b_0 \geq 0,
\end{equation}
one obtains a drift $h(\cdot)$ that is {\it discontinuous} at zero. More precisely, this choice   of parameters $\nu^i, b_i, i=1,2$, in~\eqref{eq:choice_disc},
 allows the drift $h$
to take any constant value at points where the solution $X$
of~\eqref{eq:SPDE1} is positive, and any constant  non-negative value at points where 
$X$
 equals zero. Furthermore, with \( \nu^2 \equiv 0 \), for any \( \alpha \in (0,1) \), by choosing an appropriate measure \( \nu^1 \) and suitable constants \( b_0, b_1 \), one can obtain \( h(x) \) such that it behaves like \( c x^\alpha \) for small values of \( x \), in the sense that
\[
\lim_{x \downarrow 0} \frac{h(x)}{x^\alpha} = c > 0.
\]
To the best of our knowledge, weak uniqueness for such irregular drifts at zero has not been previously established. In fact, even the existence of a solution to equation~\eqref{eq:SPDE1} with a drift discontinuous at zero was not known.
It is worth noting that standard techniques for handling drifts, such as the Dawson–Girsanov theorem, can no longer be directly applied to establish weak well-posedness when the drift has a H\"older exponent $\alpha<1/2$  at zero or exhibits a  discontinuity at zero. \\
It is also noteworthy that we establish 
weak well-posedness of \eqref{eq:SPDE1}  for discontinuous drift terms $h$ that lead to non-uniqueness or non-existence of solutions for the one-dimensional counterpart of~\eqref{eq:SPDE1} given by
\begin{equation}
\label{eq:SDE1_nue}
	dx_t =  h(x_t)dt + \sqrt{x_t}dB_t, \quad x_0 = 0,  \quad t \geq 0, 
\end{equation}
where 
$(B_t)_{t\geq 0}$ 
is a standard Brownian motion. To be more precise, we show that this equation admits non-unique solutions when  $h(x)=\mathbbm{1}_{x>0}$ and no solution exists when
$h(x)=\mathbbm{1}_{x=0}$.  Notably, these are exactly the cases for which we prove weak well-posedness of \eqref{eq:SPDE1} in this paper. 
The above examples in SPDE and SDE settings are discussed in more detail  in Section~\ref{sec:examples}.

\subsection{SPDEs with irregular coefficients}
The question addressed in this paper belongs to a broader area of well-posedness of SPDEs with irregular coefficients. In particular the SPDE 
\begin{equation}\label{eq:SPDE_multi}
		d_tX_t(x) = \frac{1}{2}\Delta X_t(x) +h(X_t(x))+ \sigma(X_t(x))\dot{W}(t,x), \quad X_0 = f,\; x\in \R, t\geq 0,
	\end{equation}
where $\dot W$ is a space-time white noise, $h$ is a drift and $\sigma$ is a noise coefficient, has been extensively studied. 
Many results on well-posedness for this equation mainly dealt with the question:\\
\textit{How irregular can the drift be and what are conditions on the  noise coefficient so that 
well-posedness for this equation holds?}\\
This problem is also known under the name of regularization by noise, as it has been shown for some irregular $h$ that the equation without noise is ill-posed 
and when the noise is added, the equation becomes well-posed. 
When the  noise is additive ($\sigma\equiv 1$), strong well-posedness was derived for some irregular function-valued drifts in~\cite{GP93a}, \cite{GP93b}, and more recently even for some measure-valued and distribution-valued drifts in~\cite{ABLM}. Later the result from~\cite{ABLM} was extended in~\cite{D24} to SPDEs driven by multiplicative noises with $\sigma$ bounded away from zero (non-degenerate) and sufficiently regular.  Weak well-posedness was established in~\cite{BM24} for~\eqref{eq:SPDE_multi} with additive noise for more irregular distributions than
 in~\cite{ABLM}. The so-called path-by-path uniqueness for~\eqref{eq:SPDE_multi} with additive noise and with bounded drifts was derived in~\cite{BM19}.  

 \paragraph{} For SPDEs driven by multiplicative noise with
non-Lipschitz and potentially degenerate noise coefficients,  results remain pretty scarce. Strong well-posedness
has been proved in~\cite{bib:MP11} for an SPDE~\eqref{eq:SPDE_multi} where $\sigma$  belongs to the class of
H\"older continuous functions with exponent greater than $3/4$ and the drift $h$ is Lipschitz continuous.
As for the SPDE~\eqref{eq:SPDE_multi} with $\sigma$  being a H\"older function with exponent less than
or equal to $3/4$ and zero drift, only weak uniqueness is known for very specific noise coefficients $\sigma$  (such
as, for example,  $\sigma(u) = u^\gamma$
 with $\gamma\geq 1/2$, or $\sigma(u) = \sqrt{u(1-u)}$);
 see~\cite{Mytnik98}, \cite{bib:Sh88}. In \cite{bib:Sh88}, weak uniqueness was established for the noise coefficient  $\sigma(u) = \sqrt{u(1-u)}$ also in the presence of  drifts of the form
$$h(x)=c_1(1-x)-c_2x+c_3x(1-x), \quad x\in[0,1],\;{\rm for}\; c_1, c_2\geq 0, c_3\in \R.$$  
Later,  in~\cite{Athreya00}, \cite{bib:MMR21}, \cite{Barnes24}, for  $\sigma(u) = \sqrt{u(1-u)}$,  the weak uniqueness results were extended to broader classes of drift terms.  
In~\cite{Mytnik98}, \cite{bib:Sh88},  \cite{Athreya00}, \cite{Barnes24} the duality technique was used to prove uniqueness, whereas in~\cite{bib:MMR21}  the Dawson-Girsanov theorem was applied. 

\subsection{Super-Brownian and duality}
Now we turn to the super-Brownian motion (SBM), corresponding to the setting $\sigma(X)=\sqrt{X}$, that is SPDE~\eqref{eq:SPDE1}. As noted earlier, the uniqueness for drifts of the form $h(X)=cX$ has been known for a long time, and follows from the duality method. It is also worth noting, that
when SBM  is subject to nonzero space-dependent immigration, meaning $h$ is continuous and depends only on the  spatial parameter $x$ and not on  $X$, pathwise {\it non-uniqueness} was shown in \cite{Chen15}.

\paragraph{} 
 Proving weak well-posedness for~\eqref{eq:SPDE1} requires establishing both weak existence and weak uniqueness. 
Weak existence, under the assumption that  $h$ is continuous, 
 follows from \cite{Shiga94}.  However, no uniqueness results were provided  in~\cite{Shiga94} for non-Lipschitz coefficients.\\
For   SBM  and other related measure-valued branching processes the duality approach is a widely used method for establishing  weak uniqueness. We will also employ this approach in our analysis. 
 Our approach is essentially based on an extension of the well-known duality relation between  SBM and the corresponding log-Laplace equation. More precisely, let $X$ denote the solution to \eqref{eq:SPDE1} with $h = 0$ and 
 let $V(\mu)$ be the solution to 
 \begin{equation}\label{eq:PDE1}
 	\begin{cases}
 		&d_t V_t(\mu) = \frac{1}{2}\Delta V_t(\mu) - \frac{1}{2}V_t(\mu)^2, \quad (t,x ) \in (0,\infty) \times \R \\
 		& V_t(\mu) \Rightarrow \mu, \text{ as } t \rightarrow 0,
 	\end{cases}
 \end{equation}
 for a finite measure $\mu$, then it is well-known (see for instance \cite[Lemma 1.2]{Konno88}) that 
 \begin{equation*}
 	\E[\exp(-\la X_t, \phi \ra)] = 	\E[\exp(-\la X_0, V_t(\phi)\ra)], \quad \forall t \geq 0,
 \end{equation*}
 for all non-negative, continuous and compactly supported  functions $\phi$. Furthermore, if $h\equiv a >0$, the duality relation becomes:  
 \begin{equation}\label{eq:duality_sbm_imm}
 	\E[\exp(-\la X_t, \mu \ra)] = 	\E\Big[\exp\Big(-\la X_0, V_t(\mu)\ra- a\int_0^t \la V_s(\mu),1\ra ds \Big)\Big], \quad \forall t \geq 0,
 \end{equation}
 which can be verified by using It\^o's formula. In this article, 
 we extend this approach for more complex 
drift terms $h$ by introducing a space-time jump 
noise into the log-Laplace  equation~\eqref{eq:PDE1}.
The resulting duality relation is presented in  Proposition~\ref{prop:duality}. In fact, Proposition~\ref{prop:duality} is  also used  in the proof of  Theorem~\ref{thm:scp}.

We outline the construction of the dual process in 
 detail in Section~\ref{sec:mainres}.  
We would like to highlight the paper~\cite{Barnes24} which provided a number of very important  insights that helped us in proving our main results. 
In~\cite{Barnes24}  the well-posedness results of \cite{Athreya00} were extended for~\eqref{eq:SPDE_multi} with $\sigma(u)=\sqrt{u(1-u)}$  for a  large class of irregular drift terms, including  some  that are H\"older continuous and discontinuous  at zero.  The duality method, based on moment duality, was again used in~\cite{Barnes24}, and the 
dual process was  a branching–coalescing system of Brownian motions.
 A key element used  in the analysis in~\cite{Barnes24} 
was the property of coming down from infinity of the dual branching–coalescing Brownian particle system. This property 
is discussed in detail in \cite{Barnes24_1} in the case of coalescing Brownian motions without branching. This phenomenon 
was crucial for establishing well-posedness in~\cite{Barnes24}  for certain  H\"older continuous and {\it discontinuous} at zero  drifts.   
Analogously, in this article, we rely on the property of coming down from infinity, applied not to a particle system, but to solutions of the log-Laplace equation, whereas we use Laplace transform duality instead of moment duality.
To successfully use this coming down from infinity property,
we employ the theory of initial trace solutions to the log-Laplace equation (see \cite{Brezis95, LeGall96, Marcus99}) which in turn allows us to extend the class of admissible drift functions $h$, particularly to those exhibiting discontinuities at zero.
More precisely, this framework of initial trace solutions enables us to introduce jumps of infinite height in the dual process. In Section~\ref{sec:mainres}, we provide some intuition on how this leads to a duality relation even when $h$ has  discontinuities at zero. 
Let us emphasize that the theory of initial trace solutions is  conceptually closely related to the coming down from infinity property for particle systems used in \cite{Barnes24}. \\
Compared to \cite{Barnes24}, we are able to control exponential moments of the number of jumps in the dual process over sufficiently small time intervals. This control allows us to extend the class of admissible drift terms $h$ to include cases with $b_0 >0$ and $b_1\leq0$, 
whose counterparts were not covered in~\cite{Barnes24}. For instance, our results cover the drift  $h(x)= b_0\mathbbm{1}_{x = 0}$, that allows for  upward drift  only at points where solution equals to zero, an intriguing case for which  no well-posedness results were previously known and which is not treated  in~\cite{Barnes24}.  \\
We believe that establishing  well-posedness for SBMs with irregular drifts, especially those discontinuous at zero, addresses an important open question in the theory of super-Brownian motions. Furthermore, we expect  these results can be extended to more general classes of superprocesses, such as $\alpha$-stable superprocesses with branching mechanisms that are not necessarily quadratic.
It would also be interesting to construct a branching particle system whose scaling limit yields a solution to equation~\eqref{eq:SPDE1} with a drift discontinuous at zero. The discontinuity of the drift at zero in~\eqref{eq:SPDE1} suggests that, in the pre-limiting particle system, there may be different subcritical or supercritical branching mechanisms  or mass immigrations, depending on whether the spatial regions have non-negligible or negligible particle "density". These questions are left for future work. 

\subsection*{Organization of the Article}
First, in Section~\ref{sec:prelim}, we introduce several function spaces and notation used throughout this article. Next, in Section~\ref{sec:mainres}, we outline the proof of Theorem~\ref{th:uniqueness}.
  \LM{Then, in Section~\ref{sec:examples}}, we discuss several examples of admissible drift terms.  
In Section~\ref{sec:parabolic_eq}, we provide a brief summary of results on certain non-linear parabolic equations, focusing in particular on the theory of initial traces and very singular solutions.
Section~\ref{sec:duality} is devoted to the construction of the dual process, its properties, and the duality relation.   
\LM{In Section~\ref{sec:artificial_branching},} we introduce a branching process that \LM{in particular} allows us to control the moments of the number of jumps in the dual process. 
We also  approximate the  dual process by truncating the size of the jumps and prove the convergence of this approximation in Section~\ref{sec:convergence_dual}. \LM{In Section~\ref{sec:appr_dual}, we derive  an approximate duality, and finally, Proposition~\ref{prop:duality}~(the duality relation) is proved.}
In Sections~\ref{sec:weak_exist} and \ref{sec:weak_unique}, \LM{we use Proposition~\ref{prop:duality} and other results} to prove the weak existence and uniqueness of solutions. 
Section~\ref{sec:scp} is dedicated to the proof of Theorem~\ref{thm:scp}, where we \LM{again} use the duality relation. Finally, in Section~\ref{sec:survival} we prove Theorem~\ref{thm:survival}.


\subsection*{Convention for Constants}
In this article, constants are denoted by $C$ and $c$, with their precise values considered irrelevant. We reserve the freedom to adjust these values without altering the notation for readability. 
When necessary, dependence on additional parameters is expressed as $C = C(\lambda)$.

\subsection*{Acknowledgments}
We want to thank Zhenyao Sun and Ed Perkins for very helpful conversations. 
The work of the authors was supported in part by the ISF grant
No. 1985/22. 

\section{Notation}\label{sec:prelim}
In this section, we introduce some notation.

We use the following notation for common function spaces.
\begin{enumerate}
	\item $\mathcal{C}(\R)$ denotes the space of real-valued, continuous functions on $\R$ and $\mathcal{C}^k(\R)$, $k \in N$ denotes the space of $k$-times continuously differentiable, real-valued functions on $\R$.
	\item $\mathcal{C}_c(\R)$ denotes the space of compactly supported functions in $\mathcal{C}(\R)$. 
	\item $	\mathcal{C}_{tem} = \{ f \in \mathcal{C}(\R) \colon |f|_{(\lambda)}<\infty \quad \forall \lambda <0\}$, where $	|f|_{(\lambda)} = \sup_{x\in \R}|\ee^{\lambda|x|}f(x)|$ for $f\in \mathcal{C}(\R)$.
	\item $\mathcal{P}(\mathcal{C}_{tem}^+)$ denotes the set of all probability measures
	on $\mathcal{C}_{tem}^+$.
	\item $	\mathcal{P}_p(\mathcal{C}_{tem}^+) = \{ \mu \in \mathcal{P}(\mathcal{C}_{tem}^+) \colon \sup_{x\in \R}\ee^{\lambda |x|}\int_{\mathcal{C}_{tem}^+}|f(x)|^p\mu(df) < \infty, \quad \forall \lambda <0\}.$
	\item $\mathcal{C}(\R_+,\mathcal{C}_{tem}^+)$ is the space of functions $f\colon \R_+ \rightarrow \mathcal{C}_{tem}^+$, which are continuous with respect to the topology induced by $\mathcal{C}_{tem}$.
	\item $L^p(A,m)$ for $p \in [1,\infty]$ denotes the usual $L^p$-space on the measure space $(A, \mathcal{A},m)$. If there is no ambiguity we write $L^p$ instead of $L^p(A,m)$ for $p \in [1,\infty]$.
	\item $L^p_{loc}((0,\infty) \times \R)$ for $p \in [1, \infty)$  denotes the set of measurable functions $f\colon [0,\infty)\times \R \rightarrow \R$ such that $\int_0^T\int_\R |f(t,x)|^p dx dt<\infty $ for all $T >0$.
	\item $\mathcal{M}_F^+$ denotes the set of finite, non-negative measures on $\R$.
	\item $\mathcal{M}_{reg}^+$ denotes the set of non-negative outer regular measures on $\R$.
	\item $\mathcal{M}_{Rad}^+(B)$ denotes the set of non-negative Radon measures on an open set $B \subseteq \R$.
	\item $\mathcal{B}(\R)$ denotes the Borel $\sigma$-algebra.
\end{enumerate}
If $A$ is a sub-space $\mathcal{C}(\R)$, we denote by $A^+$ the subset of non-negative functions in $A$. Furthermore, we denote by $A_c$ the subspace of compactly supported functions in $A$. 
We also introduce the following notation
\begin{equation*}
	\la \mu,f\ra = \int_\R f(x)\mu(dx), 
\end{equation*}
 for a measure $\mu \in \mathcal{M}_F^+$ and a function $f\in \mathcal{C}(\R)$. Furthermore, if $g\in L^1(\R)$ we also write 
 \begin{equation*}
 	\la f,g\ra = \int_\R f(x)g(x)dx.
 \end{equation*}
 Furthermore,  for a measurable subset $D\subset \R$ we denote by $\mathrm{Leb}(D)$ the Lebesgue measure of $D$.\\
 Finally, if $I$ is a set or a tuple, $|I| = \#I$ denotes the number of elements in the set or the length of the tuple, respectively. More precisely, if we have a tuple $a = (a_1, \ldots a_n)$, then $|a| = n$.

\section{Proof Outline}\label{sec:mainres}
	In this section, we outline the proof of Theorem~\ref{th:uniqueness} and, in particular, motivate the duality relation, which is the main tool in our analysis.
	Before we outline our arguments in more detail, we give a definition of weak solution to \eqref{eq:SPDE1}. 
	\begin{dfn}[Weak Solution]\label{def:weak_sol}
		We say that \eqref{eq:SPDE1} has a weak solution $(X_t)_{t\geq 0}$, if there exists a filtered probability space $(\Omega, \mathcal{F}, (\mathcal{F}_{t})_{t\geq 0}, \P)$ such that $(X_t)_{t\geq 0}$ is an adapted $\mathcal{C}_{tem}^+$ - valued 
		continuous process with $X_0 = \mu$ a.s. and such that $\dot{W}$ is a space-time white noise and for every $(t,x) \in (0,\infty)\times \R$ almost surely
		\begin{equation*}
		\begin{split}
			&X_t(x) = \int_\R p_t(x-y) X_0(y)dy + \int_\R \int_0^t p_{t-s}(x-y)h(X_s(y))dsdy \\
			& \hspace{2cm}+ \int_\R \int_0^t p_{t-s}(x-y) \sqrt{X_s(y)}\dot{W}(dsdy).
		\end{split}
		\end{equation*}
	Hereby, 
	\begin{equation}\label{eq:heat_dens}
		p_t(x) := \frac{1}{\sqrt{2\pi t}}\ee^{-x^2/(2t)}.
	\end{equation}
	\end{dfn}\par
	\noindent The proofs of weak uniqueness and of particular cases of weak existence in Theorem~\ref{th:uniqueness} are based on a duality approach.  \\ 
	We first outline this approach for a simplified setting and then step by step add additional modifications to cover more general drift terms.
	Let us first consider the case where 
	\begin{equation}\label{eq:simplified}
			h_{\infty}(x) = \la\nu^1,1\ra, \quad \forall x\in \R_+  \text{ and } \nu^2 \equiv 0,
	\end{equation}
	which means that 
		\begin{equation}\label{eq:h_simpler_1}
		h(x) = \int_0^\infty (1-\ee^{-\lambda x})\nu^1(d\lambda).
	\end{equation}
	If we  apply It\^o's formula to $\exp(-\langle \phi, X_t \rangle)$, for $\phi \in \mathcal{C}^2(\R)$ with $\phi \geq 0$. We obtain that
	\begin{equation}\label{eq:martingale_X}
		\exp(-\langle \phi, X_t \rangle) - \int_0^t \exp(-\la\phi, X_s\ra)\Big(-\la \frac{1}{2}\Delta \phi, X_s\ra + \frac{1}{2}\la \phi^2, X_s\ra -\la \phi, h(X_s) \ra\Big)ds, \quad t \geq 0,
	\end{equation}
	is an $\F_t^X$-martingale. 
	Suppose we can show the existence of a measure-valued process $Y$ such that 
	\begin{equation}\label{eq:martingale_Y}
		\exp(-\langle Y_t, \psi \rangle) - \int_0^t \exp(-\la Y_{s-}, \psi\ra)\Big(-\la \frac{1}{2}\Delta Y_{s-}, \psi\ra + \frac{1}{2}\la Y_{s-}^2, \psi\ra -\la Y_{s-}, h(\psi) \ra\Big)ds, \quad t \geq 0,
	\end{equation}
	is an $\mathcal{F}^Y_t$-martingale for $\psi \in \mathcal{C}_{tem}^+$. Then, if $X,Y$ are independent, under certain technical conditions it follows that 
	\begin{equation*}
		\E[\exp(-\la X_t, Y_0\ra)] = 	\E[\exp(-\la X_0, Y_t\ra)], \quad \forall t \geq 0,
	\end{equation*}
	by \cite[Theorem 4.4.11]{Ethier05}.
	This would allow us to infer that the one-dimensional laws are unique and thus by standard arguments the finite dimensional laws of $X$ are unique as well.
	 Consider  the SPDE
	\begin{equation}\label{eq:SPDE2}
		d_tY_t = \frac{1}{2}\Delta Y_t -\frac{1}{2} Y_t^2 +  \int_0^\infty \int_0^\infty \lambda\mathbbm{1}_{r \leq Y_{t-}(x)}\mathcal{N}(dr,d\lambda,dx, dt), 
	\end{equation}
	where $\mathcal{N}$ is a Poisson random measure with intensity 
	\begin{equation*}
		drdx \nu^1(d\lambda) dt.
	\end{equation*}
	A formal application of It\^o's formula shows that a solution to \eqref{eq:SPDE2} solves \eqref{eq:martingale_Y} 
	and thus the solution to \eqref{eq:SPDE2}  is our candidate for the dual process. More precisely, we can derive the following formal identity, which demonstrates how the jump-type noise in \eqref{eq:SPDE2} allows us to derive the duality relation for drifts of the form \eqref{eq:h_simpler_1}. Denote by $(\bar{T_k}, Z_k, U_k)_{k \geq 1}$ the atoms of 
    \begin{equation*}
        \int_0^\infty\mathbbm{1}_{r \leq Y_{t-}(x)}\mathcal{N}(dr,dx,d\lambda, dt).
    \end{equation*}
    Then we formally get
    \begin{equation}\label{eq:formal_1_ito}
    \begin{split}
        &\E\Big[\sum_{\bar{T}_k\leq t}  \exp(-\langle Y_{\bar{T}_k},\psi \rangle) - \exp(-\langle Y_{\bar{T}_k-},\psi \rangle) \Big] = -\E\Big[\sum_{\bar{T}_k\leq t} \exp(-\langle Y_{\bar{T}_k-},\psi \rangle)(1-\exp(-\langle Z_k \delta_{U_k}, \psi \rangle ))  \Big]\\
        &=- \E\Big[\int_0^t\exp(-\langle Y_{s-},\psi \rangle) \int_\R\int_0^\infty \left( 1-\exp(-\lambda \psi(x))\right)  \nu^1(d\lambda) Y_{s-}(x) dx  ds \Big] \\
        & = -\E\Big[\int_0^t\exp(-\langle Y_{s-},\psi \rangle) \langle Y_{s-}, h(\psi) \rangle ds \Big],
    \end{split}
    \end{equation}
    which allows us to show that \eqref{eq:martingale_Y} is a martingale. 
    \\
	Now, let us consider the case where the $\nu^1,\nu^2 \neq 0$ and $h_\infty \equiv \la\nu^1+\nu^2, 1\ra$. That is, we now consider the case where 
	\begin{equation}\label{eq:h_simplified_1}
		h(x) = \int_0^\infty(1-\ee^{-\lambda x}) \nu^1(d\lambda)+  \int_0^\infty(1+\ee^{-\lambda x})\nu^2(d\lambda).
	\end{equation}
	We will see that the duality function changes in this setting.  In order to define the modified duality function needed for this setting 
	we also introduce an independent sequence of i.i.d. random variables $(M_{k})_{k \geq 1}$ taking values in $\{1,2\}$  such that 
	\begin{equation}\label{eq:def_M}
		\P(M_{k} = i) = \frac{\nu^i((0,\infty))}{\nu^1((0,\infty))+ \nu^2((0,\infty))}, \quad i \in \{1,2\},
	\end{equation}
	We furthermore define
	\begin{equation}\label{eq:def_J_1}
		J_t = \#\{ k \in \N \colon \bar{T}_{k} \leq t , M_{k} = 2\},
	\end{equation}
	where $\bar{T}_{k}$ denotes the $k$-th jump of the process $(\la Y_{t},1\ra)_{t \geq 0}$.
	Hereby, $Y$ is constructed as before, but with $\nu^1$ replaced by $\nu = \nu^1+\nu^2$.
	Then, we will see in Section~\ref{sec:duality} that  the duality relation is given by
	\begin{equation}\label{eq:approximate_sign_change}
		\E[\exp(-\la X_t, Y_0\ra)] = 	\E[(-1)^{J_t}\exp(-\la X_0, Y_t\ra)], \quad \forall t \geq 0.
	\end{equation}
The role of the term $(-1)^{J_t}$
can be seen as follows.
For simplicity, assume $\nu^1 \equiv 0$. 
Then, similarly to \eqref{eq:formal_1_ito} we formally obtain
    \begin{equation*}
    \begin{split}
        &\E\Big[\sum_{\bar{T}_k\leq t}  (-1)^{J_{\bar{T}_k}}\exp(-\langle Y_{\bar{T}_k},\psi \rangle) - (-1)^{J_{\bar{T}_k-}}\exp(-\langle Y_{\bar{T}_k-},\psi \rangle) \Big] \\
        &= -\E\Big[\sum_{\bar{T}_k\leq t} (-1)^{J_{\bar{T}_k-}}\exp(-\langle Y_{\bar{T}_k-},\psi \rangle)(1+\exp(-\langle Z_k \delta_{U_k}, \psi \rangle ))  \Big]\\
        &=- \E\Big[\int_0^t(-1)^{J_{s-}}\exp(-\langle Y_{s-},\psi \rangle) \int_\R\int_0^\infty \left(1+\exp(-\lambda \psi(x))\right) \nu^2(d\lambda) Y_{s-}(x) dx ds \Big] \\
        & = -\E\Big[\int_0^t(-1)^{J_{s-}}\exp(-\langle Y_{s-},\psi \rangle) \langle Y_{s-}, h(\psi) \rangle ds \Big],
    \end{split}
    \end{equation*}
    which is the mechanism that allows us to derive the duality relation \eqref{eq:approximate_sign_change}.

\begin{rem}
Because of the factor $(-1)^{J_t}$, the random variable appearing on the right-hand side of the duality relation~\eqref{eq:approximate_sign_change}  need not be nonnegative. Nevertheless, once the duality relation is established, its expectation is nonnegative, since it equals the nonnegative quantity on the left-hand side.
A simple Poisson example may help illustrate this phenomenon. 
If $(N_t)_{t\geq 0}$, with $N_0 = 0$, is a Poisson process with rate $\lambda >0$, then 
        \begin{equation*}
            	\E\bigl[(-1)^{N_t}\bigr]
		= \sum_{n=0}^{\infty}(-1)^n\mathbb{P}(N_t=n) 
		= \sum_{n=0}^{\infty}(-1)^n e^{-\lambda t}\frac{(\lambda t)^n}{n!} 
		= e^{-\lambda t}\sum_{n=0}^{\infty}\frac{(-\lambda t)^n}{n!}  = e^{-2\lambda t} >0.
        \end{equation*}
This calculation does not apply directly to the present setting, because $J_t$  is not, in general, a Poisson random variable with a
deterministic parameter and is coupled to $Y_t$. It nevertheless 
illustrates that the presence of an alternating sign does not, by itself, prevent the expectation from being nonnegative.
\end{rem}

	Now, we turn to the discontinuous part $h_\infty$.
	If we formally consider $\nu^1  = d_1\delta_\infty$ and $\nu^2 = d_2\delta_\infty$ with $d_1,d_2 \geq 0$ and plug it into 
	\begin{equation}\label{eq:h_simplified}
		h(x) =  \int_0^\infty (1-\ee^{-\lambda x})\nu^1(d\lambda) +\int_0^\infty (1+\ee^{-\lambda x})\nu^2(d\lambda),
	\end{equation}
	we arrive at the function 
	\begin{equation*}
		h_\infty(x) =d_1\mathbbm{1}_{x>0} + 2d_2\mathbbm{1}_{x = 0} + d_2\mathbbm{1}_{x>0} =  \lim_{n \rightarrow \infty}d_1(1-\ee^{-nx}) +\lim_{n \rightarrow \infty}d_2(1+\ee^{-nx}).
	\end{equation*}
	Thus, for each fixed $n$, if we set 
	\begin{equation*}
		\nu_n^1 = \nu^1_{|(0,n]}+d_1\delta_{n} , \quad \nu_n^2 =  \nu^2_{|(0,n]} + d_2\delta_{n},
	\end{equation*}
	 we are again in the setting of \eqref{eq:h_simplified_1} and \eqref{eq:h_simplified} with
	 $\nu^i$ replaced by $\nu^i_n$, $i = 1,2$. Hereby, $\nu^i_{|(0,n])}$, $i = 1,2$ denotes the restriction of $\nu^i$ to the interval $(0,n]$. 
	 Inspired by \eqref{eq:approximate_sign_change}, this will allow us to establish approximate duality. That is, if
	\begin{equation}\label{eq:apprximate_duality_1}
			\E[\exp(-\la X_t, Y_0\ra)] = \lim_ {n \rightarrow \infty}	\E[(-1)^{J_t^n}\exp(-\la X_0, Y_t^n\ra)], \quad \forall t \geq 0,
	\end{equation}
	for each $n\geq 1$, $Y^n$ solves \eqref{eq:SPDE2} with $\nu^1$ replaced by $\nu^1_n+\nu^2_n$ and $J^n$ is defined accordingly as in \eqref{eq:def_J_1}.
	 Then, we will be able to 
	construct  processes $Y, J$, that are limits of $Y^n$ and $J^n$, respectively, such that \eqref{eq:apprximate_duality_1} holds without limit. Loosely speaking, such a $Y$ is a solution to \eqref{eq:SPDE2} with $\nu^1$ formally replaced by $\nu^1+\nu^2 +d_1\delta_\infty+d_2\delta_\infty$ and $J$ is constructed accordingly (see Section~\ref{sec:duality}).\\
With the duality function as described above, we cannot cover the full range of parameters $b_0\geq 0$ and $b_1 \in \R$ that are admissible for $h_\infty$. In fact, with \eqref{eq:approximate_sign_change} we are only able to prove weak well-posedness for drifts of the form 
\begin{equation*}
\begin{split}
	h_1(x) &= \int_0^\infty \ee^{-\lambda x}(\nu^2-\nu^1)(d\lambda),\\
	 h_\infty(x) &= 2d_2\mathbbm{1}_{x = 0} + (d_1+d_2)\mathbbm{1}_{x>0} + \la \nu^1+\nu^2,1\ra.
\end{split}
\end{equation*}
Since $d_1,d_2 \geq 0$, for the case of $\nu^1,\nu^2 \equiv0$ this, for example, does not cover the regime of $b_0 > 0$ and $b_1\leq 0$.
Introducing an additional parameter $a\in\R$ and setting 
\begin{equation}\label{eq:hinftyd}
	h_\infty(x) =  2d_2\mathbbm{1}_{x = 0} + (d_1+d_2)\mathbbm{1}_{x>0} +  \la \nu^1+\nu^2,1\ra+a,
\end{equation}  
allows us to cover the complete regime of parameters $b_0 \geq \langle \nu^1-\nu^2,1 \rangle $ and $b_1 \in \R$ for $h_ \infty$ given by \eqref{eq:hinfty}. 
More precisely, we set 
\begin{equation}\label{eq:def_d}
\begin{split}
		d_1 = (b_1-b_0)^+, \quad d_2 = (b_0-b_1)^+ , \quad a = b_0 - \langle \nu^1+\nu^2,1\rangle -2(b_0-b_1)^+.
	\end{split}
\end{equation}
It is easily verified that by plugging these definitions into \eqref{eq:hinftyd} we recover \eqref{eq:hinfty}.
\begin{rem}
	There is not necessarily a unique choice of $\nu^1,\nu^2$ and $d_1$, $d_2$ and $a$ that allows to represent a given function $h$.
\end{rem}
In order to handle this class of drift functions, we need to modify the duality function again. In this case we will see in Section~\ref{sec:duality} that the approximate duality relation is given by 
\begin{equation}\label{eq:approx_duality_final}
\begin{split}
		&\E\Big[\exp\Big(-\la X_t, Y_0\ra\Big)\Big] \\
		&\quad = \lim_ {n \rightarrow \infty}	\E\Big[(-1)^{J_t^n}\exp\Big(-\la X_0, Y_t^n\ra\Big) \exp\Big(-a\int_0^t \la Y_s^n,1\ra ds\Big)\Big], \quad \forall t \geq 0,
\end{split}
\end{equation}
where $Y^n$ and $J^n$ are defined as in \eqref{eq:apprximate_duality_1} with $d_1,d_2$ as in \eqref{eq:def_d}.
The relation \eqref{eq:approx_duality_final} can  again be easily motivated by a formal application of Ito's formula. Hereby, we want to point out that $a$ may be negative. This means, one of our main tasks is to show that the expectation on right hand side of \eqref{eq:approx_duality_final} is well defined.
As before, we are able to construct the limits of $Y^n$ and $J^n$ and thus  do away with the limit in \eqref{eq:approx_duality_final}.
These considerations lead us to the following proposition, which we will prove in Section~\ref{sec:duality}.
\begin{proposition}\label{prop:duality}
	Assume that $(X_t)_{t \geq 0}$ is a solution to \eqref{eq:SPDE1} with $X_0 \sim\eta \in \mathcal{P}_p(\mathcal{C}_{tem}^+)$, $p \geq 2$ and such that \eqref{eq:cond_x_0} holds. Let $Y_0 = \mu $ for $\mu \in \mathcal{M}_F^+$. Then, there exists $T_0>0$ such that there exist processes $(J_t)_{t\geq0}$ and $(Y_t)_{t\geq0}$, taking values in $\N\cup \{ 0\}$ and $\mathcal{M}_{reg}^+$, respectively, that are independent of $X$ and such that
	\begin{equation}\label{eq:duality}
				\E\Big[\exp\Big(-\la X_t, Y_0\ra\Big)\Big] = 	\E\Big[(-1)^{J_t}\exp\Big(-\la X_0, Y_t\ra\Big) \exp\Big(-a\int_0^t \la Y_s,1\ra ds\Big)\Big], 
	\end{equation}
for all $t \in [0,T_0]$. Hereby, $a$ is given by \eqref{eq:def_d} and $T_0$ depends on $a$ and $\la 1,\mu\ra$.
\end{proposition}
\begin{rem}\label{rem:T0}
	Note that in the case $a >0$, we can prove \eqref{eq:duality} for all $t>0$. 
\end{rem}
Note that $(Y,J)$ are the limits of $(J^n,Y^n)$ in \eqref{eq:approx_duality_final} that we discussed above with $d_1,d_2$ as in~\eqref{eq:def_d}.
Finally,  when $\nu^2 \equiv0$ and $ h_\infty \equiv \la\nu^1,1\ra$ we can use Proposition~\ref{prop:duality}
 \LM{ to prove Theorem~\ref{thm:scp}.}

\section{Examples}\label{sec:examples}
As already mentioned above, Theorem~\ref{th:uniqueness} establishes weak existence and uniqueness of solutions to \eqref{eq:SPDE1}
for discontinuous and non-Lipschitz drifts. Below we give several examples. We are particularly interested
in examples that allow us to investigate the behavior of solutions to \eqref{eq:SPDE1} at zero and compare our findings to the corresponding one-dimensional stochastic differential equations.

\subsection{Continuous Non-Lipschitz Drift}
Set $\nu^2 \equiv 0$, $\nu^1(d\lambda) = \lambda^{-1-\alpha}\mathbbm{1}_{\lambda \geq 1}d\lambda$ with $\alpha \in (0,1)$ and 
\begin{equation*}
	b_0 = b_1 = \langle \nu^1,1\rangle. 
\end{equation*}
Then it follows that 
\begin{equation*}
	h(x) = \int_1^\infty (1-\ee^{-\lambda x}) \lambda^{-1-\alpha}d\lambda = |x|^\alpha \int_{x}^\infty(1-\ee^{-\lambda}) \lambda^{-1-\alpha}d\lambda. 
\end{equation*}
Clearly,
\begin{equation*}
C_1(\alpha)|x|^\alpha\leq	|h(0)-h(x)| = h(x)  \leq C_2(\alpha)|x|^\alpha, \quad x\in [0,1],
\end{equation*}
is not Lipschitz continuous at $0$, but only $\alpha$-Hölder continuous.
\subsection{Completely Monotone Functions}
Let us now consider the case where $b_0 = b_1 = c \in\R$ and $\nu^1,\nu^2$ are finite measures on $(0,\infty)$ so that $c+\langle \nu^2-\nu^1,1\rangle \geq 0$ . That means we have 
\begin{equation}\label{eq:h_mon}
	h(x) = c+ \int_0^\infty \ee^{-\lambda x} (\nu^2-\nu^1)(d\lambda).
\end{equation}
By Bernstein's theorem on completely monotone functions (see \cite{Bernstein29})  all functions in the span of bounded completely monotone functions that are nonnegative at zero take the form \eqref{eq:h_mon}.
Recall that a completely monotone function $f$ on $[0, \infty)$ is continuous everywhere and in $\mathcal{C}^\infty((0,\infty))$ and satisfies
\begin{equation*}
	 (-1)^nf^{(n)}(x) \geq 0,  \quad x>0,
\end{equation*}
where $f^{(n)}$ denotes the $n$-th derivative of $f$. 
\subsection{Discontinuous Drift }
In this section we discuss the discontinuous part of the drift. That is, we set  
\begin{equation*}
	 h(x) = h_{b_0,b_1}(x)= b_1\mathbbm{1}_{x > 0} + b_0\mathbbm{1}_{x =0}, \quad x \geq 0,
\end{equation*}
for $b_0 \geq 0$ and $b_1 \in \R$.
This example allows us to examine some differences that occur in stochastic partial differential equations of the form \eqref{eq:SPDE1} compared to their scalar counterparts of the form 
\begin{equation}\label{eq:SDE_sbm}
	dx_t =  h_{b_0,b_1}(x_t)dt + \sqrt{x_t}dB_t, \quad x_0 = y,  \quad t \geq 0, 
\end{equation}
where $(B_t)_{t\geq 0}$ is a standard Brownian motion. 
It turns out that if we choose $b_0 = 0$ and $b_1 = 1$, weak uniqueness fails for \eqref{eq:SDE_sbm} when $x_0 =0$. However, \eqref{eq:SPDE1} has a unique weak solution for $h \equiv h_{0,1}$.  
Furthermore, if $b_0 = c \in (0,1)$ and $b_1 = 1$, then the laws of the solutions to \eqref{eq:SDE_sbm} are identical for each choice of $c$. The solutions to \eqref{eq:SPDE1} with $h\equiv h_{c,1}$, on the other hand, have distinct laws for distinct values of $c$.   
Finally, if we consider the drift $h_{1,0}$, then \eqref{eq:SDE_sbm} has no solution, while \eqref{eq:SPDE1} is weakly well-posed.
These observations are the content of the following two lemmas.
\begin{lemma}
\begin{enumerate}[label=(\roman*)]
	\item Weak uniqueness does not hold for 
		\begin{equation}\label{eq:sde_non_unique}
					dx_t =  \mathbbm{1}_{x_t>0}dt + \sqrt{x_t}dB_t,  \quad x_0 =0, \quad t \geq 0.
	\end{equation}
	\item  The solutions to 
	\begin{equation}\label{eq:sde_same_dist}
		dx_t =  h_{c,1}(x_t)dt + \sqrt{x_t}dB_t,  \quad t \geq 0, 
	\end{equation}
	have the same law for each $c \in (0,1)$. 
	\item   There exists no solution to the equation
	\begin{equation}\label{eq:sde_nonexist}
		dx_t = \mathbbm{1}_{x_t =0}dt + \sqrt{|x_t|}dB_t, \quad x_0 = 0, \quad t \geq 0.
	\end{equation} 
\end{enumerate}
\end{lemma} 
\begin{proof}
	First, we prove $(i)$. Clearly, $x_t = 0$ for all $t \geq 0$ is a solution to \eqref{eq:sde_non_unique}. Denote by $(y_t)_{t\geq0}$ the unique strong solution to \eqref{eq:SDE_sbm} with $b_0 = 1$, $b_1 =1$ and $y_0 = 0$. We now show that  $(y_t)_{t\geq0}$
 also solves \eqref{eq:sde_non_unique}. To this end it is enough to prove that 
	\begin{equation*}
		\int_0^t \mathbbm{1}_{y_s = 0}ds =0.
	\end{equation*}
	However, this can be proven along the same lines as the proof of \cite[Lemma 1.3]{Barnes24}. That is, by \cite[Theorem VI.1.7]{Revuz99} we get that 
	\begin{equation*}
		L_t^0(y) = 2\int_0^t \mathbbm{1}_{y_s =0}ds,
	\end{equation*}
	where $L_t^0(y)$ is the local time of $y$ at $0$ that arises from Tanaka's formula as defined in \cite[Theorem VI.1.2]{Revuz99} .
	By the same arguments as in  the proof of  \cite[Proposition XI.1.5]{Revuz99} we get that $L_t^0(y) \equiv 0$.
	 This finishes the proof of $(i)$. 
	 
	Next we prove $(ii)$. In order to prove the claim, we show that any solution to \eqref{eq:sde_same_dist} with $c\in (0,1)$ solves 
	\begin{equation}\label{eq:sde_1}
			dy_t =  1dt + \sqrt{y_t}dB_t,  \quad t \geq 0.
	\end{equation}
	Note that \eqref{eq:sde_1} has pathwise unique solutions and thus also weakly unique solutions. In order to show that any solution $x$ to \eqref{eq:sde_same_dist} also solves \eqref{eq:sde_1} it is enough to prove 
	\begin{equation*}
		\int_0^t \mathbbm{1}_{x_s =0}ds = 0.
	\end{equation*} 
	This can once again be achieved in a similar manner to the above.\\
	Finally, we prove $(iii)$. Assume that there exists a weak solution $x$ to \eqref{eq:sde_nonexist} and note that this solution has to be nonnegative by  a slight modification of \cite[Theorem 1.4]{LeGall84}. Then, similar to above we 
	get that 
	\begin{equation*}
		\int_0^t \mathbbm{1}_{x_s =0}ds =0.
	\end{equation*}
	However, this means that $(x_t)_{t \geq 0}$ is also a solution to
	\begin{equation*}
		dx_t = \sqrt{x_t}dB_t, \quad x_0 = 0.
	\end{equation*}
	However, this equation has a pathwise unique solution by the Yamada-Watanabe Theorem, which is given by $x_t \equiv 0$. This would imply that 
	\begin{equation*}
		0 = 	\int_0^t \mathbbm{1}_{x_s =0}ds  = t,
	\end{equation*}
	which is a contradiction. Thus, there exists no solution to this equation.
\end{proof}
\begin{lemma}
	\item The laws of  solutions to 
	\begin{equation}\label{eq:spde_distinct}
		dX_t =  \frac{1}{2}\Delta X_t + h_{c,1}(X_t)dt + \sqrt{X_t}\dot{W}, \quad X_0 =f\in \mathcal{C}_b^+,  \quad t \geq 0, 
	\end{equation}
	are different  for different values of $c \in (0,1)$. 
\end{lemma}
\begin{proof}
	 The proof is similar to the proof of \cite[Lemma 1.2]{Barnes24}. Fix $c_1,c_2 \in (0,1)$ with $c_1\neq c_2$ and denote by $X^1, X^2$ the solution to \eqref{eq:spde_distinct} with drift $h_{c_1,1}$ and $h_{c_2,1}$, respectively. 
	 We prove the claim by contradiction. To this end assume that $X^1 \overset{d}{=} X^2$. Then it follows along the same lines as in the proof of \cite[Lemma 1.2]{Barnes24}
	 that 
	 \begin{equation*}
	 	\begin{split}
	 		c_1 \E\Big[ \int_\R\int_0^t \phi(x) \mathbbm{1}_{X_s^1(x) =0}dsdx\Big] = 	c_2 \E\Big[ \int_\R\int_0^t \phi(x) \mathbbm{1}_{X_s^2(x) =0}dsdx\Big], 
	 	\end{split}
	 \end{equation*}
for all nonnegative $\phi \in \mathcal{C}_c^\infty$ and all $t \geq 0$. In order to finish the proof it is enough to show that 
\begin{equation*}
	\E\Big[ \int_\R\int_0^t \phi(x) \mathbbm{1}_{X_s^1(x) =0}dsdx\Big]  >0,
\end{equation*}
since we assumed that $X^1$ and $X^2$ have the same law. To this end, we prove 
\begin{equation*}
	\E[ \mathbbm{1}_{X_t^{i}(x) =0}] >0,\quad \forall t >0,x\in \R , i =1,2,
\end{equation*}
by a comparison argument, which follows along similar lines to the proof of \cite[Lemma 1.2]{Barnes24}. Let $X$ be the unique weak solution to 
\begin{equation*}
	 d_tX_t  = \frac{1}{2}\Delta X_t + 1dt + \sqrt{X_t}\dot{W}, \quad x_0 =f,  \quad t \geq 0,
\end{equation*}
which is the density of super-Brownian motion with constant immigration.
Furthermore, recall that $h_{1,1}\equiv 1$ and note that
\begin{equation*}
h_{c_i,1}(x) \leq h_{1,1}(x), \quad \forall x\in \R.
\end{equation*}
Moreover, $X^i$ as well as $X$ can be constructed as weak limits of an SPDE with Lipschitz continuous coefficients. 
Thus, by \cite[Corollary 2.4]{Shiga94} and weak uniqueness we get that $X$ stochastically dominates $X^i$, which implies 
\begin{equation*}
		\E[ \mathbbm{1}_{X_t^{i}(x) =0}]  \geq 	\E[ \mathbbm{1}_{X_t(x) =0}], \quad \forall t >0,x\in \R , i =1,2.
\end{equation*}
Using the duality formula \eqref{eq:duality_sbm_imm} we get for $t>0$ that
\begin{equation*}
\begin{split}
		&\E[ \mathbbm{1}_{X_t(x) =0}] = \lim_{m \rightarrow \infty} \E[ \exp(-m X_t(x))] \\
		&= \lim_{m \rightarrow \infty} \exp\Big(- \la V_t(m\delta_x), X_0 \ra - \int_0^t\la V_s(m\delta_x),1\ra ds\Big).
\end{split}
\end{equation*}
Using Proposition~\ref{prp:initial_trace_prop} we get that
\begin{equation*}
	\begin{split}
		&\lim_{m \rightarrow \infty} \exp\Big(- \la V_t(m\delta_x), X_0 \ra - \int_0^t\la V_s(m\delta_x),1\ra ds\Big) \\
		& = \exp\Big(- \la \mathcal{W}_t(x, 0) , X_0 \ra - \int_0^t\la \mathcal{W}_s(x,0 ),1\ra ds\Big) >0,
	\end{split}
\end{equation*} 
where $\mathcal{W}(x,0)$ is defined as in \eqref{eq:very_sing_def}.
This completes the proof.
\end{proof}

\section{Non-linear Parabolic Equations and Initial Trace}\label{sec:parabolic_eq}
In this section, we discuss the nonlinear parabolic equation
\begin{equation}\label{eq:parabolic_eq}
	d_t V_t = \frac{1}{2}\Delta V_t - \frac{1}{2}V_t^2, \quad t\geq 0
\end{equation}
on $\R$. This partial differential equation plays a central role in our proofs due to its close relation
with the dual process. 
It is possible to derive well posedness results for \eqref{eq:parabolic_eq} for highly irregular initial data.
For a more thorough discussion we refer the reader to \cite{Marcus99} and \cite{LeGall96} . Here we shall only give a brief overview.
More precisely, we discuss some basic properties of solutions to 
\begin{equation}\label{eq:PDE}
	\begin{cases}
		&d_t V_t(\mu) = \frac{1}{2}\Delta V_t(\mu) - \frac{1}{2}V_t(\mu)^2, \quad (t,x ) \in (0,\infty) \times \R \\
		& V_t(\mu) \Rightarrow \mu, \text{ as } t \rightarrow 0,
	\end{cases}
\end{equation}
where $\mu$ is a finite measure as well as solutions to the initial trace problem
\begin{equation}\label{eq:initial_trace}
	\begin{cases}
		&d_t \mathcal{W}_t(A,\nu) = \frac{1}{2}\Delta \mathcal{W}_t(A,\nu) - \frac{1}{2}\mathcal{W}_t(A,\nu)^2, \quad (t,x ) \in (0,\infty) \times \R\\
		& \{ y \in \R \colon \lim_{t \rightarrow 0}\int_{y-r}^{y+r} \mathcal{W}_t(A,\nu)(x)dx = \infty, \forall r >0\} = A\\
		& \lim_{t \rightarrow 0} \int \phi \mathcal{W}_t(A,\nu)(x) dx = \int \phi \nu(dx), \forall \phi \in \mathcal{C}_c(A^c), 
	\end{cases}
\end{equation}
where $A \subset \R$ is a closed set and $\nu$ is a non-negative Radon measure on $A^c$. If $A = \{x\}$ for some $x\in \R$ we sometimes write 
\begin{equation}\label{eq:very_sing_def}
	\mathcal{W}_t(x,\nu) = V_t(\nu +\infty\delta_x), \quad t>0, x\in \R.
\end{equation}
Furthermore, for $t < 0$ we use the convention that $\mathcal{W}_t(A,\nu) = V_t(\mu) = 0$ and $V_0(\mu) = \mu$, for $\mu \in \mathcal{M}_F^+$ and $\nu \in \mathcal{M}_{Rad}^+(A^c)$.
\subsection{Finite Measures}
In \cite{Brezis83} and \cite{Fleischmann88, Mytnik02} the authors prove existence and uniqueness of solutions to \eqref{eq:PDE} as well as several estimates.  We summarized their results in this section.
First, we get existence and uniqueness by \cite[Lemma 2.1]{Mytnik02}, which we restated in the following proposition
\begin{proposition}
	Let $\mu$ be a finite measure. Then there exists a unique solution $V(\mu)$ \eqref{eq:PDE} in  $L^1_{loc}((0,\infty) \times \R)$ satisfying
	\begin{equation}\label{eq:sol_rep_V}
		V_t(\mu)(x) = (S_t\mu)(x) -\int_0^t \frac{1}{2}\Big(S_{t-s}V_s(\mu)^2\Big)(x)ds, \quad t>0, x\in \R,
	\end{equation}
	where $(S_t)_{t\geq 0}$ is the heat semi-group with transition density $(p_t)_{t \geq 0}$ given by \eqref{eq:heat_dens}. 
\end{proposition}
The following lemma is a well known fact about \eqref{eq:PDE} (see for instance \cite[Lemma 2.1]{Mytnik02}).
\begin{lemma}\label{lem:est_pde}
	Let $\mu$ be a finite measure. Then it holds that 
	\begin{equation*}
		V_t(\mu)(x) \leq S_t\mu(x),
	\end{equation*}
	for all $t \geq 0$ and $x\in \R$.
\end{lemma}
The following lemma is \cite[Lemma 2.6]{Mytnik02} and provides estimates that we shall use frequently. 
\begin{lemma}\label{lem:pde_est}
	Let $\mu$ and $\eta$ be two finite measures. Then the following holds:
	\begin{enumerate}
		\item[(i)]
		For all $t \geq 0$ it holds that
		\begin{equation*}
			V_t(\mu+ \eta)(x) \leq V_t(\mu)(x) + V_t (\eta)(x), \quad x\in\R.
		\end{equation*}
		\item[(ii)] If $\mu \leq \eta$, then 
		\begin{equation*}
			V_t(\mu)(x) \leq V_t(\eta)(x),\quad x\in \R,
		\end{equation*}
		for all $t \geq 0$.
	\end{enumerate}
\end{lemma}
Finally, using \cite[Lemma 2.3]{Klenke98} we can infer the following lemma.
\begin{lemma}\label{lem:V_asymp}
	Let $\mu$ be a finite, non-zero measure. Then it follows that 
	\begin{equation*}
		\int_0^\infty \la V_t(\mu),1\ra dt = \infty.
	\end{equation*}
\end{lemma}

\subsection{Initial Trace}
In this section we shall give a brief overview of results about solutions to \eqref{eq:initial_trace} based on \cite{Marcus99}, \cite{LeGall96} and \cite{Barnes24}.
We denote by $\mathcal{T}$ the set of pairs $(A,\nu)$, where $A$ is a closed subset of $\R$ and $\nu$ is a Radon measure supported on $A^c$.
We  introduce the partial order $\preceq$ on $\mathcal{T}$, where 
\begin{equation*}
	(A,\nu) \preceq (\hat{A},\hat{\nu}), \text{ if } A \subset \hat{A} \text{ and } \nu \leq \hat{\nu}.
\end{equation*}
We furthermore denote by $\mathcal{M}_{reg}^+$ the set of outer regular Borel measures on $\R$ and define the map 
\begin{equation*}
	\eta \colon \mathcal{T} \rightarrow \mathcal{M}_{reg}^+ , (A,\nu) \mapsto \eta^{(A,\nu)}, 
\end{equation*}
where 
\begin{equation*}
	\eta^{A,\nu}(B) = \begin{cases}
		\infty, B\cap A \neq \emptyset, \\
		\nu(B), B\cap A = \emptyset,
	\end{cases}
\end{equation*}
for all Borel-measurable sets $B \subset \R$. 
Since $\nu \in \mathcal{M}_{reg}^+$ is not necessarily locally bounded, we introduce the set of singular points $\mathcal{S}_{\eta}$ associated with $\eta$ given by  
\begin{equation*}
	\mathcal{S}_\eta = \{ x \in \R \colon \eta(U) = \infty \text{ for every open neighborhood } U \text{ of } x\},
\end{equation*}
as well as the set of regular points $R_\eta$ given by $R_\eta = \mathcal{S}_\eta^c$.
In order to formulate the results of this section, we furthermore need the notion of $m$-weak convergence.
\begin{dfn}
	A sequence $(\eta_n)_{n\geq 1} \subset \mathcal{M}_{reg}^+$ is said to converge $m$-weakly to $\eta \in \mathcal{M}_{reg}^+$ if the following is satisfied:
	\begin{enumerate}
		\item[(i)]  For all open sets $U \subset \R$ such that $\eta(U) = \infty$ it follows that $\lim_{ n\rightarrow \infty}\eta_n(U) = \infty$.
		\item[(ii)] For all compact sets $K \subset \mathcal{R}_\eta$ there exists $M>0$ and $N \in \N$ such that $\eta_n(K) <M$ for all
		$n \geq N$.
		\item[(iii)] For all $\phi \in \mathcal{C}_c(\mathcal{R}_\eta)$ it holds that 
		\begin{equation*}
			\lim_{n \rightarrow \infty}\int \phi d\eta_n = \int \phi d\eta.
		\end{equation*} 
	\end{enumerate}
\end{dfn}
In the next sections we shall rely heavily on  the following result, which is taken from \cite[ Proposition 3.10]{Marcus99}, \cite[Theorem 4]{ LeGall96} and \cite[Section 2.1]{Barnes24}.
\begin{proposition}\label{prp:initial_trace_prop}
	Let $(A,\nu) \in \mathcal{T}$. Then there exists a unique solution to \eqref{eq:initial_trace}. Furthermore, the following holds:
	\begin{enumerate}
		\item Let $(A,\nu), (\hat{A},\hat{\eta})\in \mathcal{T}$ such that $(A,\nu) \preceq(\hat{A},\hat{\eta})$. Then 
		\begin{equation*}
			\mathcal{W}_t(A,\nu)(x) \leq 	\mathcal{W}_t(\hat{A},\hat{\eta})(x), \quad \forall t>0, x\in \R. 
		\end{equation*}
		\item Let $((A_n,\nu_n))_{n\geq 1} \subset \mathcal{T}$ such that $\eta^{A_n,\nu_n}$ converges $m$-weakly to $\eta^{(A,\nu)}$, then it follows that $\mathcal{W}(A_n,\nu_n)$ converges to $\mathcal{W}(A,\nu)$ uniformly on compacts of $(\R_+\backslash \{0\})\times \R$.
	\end{enumerate}  
\end{proposition}
In particular, by Lemma~\ref{lem:pde_est} and Proposition~\ref{prp:initial_trace_prop} it follows that $V_t(n\delta_x)(y)$ converges monotonically to $\mathcal{W}_t(x,0)(y)$ for all $t > 0$
and $x,y\in \R$. Therefore, 
\begin{equation}\label{eq:VW_ineq}
	V_t(n\delta_x)(y) \leq \mathcal{W}_t(x,0)(y), \forall t > 0, x,y\in\R, \quad n \geq 1
\end{equation}
which we shall frequently use in subsequent sections.
The following proposition is a version of  Theorem~1 from \cite{Brezis95}.
\begin{proposition}\label{prp:very_sing_asymp}
	It holds that 
	\begin{equation*}
		\mathcal{W}_t(\{x_0\}, 0)(x) = t^{-1} f(t^{-1/2} |x-x_0|),
	\end{equation*}
	where $f \colon [0,\infty) \rightarrow [0,\infty)$ is a smooth function satisfying
	\begin{equation*}
		f(x)  \asymp C\ee^{-\frac{ x^2}{2}}x\Big( 1+ o(x^{-2})\Big), 
	\end{equation*}
	as $x \rightarrow \infty$.
\end{proposition}

\section{Duality}\label{sec:duality}

In this section, we construct the dual process $Y$ and prove the
duality relation stated in Proposition~\ref{prop:duality}. The
construction is in many ways inspired by ideas from  \cite{Mytnik98} and \cite{Mytnik02}. We briefly describe the structure of the section.
We first construct the dual process $Y$ by alternating evolution
according to the log-Laplace equation with random jumps, some of which
may have infinite size. In Section~\ref{sec:artificial_branching}, we introduce an auxiliary
branching process that dominates the jump mechanism and provides moment
estimates for the number of jumps and for the time-integrated total
mass of $Y$. In Section~\ref{sec:convergence_dual}, we approximate the dual process by
truncating the jump sizes and prove convergence of the resulting
processes and duality functionals. Finally, in Section~\ref{sec:appr_dual}, we derive
the duality relation for the truncated processes and pass to the limit,
thereby proving Proposition~\ref{prop:duality}.
Since the jumps of $Y$ may take infinite values, we consider measures
on $\overline{\R}_+ = \R_+ \cup \{ \infty\}$, which is
endowed with the metric

	\begin{equation*}
		d(x,y) = |\arctan(x)-\arctan(y)|,
	\end{equation*}
	with the convention that $\arctan(\infty) = \pi/2$. 
	We denote
	\begin{equation*}
		\nu = \nu^1+ \nu^2+d_1\delta_{\infty} + d_2\delta_{\infty},
	\end{equation*}
	where $d_1,d_2$ are defined as in \eqref{eq:def_d} and 
	\begin{equation}\label{eq:def_nub}
		\nu([b,\infty]) = (\nu^1+\nu^2)([b,\infty)) + d_1+d_2, \quad b\in \R_+.
	\end{equation}
	 Furthermore,  assume that $Y_0 = \mu + m \delta_{x_0}$, 
	where $\mu$ is a nonnegative finite measure, $m \in \R_+ \cup\{\infty\}$, and $x_0 \in \R$.
	We construct the process iteratively. For this purpose, let $S_{1}$ be an exponentially distributed random variable with mean $1$ and denote 
	\begin{equation*}
		T_{1}= \inf\{ t> 0 | \nu((0,\infty])\int_0^t \la V_s(Y_0),1\ra ds  > S_{1}\}.
	\end{equation*}
	Furthermore, let $M_{1}$ be an independent random variable with values in $\{1,2\}$ given by
	\begin{equation}\label{eq:def_M_K_2}
			\P(M_{1} = i) = \frac{\nu^i((0,\infty)) +d_i}{\nu((0,\infty])}, \quad i \in \{1,2\},
	\end{equation}
	and let $Z_{1}$ be a random variable such that 
	\begin{equation*}
		\P(Z_{1} \geq b|M_{1}) = \frac{\nu^{M_{1}}([b,\infty))\mathbbm{1}_{b<\infty} + d_{M_1}}{\nu^{M_{1}}((0,\infty))+d_{M_1}}, \quad b\in \overline{\R}_+.
	\end{equation*}
	Finally we also define
	\begin{equation*}
		\P(U_{1} \in A| Y_{T_{1}-}) = \frac{\int_A Y_{T_{1}-}(x)dx}{\int_\R Y_{T_{1}-}(x)dx},
	\end{equation*}
	for $A\in \mathcal{B}(\R)$. Then we set 
	\begin{equation*}
		\begin{cases}
			Y_t = V_t(Y_0), \quad 0 <t < T_{1}, \\
			Y_{T_{1}} = V_{T_{1}-}(Y_0) + Z_{1} \delta_{U_{1}}, \quad \text{if } T_{1} <\infty.
		\end{cases}
	\end{equation*}
	Note that the process takes values in $\mathcal{M}_{reg}^+$ at time $T_1$.
	We then continue this construction, by letting $(M_{k})_{k = 2}^\infty$ be i.i.d. copies of $M_{1}$, letting $(S_{k})_{k = 1}^\infty$ be i.i.d. exponentially distributed random variables with mean $1$, independent of $(M_{k})_{k = 1}^\infty$. Now we define
	\begin{equation}\label{eq:def_time}
		T_{k} = \inf\{t > 0 | \nu((0,\infty]) \int_0^t\int_\R V_s(Y_{T_{k-1}})dxds  > S_{k}\}, \quad \bar{T}_{k} = \sum_{i = 1}^{k}T_{i},
	\end{equation}
	for $k  \geq 2$ and set $\bar{T}_0 = T_0 = 0$. Moreover, we define random variables $(U_{k})_{k = 1}^\infty$ as
	\begin{equation}\label{eq:pos_yn}
		\P(U_{k} \in A| Y_{\bar{T}_{k}-}) = \frac{\int_A Y_{\bar{T}_{k}-}(x)dx}{\int_\R Y_{\bar{T}_{k}-}(x)dx},
	\end{equation}
	for $A \in \mathcal{B}(\R)$ and
	\begin{equation}\label{eq:def_height}
				\P(Z_{k} \geq b | M_{k}) = \frac{\nu^{M_{k}}([b,\infty))\mathbbm{1}_{b<\infty} + d_{M_{k}}}{\nu^{M_{k}}((0,\infty)) + d_{M_{k}}} , \quad b\in \overline{\R}^+.
	\end{equation}
	Hereby,   $(Z_{k}, M_{k})$ and $(T_{k}, U_{k})$ are independent for each $k \geq 1$.
	Finally, we set 
	\begin{equation*}
		\mathcal{F}_t = \cap_{\epsilon >0} \sigma(Y_s(x), s \leq t+\epsilon, x\in \R)
	\end{equation*}
	and we choose  $(Z_{k})_{k \geq 1}$ and $(S_{k})_{k\geq 1}$ and $(M_{k})_{k\geq1}$ in such a way that $(Z_{\ell}, S_{\ell}, M_{\ell})_{\ell \geq k}$ is  independent of $\mathcal{F}_{\bar{T}_{k-1}}$ for all $k \geq 1$.\\
	We iteratively set
	\begin{equation}\label{eq:def_procc}
		\begin{cases}
			Y_t = V_{t-\bar{T}_{k-1}}(Y_{T_{k-1}}), \quad \bar{T}_{k-1}<t < \bar{T}_{k}, \\
			Y_{\bar{T}_{k}} = V_{T_{k}}(Y_{\bar{T}_{k-1}}) + Z_{k}\delta_{U_{k}}, \quad \text{ if } \bar{T}_{k} < \infty
		\end{cases}
	\end{equation}
	and continue this construction for all $k \geq 1$. Again, note that this process takes values in $\mathcal{M}_{reg}^+$.
	The next section  demonstrates that $\bar{T}_{k} \rightarrow \infty$ as $k \rightarrow \infty$.

	\subsection{An Artificial Branching Process}\label{sec:artificial_branching}
	The aim of this section is to derive bounds for $Y$, which  only depend on $Y_0$ and $\nu((0,\infty])$.
	 These bounds are crucial in order to prove that $Y$  only has finitely many jumps in each finite time interval and in order to show that all the terms in \eqref{eq:duality} have finite moments. 
	To this end, we introduce an artificial branching mechanism that will allow us to control the number of jumps of $Y$. This approach is inspired by the ideas from \cite{Barnes24}.\\
	\begin{enumerate}[label=\textup{(\roman*)},leftmargin=*,itemsep=0.9em]

\item \label{itm:br_process_1}
The process is a pure-birth branching process.  Once a particle is born, it
remains alive forever.  In every branching event exactly one new child is born, and the parent is not removed or relabelled.  

\item \label{itm:br_process_2}
Let
\[
    \mathcal{U}:=\bigcup_{m\ge 1}\N^m
\]
be the space of possible labels.  The initial particle has label $(1)$.  If
$\alpha=(\alpha_1,\ldots,\alpha_m)\in\mathcal{U}$, then the $i$-th child of the particle
with label $\alpha$ has label
\[
    (\alpha,i):=(\alpha_1,\ldots,\alpha_m,i).
\]
 Each actual particle is
specified by its label $\alpha$, its birth time $\tau_\alpha$, its profile
$Y^{(\alpha)}$, and its driving Poisson random measure $\mathfrak{M}^{(\alpha)}$.  Labels
in $\mathcal{U}$ that are not descendants of $(1)$ are merely potential labels; we set
their birth times equal to $\infty$, and they never enter the alive set defined
below.
For a label $\alpha = (\alpha_1, \ldots \alpha_n)$ we denote its length by 
\begin{equation*}
    |\alpha| = |(\alpha_1, \ldots \alpha_n)| = n.
\end{equation*}

\item \label{itm:br_process_3}
For each potential label $\alpha\in\mathcal{U}$, let $\mathfrak{M}^{(\alpha)}$ be a Poisson
random measure on
\[
    (0,\infty]\times(0,\infty)\times(0,\infty)\times\mathbb R
\]
with intensity
\begin{equation}\label{eq:base-intensity-enumerated}
    \nu(d\lambda) dt dr dz .
\end{equation}
The random measures $\{\mathfrak{M}^{(\alpha)}:\alpha\in\mathcal{U}\}$ are independent. 

\item \label{itm:br_process_4}
Let $\beta=(\alpha,i)$ be a child with finite birth time $\tau_\beta$.  Its
profile is defined by
\begin{equation}\label{eq:child-profile-enumerated}
	Y_t^{(\beta)}
	:= \mathcal{W}_{t-\tau_\beta}(0,0)\,\mathbbm{1}_{\{t>\tau_\beta\}},
	\qquad t\ge 0,
\end{equation}
where $\mathcal{W}(0,0)$ is defined as in \eqref{eq:very_sing_def}.  If $\tau_\beta=\infty$, set $Y^{(\beta)}\equiv 0$. We then denote the particle with label $(\beta)$ by $Z^{(\beta)} = (Y^{(\beta)}, \mathfrak{M}^{(\beta)})$.

\item \label{itm:br_process_5}
The particle with label $(1)$ is born at time
\[
    \tau_{(1)}:=0
\]
and carries the profile
\[
    Y_t^{(1)}=V_t(Y_0), \qquad t\ge 0 .
\]
Equivalently, the initial particle is $Z^{(1)}=(Y^{(1)},\mathfrak{M}^{(1)})$.

\item \label{itm:br_process_6}
Suppose that a particle with label $\alpha$ has already been born at a finite
time $\tau_\alpha$. Then we define the random measure
$\mathfrak{N}^{(\alpha)}$  by
\begin{equation}\label{eq:accepted-measure-enumerated}
\begin{aligned}
    \mathfrak{N}^{(\alpha)}(G_1,G_2,G_3,G_4)
    := {}& \int_{G_1}\int_{G_2}\int_{G_3}\int_{G_4}
    \mathbbm{1}_{\{t>\tau_\alpha\}}
    \mathbbm{1}_{\{r\le Y_{t-}^{(\alpha)}(z)\}}
    \,\mathfrak{M}^{(\alpha)}(d\lambda, dt, dr, dz),
\end{aligned}
\end{equation}
for Borel sets $G_1\subset   (0,\infty],G_2,G_3\subset(0,\infty)$ and $G_4\subset\mathbb R$.
The branching times of the particle $(\alpha)$ are the arrival times of 
\begin{equation*}
	\mathfrak{N}^{\alpha}_t = 	\mathfrak{N}^{\alpha}((0,\infty], (0,t] , (0,\infty), \mathbb{R}),\quad t \geq 0. 
\end{equation*}

\item \label{itm:br_process_7}
For $t\ge\tau_\alpha$, let
\begin{equation}\label{eq:birth-counting-process-enumerated}
    \mathfrak{N}_t^{(\alpha)}
    := \mathfrak{N}^{(\alpha)}\bigl( (0,\infty],(\tau_\alpha,t],(0,\infty),\mathbb R\bigr)
\end{equation}
be the number of children born to $\alpha$ by time $t$ and  let
\[
    \tau_{(\alpha,1)}\le \tau_{(\alpha,2)}\le\cdots
\]
be the arrival times of $(\mathfrak{N}^{(\alpha)}_{s})_{s\in (\tau_\alpha, \infty)}$. Note that the probability that $\tau_{\alpha,i } = \tau_{\alpha, j}$ for $i \neq j$ is zero.
Furthermore, define 
\begin{equation}\label{eq:compensator-enumerated}
    A_t^{(\alpha)}
    :=  \nu((0,\infty])\int_{\tau_\alpha}^{t}\int_{\mathbb R}
       Y_{s-}^{(\alpha)}(z) dzds,
    \qquad t\ge\tau_\alpha.
\end{equation}
After conditioning on
$\tau_\alpha$ and $Y^{(\alpha)}$, the children of $\alpha$ are born according to
an inhomogeneous Poisson process with cumulative rate $A_t^{(\alpha)}$.  Hence,
for $t\ge\tau_\alpha$,
\begin{equation}\label{eq:birth-time-tail-enumerated}
    \mathbb P\!\left(\tau_{(\alpha,i)}>t\,\middle|\,\tau_\alpha,Y^{(\alpha)}\right)
    = \exp\{-A_t^{(\alpha)}\}
      \sum_{k=0}^{i-1}\frac{\bigl(A_t^{(\alpha)}\bigr)^k}{k!}.
\end{equation}

\item \label{itm:br_process_8}
We denote the set of labels of particles alive at
time $t$ by
\begin{equation}\label{eq:defI_t}
    I_t^\nu(Y_0):=\{\alpha\in\mathcal{U}:\tau_\alpha\le t\} .
\end{equation}
Since particles never die, this is also the set of labels of particles born by
time $t$.  Its cardinality is denoted by
\[
    |I_t^\nu(Y_0)|:=\# I_t^\nu(Y_0).
\]
When there is no ambiguity, we suppress $\nu$ or $Y_0$  from the notation and write
$I_t(Y_0)$ or $I_t$ instead of $I_t^\nu(Y_0)$.
\end{enumerate}
	The next lemmas allow us to control the moments of the number of jumps of the dual process.
	\begin{lemma}\label{lem:number_jumps}
	For all $T>0$ there exists $C =C(T)>0$ such that 
		\begin{equation*}
			\E[|I_t(Y_0)|] < C,
		\end{equation*}
	\end{lemma}
	for all $t \leq T$.
	\begin{proof}
		Define 
		\begin{equation*}
			I_t^0 = \{ (1)\},\quad  I_t^{i} = \{ \alpha \in \mathcal{U} | |\alpha | = i+1, \la Y_t^{(\alpha)}, 1\ra >0\}, \quad i \geq 1,
		\end{equation*}
		where $|\alpha|$ denotes the length of the vector $\alpha$.
		Let $T_0>0$ be chosen sufficiently small and $t \leq T_0$. The precise value of $T_0$ will be determined later.  We will show that
		\begin{equation}\label{eq:number_gen_i}
			 \E{|I_t^i|} \leq 2^{-i}, \quad 0 \leq t \leq T_0.
		\end{equation}
		First, we prove \eqref{eq:number_gen_i} for the case $i = 1$. The expected number of children of the particle 
		$Z^{(1)}$ until time $t \in [0,T_0]$ is given by the expectation of $\mathfrak{N}^{(1)}_t =\mathfrak{N}^{(1)}((0,\infty] \times (0,t] \times (0,\infty)\times \R) $, that is 
		\begin{equation*}
			\begin{split}
				& \E{|I_t^1|} =  \E[\mathfrak{N}^{(1)}_t] \\
				&= \nu((0,\infty])\E\Big[\int_0^t \int_\R Y_{s}^{(1)}(x)dxds\Big] = \nu((0,\infty])\int_0^t \langle V_s(Y_0), 1 \rangle ds. 
			\end{split}
		\end{equation*}
		Since, by Lemma~\ref{lem:est_pde}, $(t,x) \mapsto V_t(Y_0)(x)$ is integrable on $[0,T] \times \R$ for all $T>0$, we can find $T_0>0$ sufficiently small so that
		\begin{equation*}
			\E{|I_t^1|} \leq 2^{-1}.
		\end{equation*}
		Now suppose that \eqref{eq:number_gen_i} holds for $i \geq 1$. 
			Denote by $(\tilde{\mathfrak{N}}^{(\alpha)})_{\alpha \in \mathcal{U}}$ a family of independent Poisson processes on $(0,\infty)$ with intensity given by   
		\begin{equation*}
			\nu((0,\infty])\int_\R \mathcal{W}_t(0,0)(x)dxdt.
		\end{equation*}
		Then it follows that 
		\begin{equation*}
		\begin{split}
		\E[|I_t^{{i+1}}|] &= \E\Big[ \sum_{\alpha \in I_t^{i}} \mathfrak{N}^{(\alpha)}((0,\infty] \times (0,t] \times (0,\infty)\times \R)\Big]\\
		& \leq \E\Big[ \sum_{\alpha \in I_t^{i}} \tilde{\mathfrak{N}}^{(\alpha)}_t\Big] \leq  \E[|I_t^{i}|] \nu((0,\infty]) \int_0^t\int_\R \mathcal{W}_s(0,0)(x)dxds \leq 2^{-i} 2^{-1}, 
		\end{split}
		\end{equation*}
		if we choose $T_0$ sufficiently small, so that 
		\begin{equation*}
			\nu((0,\infty]) \int_0^{T_0}\int_\R \mathcal{W}_s(0,0)(x)dxds \leq 2^{-1}.
		\end{equation*}
		Hereby, we used that the intensity of $\tilde{\mathfrak{N}}^{(\alpha)}$ bounds  the intensity of $\mathfrak{N}^{(\alpha)}((0,\infty], dt, (0,\infty), \R)$ from above for each $\alpha$.
		Thus, we arrive at 
		\begin{equation*}
			\E[|I_t(Y_0)|] = \sum_{i = 0}^\infty \E[|I_t^i|] \leq 1 + \sum_{i = 1}^\infty 2^{-i},
		\end{equation*}
		 for all $0 \leq t \leq T_0$.\\
		Now, we need to prove the claim for general $T>0$. For every $T >0$ we can find $M \in \N$ such that $T \leq MT_0$. 
		Then, using the fact that no particle dies, we can estimate
		\begin{equation*}
			\E[|I_T(Y_0)|] \leq \E[|I_{MT_0}(Y_0)|] \leq \E[\sum_{\alpha \in I_{T_0}}|\tilde{I}_{T_0,MT_0}^{\alpha}|],
		\end{equation*} 
		where 
		\begin{equation*}
			\tilde{I}^{\alpha}_{s,t} = \{ \beta \in \mathcal{U} : \beta \text{ is a descendant of } \alpha \text{ born in } [s,t] \} \cup \{\alpha\}.
		\end{equation*}
		Note that for $\alpha \neq (1)$ it holds that
		\begin{equation*}
			\E[|I^{\alpha}_{T_0,MT_0}|] \leq \E[|I_{(M-1)T_0}(\mathcal{W})|], 
		\end{equation*}
		where $I_{t}(\mathcal{W})$ denotes the particles that are alive at time $t$ in the branching process described in \ref{itm:br_process_1}-\ref{itm:br_process_8}, where $Y^{(1)}$ is replaced by $\mathcal{W}(0,0)$.  This simply follows from the fact that 
		\begin{equation*}
			\int_0^t\la \mathcal{W}_{s+r}(0,0),1\ra dr\leq \int_0^t\la \mathcal{W}_{r}(0,0),1\ra dr, \quad \forall s,t \geq 0, 
		\end{equation*}
		which follows from Proposition~\ref{prp:very_sing_asymp}.
		Thus, we iteratively get that
		\begin{equation*}
		\begin{split}
				&\E[|I_T(Y_0)|] \leq \E[|I_{MT_0}(Y_0)|] \\
				&\leq \E[|I_{T_0}(Y_0)|-1] \E[|I_{(M-1)T_0}(\mathcal{W})|] + \E[|I_{(M-1)T_0}(V_{T_0}(Y_0))|] \leq \ldots\leq C(T).
		\end{split}
		\end{equation*}

This finishes the proof.
	\end{proof}

	\begin{lemma}\label{lem:exponential1}
		Let $Y_0 = \mu + m\delta_{x_0}$ with $\mu \in \mathcal{M}_F^+$, $x_0 \in \R$ and $m \in \N\cup\{ \infty\}$. Furthermore, let $I_t(Y_0)$ be defined as in \eqref{eq:defI_t}. Then for every $\gamma > 0$, $K>0$, there exists $T_0 >0$ such that 
		\begin{equation}\label{eq:exponI}
			\sup_{\nu: \nu((0,\infty]) \leq K} \E[ \ee^{\gamma |I_t(Y_0)|}] < \infty,
		\end{equation}
	for all $t \leq T_0$ . Hereby, $T_0>0$ depends only on $K$, $\la \mu,1\ra$ and $\gamma$.
	\end{lemma}

	\begin{proof}
		 Fix arbitrary $\gamma>0$, $K>0$ and $t \leq T_0$, where $T_0$ is sufficiently small, the precise value will be chosen later. Define $\eta_{t}^{i,k}$ for $i,k \in \N$ to be i.i.d. Poisson distributed random variables with parameter 
		 \begin{equation}\label{eq:lambda_intensity}
		 	\lambda  = K\Big( \int_0^t \la \mathcal{W}_s(0,0),1 \ra ds +\int_0^t \la V_s(Y_0),1 \ra ds\Big).
		 \end{equation} 
	Furthermore,  set $Z_0 =1$ and
	\begin{equation*}
		 Z_n^t= \sum_{k = 1}^{Z_{n-1}^t} \eta_t^{n,k}, \quad n \geq 1.
	\end{equation*}
	Let $\nu$ be a measure such that $\nu((0,\infty]) \leq K$.
	Then it easily follows that $|I_t^i|$ is stochastically dominated by $Z_i^t$ for each $i\geq 1$. Thus, we can also infer that 
	\begin{equation*}
		|I_t(Y_0)| = \sum_{i = 0}^\infty |I_t^i| 
	\end{equation*}
	is stochastically dominated by 
	\begin{equation*}
		Z^t = \sum_{i= 0}^\infty Z_i^t.
	\end{equation*}
	The random variable $Z^t$ is nothing else but the total progeny of $(Z^t_n)_{n\geq 0}$ and its law can be made explicit by using the Otter-Dwass formula (see \cite{Dwass69}), whenever $\lambda \leq1$. By \eqref{eq:lambda_intensity} we can choose $T_0$ small enough such that 
	\begin{equation}\label{eq:cond_lambda}
		\lambda \leq 1, \quad \log(\lambda) < \lambda -\gamma -1, \quad \forall t \leq T_0.
	\end{equation}
	Then we get 
	\begin{equation*}
		\P(Z^t = k) = \ee^{-\lambda k} \frac{(\lambda k)^{k-1}}{k!},  \quad k \geq 1,
	\end{equation*}
	that is $Z^t$ has the Borel-Tanner distribution.
	Using Stirling's approximation, we thus  get that 
	\begin{equation*}
		\P(Z^t = k)   \asymp \sqrt{2\pi}\ee^{(1-\lambda) k} \lambda^{k-1} k^{-3/2}  \asymp \sqrt{2\pi}\ee^{(1-\lambda +\log(\lambda))k - \log(\lambda)}  k^{-3/2}. 
	\end{equation*}
	Thus, 
	\begin{equation*}
	\sqrt{2\pi} \ee^{-\log(\lambda)}\sum_{k = 1}^\infty k^{-3/2} \ee^{(\gamma +1-\lambda +\log(\lambda)) k } <\infty \quad \Rightarrow \quad 	\sum_{k = 1}^\infty e^{\gamma k} \P(Z^t = k) <\infty,
	\end{equation*}
	and 
	\begin{equation*}
		\sum_{k = 1}^\infty k^{-3/2} \ee^{(\gamma +1-\lambda +\log(\lambda)) k } <\infty,
	\end{equation*}
	by \eqref{eq:cond_lambda}.
	Therefore, we obtain that 
	\begin{equation*}
		 \E[\ee^{\gamma Z^t}] < \infty.
	\end{equation*}
	Note that the function $x \mapsto \ee^{\gamma x}$ is non-decreasing. Thus, the fact that $|I_t(Y_0)|$ is stochastically dominated by $Z^t$ implies that 
	\begin{equation*}
		\E[\ee^{\gamma |I_t(Y_0)|}] \leq	 \E[\ee^{\gamma Z^t}] < \infty.
	\end{equation*}
	Since $\nu$ such that $\nu((0,\infty]) \leq K$ was arbitrary, this completes the proof.
	\end{proof}
	\noindent Now, we need to relate the branching process to the process $Y$. 
	Since the birth times of new particles in the branching process happen at a higher rate than the jumps of $ t \mapsto \la Y_t,1\ra$, we aim to control $\int_0^t\la Y_s,1\ra ds$, which is sufficient for our arguments. 
	To this end, we define a coupling between the branching process and $Y$.
    Denote by $\bar{R}_i$ the time of the $i$-th branching event for $i \geq 1$ and set $\bar{R}_0 = 0$. It follows that 
    \begin{equation*}
       \nu((0,\infty]) \int_{\bar{R}_{i-1}}^{\bar{R}_i}\sum_{\alpha \in I_{R_{i-1}(Y_0)}} \langle Y^{(\alpha)}_s,1 \rangle ds, \quad i \geq 1
    \end{equation*}
    are independent exponential random variables with mean $1$.
    Recall that $(S_i)_{i = 1}^\infty$ is the sequence of independent $\exp(1)$-distributed random variables used for the construction of $Y$ as defined at the beginning of Section~\ref{sec:duality}.
	   That means we can define the branching process and $Y$ on a common probability space in such a way that  $(\bar{R}_i)_{i = 1}^\infty$ are given as follows 
	\begin{equation*}
		\begin{split}
			&\bar{R}_{0} = 0,\\
			&\bar{R}_{i} = \inf\{ t > \bar{R}_{i-1} \colon \sum_{\alpha \in I_{\bar{R}_{i-1}}}\int_{\bar{R}_{i-1}}^t \la Y^{(\alpha)}_s,1\ra ds \nu((0,\infty]) > S_{i}\}, \quad i \geq 1.
		\end{split}
	\end{equation*}
	Furthermore, we set 
	\begin{equation*}
		R_i = \bar{R}_i-\bar{R}_{i-1}, \quad i \geq 1. 
	\end{equation*}
	and
	\begin{equation*}
		\hat{Y}_t = \sum_{\alpha \in I_t} Y_t^{(\alpha)}, \quad t \geq 0.
	\end{equation*}
	Then, by construction, we get the following identity
	\begin{equation}\label{eq:branchin_identity}
		\hat{Y}_t = V_t(Y_0) + \sum_{k = 1}^{|I_t(Y_0)|-1} \mathcal{W}_{t-\bar{R}_k}(0,0).
	\end{equation}
	Using this construction, we can prove our first estimates. 
	\begin{lemma}\label{lem:coupling}
		For all $i \geq 1$ the following inequalities hold.
		\begin{enumerate}
			\item[(i)]  It holds almost surely that
			\begin{equation}\label{eq:ineq_times}
				R_{i} \leq T_{i}.
			\end{equation}
			\item[(ii)] It holds almost surely that
			\begin{equation}\label{eq:ineq_V}
				\la Y_{\bar{T}_{i-1}+s},1\ra \leq \la \hat{Y}_{\bar{R}_{i-1}+s},1\ra, \quad \forall 0 <s \leq T_{i}.  
			\end{equation}
			\item[(iii)] \label{itm:ineq_Y_int}
			It holds almost surely that
			\begin{equation}\label{eq:sec_step}
				\int_0^t \la Y_s,1\ra ds \leq \int_0^t\la \hat{Y}_s,1 \ra ds, \quad \forall t \in (\bar{T}_{i-1}, \bar{T}_{i}).
			\end{equation}
		\end{enumerate}
	\end{lemma}
	\begin{proof}
		We prove the claim by induction over $i\geq 1$. 
		It is easy to see that $R_1 = T_1$ and then \eqref{eq:ineq_V} and \eqref{eq:sec_step} are immediate for $i = 1$. 
		Suppose that \eqref{eq:ineq_times}, \eqref{eq:ineq_V} and \eqref{eq:sec_step} hold for all $ j \leq i$ for some $i \geq 1$. We now show that this implies the claim for $i+1$. 
		First, we prove \eqref{eq:ineq_V}. We apply Lemma~\ref{lem:pde_est} repeatedly to arrive at
		\begin{equation*}
			\begin{split}
				\la Y_{\bar{T}_{i}+s},1\ra &\leq \la V_{s}(Y_{\bar{T}_{i}-}),1\ra + \la\mathcal{W}_s(0,0),1\ra\\
				& \leq \la V_{s+T_i}(Y_{\bar{T}_{i-1}-}),1\ra + \la\mathcal{W}_{s+\bar{T}_i}(0,0),1\ra  + \la\mathcal{W}_{s}(0,0),1\ra\\
				&\leq \la V_{s+\bar{T}_i}(Y_0),1\ra + \sum_{k = 1}^i \la \mathcal{W}_{s+\bar{T}_i-\bar{T}_k}(0,0), 1\ra, \quad a.s.
			\end{split}
		\end{equation*}
	By the induction hypothesis, we know that $R_j \leq T_j$ for all $j \leq i$. Thus, we have 
	\begin{equation*}
		\bar{T}_{i}- \bar{T}_{k} = \sum_{\ell = k+1}^i T_{\ell} \geq \bar{R}_i-\bar{R}_k,
	\end{equation*}
	for all $0 \leq k \leq i $. 
	We have the monotonicity relations
	\begin{align}
			&\la \mathcal{W}_s(0,0), 1\ra \leq 	\la \mathcal{W}_t(0,0), 1\ra \quad 0< t\leq s, \label{eq:monotonicity1}\\
				&\la V_s(Y_0), 1\ra \leq 	\la V_t(Y_0), 1\ra \quad 0< t\leq s. \label{eq:monotonicity2}
	\end{align}
	 Hereby, \eqref{eq:monotonicity1} follows from Proposition~\ref{prp:very_sing_asymp}  and \eqref{eq:monotonicity2} follows from integrating \eqref{eq:sol_rep_V}. 
	We arrive at
	\begin{equation}\label{eq:est_y_t_s}
		\begin{split}
			&\la Y_{\bar{T}_{i}+s},1\ra \leq \la V_{s+\bar{R}_i}(Y_0),1\ra + \sum_{k = 1}^i \la \mathcal{W}_{s+\bar{R}_i-\bar{R}_k}(0,0), 1\ra\\
			& \leq \la \hat{Y}_{\bar{R}_i+s},1\ra \quad a.s.,
		\end{split}
	\end{equation}
	for all $0<s <T_{i+1}$.  \\
	Now, we will show
	\begin{equation}\label{eq:induc_times_i+1}
	 R_{i+1} \leq T_{i+1}.
	 \end{equation}
	  Recall that
	\begin{equation*}
		T_{i+1} = \inf\{t \geq 0 \colon \nu((0,\infty]) \int_0^t \la Y_{\bar{T}_i +s},1\ra ds > S_{i+1} \},
	\end{equation*}
	and 
	\begin{equation*}
		R_{i+1} = \inf\{t \geq 0 \colon \nu((0,\infty]) \int_0^t \la \hat{Y}_{\bar{R}_i +s},1\ra ds > S_{i+1} \}.
	\end{equation*}
	Therefore, \eqref{eq:induc_times_i+1} is immediate by \eqref{eq:ineq_V} for $i+1$. \\
	Finally, we need to verify
	\begin{equation*}
		\int_0^t \la Y_s,1\ra ds \leq \int_0^t\la \hat{Y}_s,1 \ra ds, \quad \forall t \in (\bar{T}_{i}, \bar{T}_{i+1}).
	\end{equation*}
	Fix arbitrary $t \in (\bar{T}_{i}, \bar{T}_{i+1})$ 
	 and define 
	\begin{equation*}
	j^\ast = \sup\{ k \geq 1 \colon \bar{R}_k \leq t\}. 
	\end{equation*}
	Since $t \in  (\bar{T}_{i}, \bar{T}_{i+1})$ and $\bar{T}_k \geq \bar{R}_k$ for all $k \leq i+1$ by \eqref{eq:induc_times_i+1} we get 
	$j^\ast \geq i$.
	Note that 
	\begin{equation*}
	\begin{split}
		&\nu((0,\infty])\int_0^t \la Y_s,1\ra ds\\
		& =\nu((0,\infty]) \Big( \int_0^{\bar{T}_i} \la Y_s,1\ra ds + \int_{\bar{T}_i}^t \la Y_s,1\ra ds \Big)= \sum_{k =1}^i S_k + \nu((0,\infty])\int_{\bar{T}_i}^t \la Y_s,1\ra ds,
	\end{split}
	\end{equation*}
	and 
	\begin{equation*}
	\begin{split}
			&\nu((0,\infty])\int_0^t \la \hat{Y}_s,1\ra ds\\
			& = \nu((0,\infty]) \Big(\int_0^{\bar{R}_{j^\ast}} \la \hat{Y}_s,1\ra ds + \int_{\bar{R}_{j^\ast}}^t \la \hat{Y}_s,1\ra ds\Big) = \sum_{k = 1}^{j\ast} S_k + \nu((0,\infty]) \int_{\bar{R}_{j^\ast}}^t \la \hat{Y}_s,1\ra ds.
	\end{split}
	\end{equation*}
	Thus, when $j^\ast > i$, the proof is finished. In the case when $j^\ast = i$ it is enough to verify that
	\begin{equation*}
		 \int_{\bar{T}_i}^t \la Y_s,1\ra ds \leq  \int_{\bar{R}_{i}}^t \la \hat{Y}_s,1\ra ds.
	\end{equation*}
	Using \eqref{eq:ineq_V} for $i+1$, which follows from \eqref{eq:est_y_t_s}, and using that $\bar{R}_i \leq \bar{T}_i$, we get
	\begin{equation*}
		\begin{split}
			&  \int_{\bar{T}_i}^t \la Y_s,1\ra ds = \int_0^{t-\bar{T}_i} \la Y_{\bar{T}_i+s},1\ra ds \leq \int_0^{t-\bar{T}_i} \la \hat{Y}_{\bar{R}_i+s},1\ra ds\\
			& \quad \leq \int_0^{t-\bar{R}_i} \la \hat{Y}_{\bar{R}_i+s},1\ra ds =   \int_{\bar{R}_i}^t \la \hat{Y}_s,1\ra ds.
		\end{split}
	\end{equation*}	
	This finishes the proof. 
	\end{proof}
With these estimates at hand, we can prove the following proposition.
\begin{proposition}\label{prop:main}
	Let $Y$ be the process constructed in \eqref{eq:def_procc} with $Y_0^m = \mu+m\delta_{x_0}$, where $\mu \in \mathcal{M}_F^+$, $x_0 \in\R$ and $m \in \N \cup \{ \infty\}$ and denote by $\mathbb{E}_{Y_0^m}$ the expectation with respect to the law of the process $Y$ started from $Y_0^m$.
  Then it follows for all $T>0$ that
	\begin{equation}\label{eq:prop_eq}
		\sup_m \E_{Y_0^m}\Big[\int_0^T\la Y_t,1\ra dt\Big] \leq C(T),
	\end{equation}
where $C(T)$ depends only on $\nu((0,\infty])$ and $\la \mu,1\ra$.
\end{proposition}
\begin{proof}
	Since $\bar{R}_i \rightarrow \infty$ as $i \rightarrow \infty$ we get by the monotone convergence theorem and using \eqref{eq:sec_step} that
	\begin{equation*}
		\begin{split}
			& \E_{Y_0^m}\Big[ \int_0^t \la Y_s, 1\ra ds\Big] = \lim_{i \rightarrow \infty}\E_{Y_0^m}\Big[ \int_0^t \la Y_s, 1\ra ds \mathbbm{1}_{t \leq \bar{T}_{i}}\Big] \\
			& \leq \E_{Y_0^m}[\int_0^{t}\la \hat{Y}_s,1\ra ds].
 		\end{split}
	\end{equation*}
	By using Lemma~\ref{lem:number_jumps},  this immediately yields the following inequality
	\begin{equation*}
		\begin{split}
			&\E _{Y_0^m}\Big[ \int_0^t \la Y_s, 1\ra ds\Big] \leq \E_{Y_0^m} \Big[ \int_0^t \sum_{\alpha \in I_s} \la Y_s^{(\alpha)},1\ra ds\Big]\\
			& \leq  \E_{Y_0^m} \Big[ \int_0^t \la Y_s^{(1)},1\ra ds\Big] + \E \Big[ \sum_{\alpha \in I_t}\int_0^t \la \mathcal{W}_s(0,0), 1\ra ds\Big] \\
			& \leq \int_0^t\la V_s(Y_0^m),1 \ra ds + C(t)\E[|I_t|] \leq C(t),
		\end{split}
	\end{equation*}
	where $C(t)>0$ depends only on $\nu((0,\infty])$ and $\la \mu,1\ra$.
\end{proof}

\begin{lemma}\label{lem:exponential_moment}
	Let $Y_0^m = \mu+m\delta_{x_0}$, where $\mu \in \mathcal{M}_F^+$, $x_0 \in\R$ and $m \in \N \cup \{ \infty\}$. Let $Y(\nu)$ be the process constructed above with measure $\nu$. Then, for every $\gamma \in \R$ and $K>0$, there exists $T_0>0$ such that for all $t\leq T_0$ we have 
	\begin{equation}\label{eq:est_lem65}
		\sup_{\nu: \nu((0,\infty]) \leq K} \sup_m\E_{Y_0^m}\Big[\exp\Big(-\gamma\int_0^t\la Y_s(\nu) ,1\ra ds\Big)\Big] < C, 
	\end{equation}
	Hereby, $T_0>0$ and $C>0$ depend only on $\la \mu,1\ra$, $K$ and $\gamma$.
\end{lemma}
\begin{proof}
	If $\gamma\geq0$, the statement follows immediately. Thus, we only need to consider the case $\gamma<0$. 
	By Lemma~\ref{lem:exponential1} and \eqref{eq:sec_step} and using that
	\begin{equation*}
		\int_0^t\la Y^{\alpha}_s(\nu),1\ra ds \leq \int_0^t\la \mathcal{W}_s(0,0),1\ra ds , \quad \forall \alpha \in I_t, t \geq0, 
	\end{equation*}
	we get that there exists $T_0>0$ such that
	\begin{equation*}
    \begin{split}
		&\E_{Y_0^m}\Big[\exp\Big(-\gamma\int_0^t\la Y_s(\nu) ,1\ra ds\Big)\Big]\\
        &\leq  \E\Big[\exp\Big(-\gamma|I_t^\nu|\int_0^t\la \mathcal{W}_s(0,0) ,1\ra ds -\gamma \int_0^t \la V_s(\mu),1\ra ds\Big)\Big]  < C,
    \end{split}
	\end{equation*}
	for all $t \leq T_0$, $m\geq 1$ and $\nu$ with $\nu((0,\infty])\leq K$.
	This finishes the proof. 
\end{proof}

\subsection{Convergence of the Dual Process}\label{sec:convergence_dual}
Let $\nu$ be as in \eqref{eq:def_nub}. For $n\in \mathbbm{N}\cup \{\infty\} $ and $m\in \R_+\cup \{\infty\}$ denote by $(Y_t^{n,m}(Y_0^m))_{t\geq 0}$ the process as constructed in \eqref{eq:def_procc} 
with initial condition
$Y_0^m = \mu+ m \delta_{x_0}$, where $\mu\in \mathcal{M}_F^+$ and the measure $\nu$ is replaced by 
\begin{equation*}
	\nu_n = \nu^1_{|[0,n]} +\nu^2_{|[0,n]} + d_1\delta_n + d_2\delta_n =  \nu^1_{n} +\nu^2_{n} + d_1\delta_n + d_2\delta_n,
\end{equation*}
for $n\in \N$. Here, $d_1$, $d_2$ and $a$ are defined as in \eqref{eq:def_d} and for $n= \infty$ we set $\nu_\infty \equiv \nu$. 
We also assume that $\la Y_0^m, 1\ra>0$.
For $i \geq 1$, we denote by $\bar{T}_{n,i}^m$ the jump times, by $T_{n,i} = \bar{T}_{n,i}^m- \bar{T}_{n,i-1}^m$ the time between jumps, by $Z_{n,i}^m$ the jump sizes and by $U_{n,i}^m$ the jump positions defined as in \eqref{eq:def_time}-\eqref{eq:def_height}. Moreover, we denote by $M_{n,i}^m$
the random variables defined as in \eqref{eq:def_M_K_2}. 
Finally, we shall denote by $(S_{n,k}^m)_{k\geq1}$ the sequence independent of $\exp(1)$-distributed random variables from the construction of the process $Y^{n,m}$ in the beginning of  Section~\ref{sec:duality}. Let $X_0$ be a bounded, continuous and nonnegative function on $\R$.\\
In the remainder of the section, we prove that 
\begin{equation*}
	\begin{split}
		&\lim_{m \rightarrow \infty} \E\Big[(-1)^{J_t^{m,m}}\exp\Big(-\la Y_t^{m,m}, X_{0}\ra-a \int_0^t\la Y_s^{m,m}, 1 \ra ds\Big)\Big]\\
		&\quad  = \E\Big[(-1)^{J_t^{\infty,\infty}}\exp\Big(-\la Y_t^{\infty,\infty}, X_{0}\ra-a \int_0^t\la Y_s^{\infty,\infty}, 1 \ra ds\Big)\Big],
	\end{split}
\end{equation*}
\begin{equation*}
	\begin{split}
		&\lim_{n \rightarrow \infty} \E\Big[(-1)^{J_t^{n,m}}\exp\Big(-\la Y_t^{n,m}, X_{0}\ra-a \int_0^t\la Y_s^{n,m}, 1 \ra ds\Big)\Big]\\
		&\quad  = \E\Big[(-1)^{J_t^{\infty,m}}\exp\Big(-\la Y_t^{\infty,m}, X_{0}\ra-a \int_0^t\la Y_s^{\infty,m}, 1 \ra ds\Big)\Big],
	\end{split}
\end{equation*}
and
\begin{equation*}
	\begin{split}
		&\lim_{m \rightarrow \infty} \E\Big[(-1)^{J_t^{\infty,m}}\exp\Big(-\la Y_t^{\infty,m}, X_{0}\ra-a \int_0^t\la Y_s^{\infty,m}, 1 \ra ds\Big)\Big]\\
		&\quad  = \E\Big[(-1)^{J_t^{\infty,\infty}}\exp\Big(-\la Y_t^{\infty,\infty}, X_{0}\ra-a \int_0^t\la Y_s^{\infty,\infty}, 1 \ra ds\Big)\Big],
	\end{split}
\end{equation*}
for $t \in [0,T_0]$. Here $T_0$ is chosen as in Lemma~\ref{lem:exponential_moment} with $\gamma =2a$ and $K = \nu((0,\infty])$.\\
Recall that by Lemma~\ref{lem:exponential_moment}, $T_0$ can be chosen independent of $n,m$ since 
\[ \sup_{n\in\mathbb N\cup\{\infty\}}
\nu_n((0,\infty])
=
\nu((0,\infty]).
\]

Before we begin the proof, we introduce the following notation.
For $\pmb{x} = (x_1, \ldots x_k)$ for some $k \geq 2$ we denote
$\overset{\leftarrow}{\pmb{x}} = (x_1, \ldots x_{k-1})$.  
Furthermore, for fixed $\mu \in \mathcal{M}_F^+$ with $\mu \neq 0$ define
\begin{equation*}
	\mathfrak{V}^1(u,z,t) = V_{t}(\mu +z\delta_u), \quad  (u,z,t) \in\R \times \overline{\R}_+ \times (0,\infty)
\end{equation*}
and for $k > 1$ we define recursively
\begin{equation*}
	\mathfrak{V}^k(\pmb{u},\pmb{z},\pmb{t}) = V_{t_k} (	\mathfrak{V}^{k-1}(\overset{\leftarrow}{\pmb{u}}, \overset{\leftarrow}{\pmb{z}} ,\overset{\leftarrow}{\pmb{t}} ) +z_k\delta_{u_k}), \quad (\pmb{u}, \pmb{z}, \pmb{t}) \in \R^k \times \overline{\R}_+^k \times (0,\infty)^k.
\end{equation*}
Recall that $ \overline{\R}_+ = \R_+ \cup \{ \infty\}$ is endowed with the metric 
\begin{equation*}
	d(x,y) = |\arctan(x)-\arctan(y)|, \quad x,y \in \bar{\R}_+,
\end{equation*}
with the convention that $\arctan(\infty) = \pi/2$. We will use $\mathfrak{V}^k$ later as follows. 
For all $k>1$, when $\pmb{u} = (x_0, U_{n,1}^m, \ldots U_{n,k-1}^m)$, $\pmb{z} = (m, Z_{n,1}^m, \ldots Z_{n,k-1}^m)$ and $\pmb{t} = (T_{n,1}^m, \ldots T_{n,k-1}^m,t)$, then 
\begin{equation}\label{eq:frakV_k_identity}
	\mathfrak{V}^k(\pmb{u}, \pmb{z},\pmb{t})(x) = Y^{n,m}_{\bar{T}_{n,k-1}^m+t}(x), \quad x\in \R,
\end{equation}
and for all $t\in (0,T_{n,k}^m)$.
\begin{lemma}
	For every $k \in \N$ the map
	\begin{equation*}
		\R^k \times \overline{\R}_+^k \times (0,\infty)^k \ni (\pmb{u}, \pmb{z}, \pmb{t})  \mapsto 	\mathfrak{V}^k(\pmb{u},\pmb{z},\pmb{t}) 
	\end{equation*}
	is continuous with respect to the topology induced by uniform convergence on compacts. 
\end{lemma}
\begin{proof}
	First, let us demonstrate that $\mu_n = \psi_n + z_n\delta_{u_n}$, $n\geq 1$ converges $m$-weakly to $\mu \in \mathcal{M}_{reg}^+$, 
	if $(\psi_n)_{n\geq 1}\subset \mathcal{C}(\R)$ converges uniformly on compacts  to $\psi \in \mathcal{C}(\R)$, $(z_n)_{n\geq 1}\subset \bar{R}_+$ converges to $z$
	with respect to $d$ and $(u_n)_{n\geq 1}\subset \R$ converges to $u$. 
	Clearly, $S_{\mu}$ is empty when $z \neq \infty$ and when $z = \infty$ we have $S_\mu = \{u\}$. 
	It is easy to verify that in the latter case $\eta^{(\emptyset, \mu_n)}$ converges $m$-weakly to $\eta^{u,\psi}$.
	By Proposition~\ref{prp:initial_trace_prop}, this demonstrates that 
	\begin{equation*}
		V(\psi_n + z_n\delta_{u_n}) \rightarrow V(\psi+ z\delta_{u}) \
	\end{equation*}
	uniformly on compacts of $(0,\infty) \times \R$. Furthermore, let $\{t_n\}_{n\geq 1} \subset (0, \infty)$ be a sequence that converges to $t \in (0, \infty)$. Then it follows that 
	\begin{equation*}
		V_{t_n}(\psi_n + z_n\delta_{u_n}) \rightarrow V_t(\psi+ z\delta_{u}), \quad \text{as } n\rightarrow \infty,
	\end{equation*}
	uniformly on compacts of $\R$. 
	Now it is straightforward to see that the continuity of $\mathfrak{V}^k$ can be proven inductively. 
\end{proof}  
\noindent We furthermore introduce, for $k \geq 1$, the map 
 \begin{equation*}
 	\Theta^k \colon \R^{k+1}\times (0,\infty)^k \times \overline{\R}_+^{k+1}  \times (0,\infty) \rightarrow \R_+,
 \end{equation*}
 given by
 \begin{equation*}
 	(\pmb{u},\pmb{t}, \pmb{z},s) \mapsto \inf \{ r \geq 0 \colon \int_0^r \la V_{\tau}(\mathfrak{V}^k( \overset{\leftarrow}{\pmb{u}}, \overset{\leftarrow}{\pmb{z}}, \pmb{t} ) +z_{k+1}\delta_{u_{k+1}}) ,1\ra d\tau >s\}.
 \end{equation*}
 Note that when $\pmb{u} = (x_0, U_{n,1}^m, \ldots U_{n,k}^m)$, $\pmb{t} = (T_{n,1}^m, \ldots T_{n,k}^m)$ and $\pmb{z} =(m,Z_{n,1}^m, \ldots Z_{n,k}^m)$ we have
 \begin{equation}\label{eq:phikexplain}
 	\Theta^k(\pmb{u}, \pmb{t}, \pmb{z}, S_{n,k+1}^m/\nu_n((0,\infty])) = T^m_{n,k+1}.
 \end{equation}
\begin{lemma}\label{lem:theta_cont}
	For every $k\geq 1$ the map $\Theta^k $
	is continuous on $\R^{k+1}\times (0,\infty)^k \times \overline{\R}_+^{k+1}  \times (0,\infty)$.
\end{lemma}
\begin{proof}
	Let  $(\pmb{u},\pmb{t}, \pmb{z},s)\in \R^{k+1}\times (0,\infty)^k \times \overline{\R}_+^{k+1}  \times (0,\infty) $, denote $\pmb{x} =  (\pmb{u},\pmb{t}, \pmb{z},s)  $ and let $B_r(\pmb{x})$ be defined as 
	\begin{equation*}
		B_r(\pmb{x}) = \{ \pmb{y} \in \R^{k+1}\times (0,\infty)^k \times \overline{\R}_+^{k+1}  \times (0,\infty)  \colon d_{k+1}(\pmb{x},\pmb{y}) \leq r\},
	\end{equation*}
	where $d_{k+1}$ is the canonical metric on $\R^{k+1}\times (0,\infty)^k \times \overline{\R}_+^{k+1} \times (0,\infty)  $.
	First, we prove that
	\begin{equation}\label{eq:supPhi}
		T := \sup_{\pmb{y} \in B_1(\pmb{x})} \Theta^k(\pmb{y}) < \infty.
	\end{equation}
	This claim can easily be verified by denoting $\pmb{y} = (\pmb{v}, \pmb{\tau},\pmb{w}, b)\in \R^{k+1}\times (0,\infty)^k \times \overline{\R}_+^{k+1}  \times (0,\infty) $ and noting that
	\begin{equation}
		\begin{split}
			&\sup_{\pmb{y} \in B_1(\pmb{x})} \Theta^k(\pmb{y}) = \sup_{\pmb{y} \in B_1(\pmb{x})}  \inf \{ r \geq 0 \colon \int_0^r \la V_t(\mathfrak{V}^k(\overset{\leftarrow}{\pmb{v}}, \overset{\leftarrow}{\pmb{w}}, \pmb{\tau})+ w_{k+1}\delta_{v_{k+1}}),1\ra dt> b\}\\
			&\leq  \sup_{\pmb{y} \in B_1(\pmb{x})}  \inf \{ r \geq 0 \colon \int_0^r \la V_t(\mathfrak{V}^k(\overset{\leftarrow}{\pmb{v}}, \overset{\leftarrow}{\pmb{w}}, \pmb{\tau})+ w_{k+1}\delta_{v_{k+1}}),1\ra dt> s+1\} \\
			& \leq  \sup_{\pmb{y} \in B_1(\pmb{x})}  \inf \{ r \geq 0 \colon \int_0^r \la V_{t+\bar{t}}(\mu),1\ra dt> s+1\} < \infty,
		\end{split}
	\end{equation}
	where 
	\begin{equation}\label{eq:bart}
		\bar{t} = \sum_{i = 1}^kt_i+1.
	\end{equation}
	Hereby, we used Lemma~\ref{lem:pde_est} and the fact that by Lemma~\ref{lem:V_asymp} we have
	\begin{equation}\label{eq:vmon}
		\la V_t(\mu),1\ra \leq \la V_s(\mu),1\ra,\quad \forall 0\leq s \leq t \text{ and } \int_0^\infty \la V_s(\mu),1\ra ds = \infty,
	\end{equation}
	for all finite, non-zero measures $\mu$.
	Now, let $\epsilon >0$ be arbitrary, $\delta>0$ sufficiently small and $\pmb{y}   = (\pmb{v},\pmb{\tau}, \pmb{w},b)\in B_\delta(\pmb{x})$.
	Then we get 
	\begin{equation*}
		\begin{split}
			&\Theta^k(\pmb{u}, \pmb{t},\pmb{z},s) \\
			& = \inf\Big\{ r \geq 0 \colon \int_0^r\la V_t ( \mathfrak{V}^k(\overset{\leftarrow}{\pmb{v}}, \overset{\leftarrow}{\pmb{w}},\pmb{\tau}  ) + w_{k+1}\delta_{v_{k+1}}),1\ra dt >s + \pmb{\Delta}\theta_r^k(\overset{\leftarrow}{\pmb{y }}, \overset{\leftarrow}{\pmb{x }}) \Big \},
		\end{split}	
	\end{equation*} 
	where
	\begin{equation*}
		\overset{\leftarrow}{\pmb{x}} = (\pmb{u}, \pmb{t},\pmb{z} ) , \quad \overset{\leftarrow}{\pmb{y}} = (\pmb{v}, \pmb{\tau},\pmb{w})
	\end{equation*}
	and
	\begin{equation*}
		\pmb{\Delta}\theta_r^k(\overset{\leftarrow}{\pmb{y}}, \overset{\leftarrow}{\pmb{x }}) = \int_0^r\la V_t ( \mathfrak{V}^k(\overset{\leftarrow}{\pmb{v}}, \overset{\leftarrow}{\pmb{w}},\pmb{\tau}  ) + w_{k+1}\delta_{v_{k+1}})-V_t ( \mathfrak{V}^k(\overset{\leftarrow}{\pmb{u}}, \overset{\leftarrow}{\pmb{z}},\pmb{t}  ) + z_{k+1}\delta_{u_{k+1}}),1\ra dt.
	\end{equation*}
	Now, we denote 
	\begin{equation*}
		\eta(\delta) = \sup_{\pmb{y}\in B_\delta(\pmb{x})}\sup_{r \in [0,T]} |\pmb{\Delta}\theta_r^k(\overset{\leftarrow}{\pmb{y}})|.
	\end{equation*}
	It follows that 
	\begin{equation*}
		\Theta^k(\overset{\leftarrow}{\pmb{y}}, s-\eta) \leq 	\Theta^k(\overset{\leftarrow}{\pmb{x}}, s) \leq 	\Theta^k(\overset{\leftarrow}{\pmb{y}}, s+\eta),
	\end{equation*}
	and thus 
	\begin{equation*}
		|\Theta^k(\overset{\leftarrow}{\pmb{x}}, s) - \Theta^k(\overset{\leftarrow}{\pmb{y}}, s) |\leq  |\Theta^k(\overset{\leftarrow }{\pmb{y}},s+\eta) - \Theta^k(\overset{\leftarrow }{\pmb{y}},s) | +  |\Theta^k(\overset{\leftarrow }{\pmb{y}},s-\eta) - \Theta^k(\overset{\leftarrow }{\pmb{y}},s) |.
	\end{equation*}
	It thus remains to prove that $\eta(\delta) \rightarrow 0 $ as $\delta \rightarrow 0$ and 
	\begin{equation}\label{eq:conv_theta_s}
		\sup_{\pmb{y}\in B_\delta(\pmb{x})} |\Theta^k(\overset{\leftarrow }{\pmb{y}},s+\tilde{\eta}) - \Theta^k(\overset{\leftarrow }{\pmb{y}},s) | \rightarrow 0,
	\end{equation}
	as $\tilde{\eta} \rightarrow 0$.
	First, denote 
	\begin{equation*}
		D_R = \bigcup_{i = 1}^{k+1}B_R(u_i), \quad R>0.
	\end{equation*}
	We get for $\bar{t} \in (0,T)$ and $R  > 0$  that
	\begin{equation}\label{eq:eta_est}
		\begin{split}
		  \eta(\delta)  &\leq \sup_{\pmb{y}\in B_\delta(\pmb{x})}  \int_0^T \la|V_t ( \mathfrak{V}^k(\overset{\leftarrow}{\pmb{v}}, \overset{\leftarrow}{\pmb{w}},\pmb{\tau}  ) + w_{k+1}\delta_{v_{k+1}})-V_t ( \mathfrak{V}^k(\overset{\leftarrow}{\pmb{u}}, \overset{\leftarrow}{\pmb{z}},\pmb{t}  ) + z_{k+1}\delta_{u_{k+1}}) |,1\ra dt\\
			& \leq 2 \sup_{\pmb{y}\in B_\delta(\pmb{x})} \int_0^{\bar{t}} \la V_t ( \mathfrak{V}^k(\overset{\leftarrow}{\pmb{v}}, \overset{\leftarrow}{\pmb{w}},\pmb{\tau}  ) + w_{k+1}\delta_{v_{k+1}}), 1\ra dt\\
			& \qquad +  \sup_{\pmb{y}\in B_\delta(\pmb{x})} \int_{\bar{t}}^T \la |V_t ( \mathfrak{V}^k(\overset{\leftarrow}{\pmb{v}}, \overset{\leftarrow}{\pmb{w}},\pmb{\tau}  ) + w_{k+1}\delta_{v_{k+1}})-V_t ( \mathfrak{V}^k(\overset{\leftarrow}{\pmb{u}}, \overset{\leftarrow}{\pmb{z}},\pmb{t}  ) + z_{k+1}\delta_{u_{k+1}}) | , \mathbbm{1}_{D_R} \ra dt \\
			& \qquad + 2(k+1)\int_{\bar{t}}^T \la \mathcal{W}_t(0,0) , \mathbbm{1}_{B_{R^c}(0)}\ra dt.  
		\end{split}
	\end{equation}
	The second term on the right hand side of \eqref{eq:eta_est} converges to zero as $\delta \rightarrow 0$ by Proposition~\ref{prp:initial_trace_prop}  and
	for $\tilde{\epsilon} >0$ arbitrary we can choose $\bar{t}\in (0,T)$ small enough and 
	$R >0$ large enough such that 
	\begin{equation*}
		\sup_{\pmb{y}\in B_\delta(\pmb{x})} \int_0^{\bar{t}} \la V_t ( \mathfrak{V}^k(\overset{\leftarrow}{\pmb{v}}, \overset{\leftarrow}{\pmb{w}},\pmb{\tau}  ) + w_{k+1}\delta_{v_{k+1}}), 1\ra dt + (k+1)\int_{\bar{t}}^T \la \mathcal{W}_t(0,0) , \mathbbm{1}_{B_{R^c}(0)}\ra dt \quad <\tilde{\epsilon}.
	\end{equation*}
	Thus, it remains to prove \eqref{eq:conv_theta_s}. To this end it is enough to prove 
	\begin{equation*}
		\begin{split}
			r_{\Delta}(\tilde{\eta})
			:=\sup &\Big\{ |r_2-r_1|, 0 \leq r_1\leq r_2 \leq T\colon \\
			&\qquad  \exists  \pmb{y} \in B_1(\pmb{x}) \text{ s.t. } \int_{r_1}^{r_2}\la V_t ( \mathfrak{V}^k(\overset{\leftarrow}{\pmb{v}}, \overset{\leftarrow}{\pmb{w}},\pmb{\tau}  ) + w_{k+1}\delta_{v_{k+1}}),1\ra dt \leq \tilde{\eta}  \Big\}
		\end{split}
	\end{equation*}
	converges to zero as $\tilde{\eta} \rightarrow 0$.
	Note that by using Lemma~\ref{lem:pde_est} and \eqref{eq:vmon} we have
	\begin{equation*}
		\begin{split}
			\int_{r_1}^{r_2} \la  V_t ( \mathfrak{V}^k(\overset{\leftarrow}{\pmb{v}}, \overset{\leftarrow}{\pmb{w}},\pmb{\tau}  ) + w_{k+1}\delta_{v_{k+1}}), 1\ra dt \geq \int_{r_1}^{r_2}\la V_{t+\bar{t}}(\mu), 1\ra dt\geq |r_2-r_1| \la V_{T+\bar{t}}(\mu ), 1\ra,
		\end{split}
	\end{equation*}
	where $\bar{t}$ is given by \eqref{eq:bart} and $T$ is given by \eqref{eq:supPhi}. 
	Thus, $r_{\Delta}(\tilde{\eta}) \rightarrow 0$ as $\tilde{\eta} \rightarrow 0$. This finishes the proof.
\end{proof}

In the following lemma we establish the weak convergence of the random vectors
\begin{align}
	&\pmb{U}_{k}^{n,m} = (x_0,U_{n,1}^m, \ldots U_{n,k}^m),\\
	&\pmb{Z}_k^{n,m} = (m,Z_{n,1}^m, \ldots Z_{n,k}^m),\\
	&\pmb{T}_k^{n,m} = (T_{n,1}^m, \ldots T_{n,k}^m),\\
	&\pmb{M}_k^{n,m} = (M_{n,1}^m,\ldots M_{n,k}^m),
\end{align}
for $k \geq 1$, which will allow us to prove the convergence towards the dual process.

\begin{lemma}\label{lem:weak_tuzm}
	For every $k \in \N$ it follows that
	\begin{equation}	
		(\pmb{U}_{k}^{n,m}, \pmb{Z}_{k}^{n,m}, \pmb{T}_{k}^{n,m}, \pmb{M}_{k}^{n,m}) \Rightarrow 	(\pmb{U}_{k}^{\infty,m}, \pmb{Z}_{k}^{\infty,m}, \pmb{T}_{k}^{\infty,m}, \pmb{M}_{k}^{\infty,m}), \quad \forall m \in \R_+,
	\end{equation}
	as $n \rightarrow \infty$ and
	\begin{equation}\label{eq:conv_mm}
		(\pmb{U}_{k}^{m,m}, \pmb{Z}_{k}^{m,m}, \pmb{T}_{k}^{m,m}, \pmb{M}_{k}^{m,m}) \Rightarrow 	(\pmb{U}_{k}^{\infty,\infty}, \pmb{Z}_{k}^{\infty,\infty}, \pmb{T}_{k}^{\infty,\infty}, \pmb{M}_{k}^{\infty,\infty}),
	\end{equation}
	as $m \rightarrow \infty$ along integers
	and 
	\begin{equation}
		(\pmb{U}_{k}^{\infty,m}, \pmb{Z}_{k}^{\infty,m}, \pmb{T}_{k}^{\infty,m}, \pmb{M}_{k}^{\infty,m}) \Rightarrow 	(\pmb{U}_{k}^{\infty,\infty}, \pmb{Z}_{k}^{\infty,\infty}, \pmb{T}_{k}^{\infty,\infty}, \pmb{M}_{k}^{\infty,\infty}),
	\end{equation}
	as $m \rightarrow \infty$ along integers,
	with respect to the canonical metric of $\R^{k+1}\times\overline{\R}_{+}^{k+1} \times \R_+^k\times\{1,2\}^k$.
\end{lemma}
\begin{proof}
	We only prove \eqref{eq:conv_mm}. The other convergences follow in a similar manner. 
	Furthermore, the first components of $\pmb{U}_k^{m,m}$ and $\pmb{Z}_{k}^{m,m}$ are convergent and are deterministic.\\
	We denote 
	\begin{equation*}
		(P_{m,k}^1, P_{m,k}^2,P_{m,k}^3,P_{m,k}^4) = ( Z_{m,k}^m, M_{m,k}^m, U_{m,k}^m, T_{m,k}^m)
	\end{equation*}
	and 
	\begin{equation*}
		D_1\times D_2\times D_3 \times D_4 = \overline{\R}_{+}\times\{ 1,2\}\times \R\times \R_+.
	\end{equation*} 
	It is enough to prove
	\begin{equation}\label{eq:weak_conv_cond}
		\lim_{m \rightarrow \infty }\E\Big[\prod_{\ell = 1}^{k}\prod_{i = 1}^4f^i_\ell(P_{m,\ell}^i)\Big] = \E\Big[\prod_{\ell = 1}^{k}\prod_{i = 1}^4f_\ell^i(P_{\infty,\ell}^i)\Big]
	\end{equation}
	for all continuous, bounded functions $f_{k}^i$ on $D_{i}$ for $i \in \{1,\ldots,4\}$ and $k \geq 1$.
	We prove the claim  via induction. 
	First set $k = 1$ and let $f^i$ for $i \in \{1, \ldots 4\}$ be arbitrary  continuous bounded functions on
	$D_i$ for $i \in \{1,\ldots 4\}$.
	We need to prove
	\begin{equation}\label{eq:weak_conv1}
		\lim_{m \rightarrow \infty }\E\Big[\prod_{i = 1}^4f^i(P_{m,1}^i)\Big] = \E\Big[\prod_{i = 1}^4f^i(P_{\infty,1}^i)\Big].
	\end{equation} 
	We get 
	\begin{equation}\label{eq:induc_weak_1}
		\E\Big[\prod_{i = 1}^4f^i(P_{m,1}^i)\Big]  = \E\Big[\E\Big[\prod_{i = 1}^4f^i(P_{m,1}^i)| T_{m,1}^m\Big]\Big] = \E\Big[\E\Big[\prod_{i = 1}^3f^i(P_{m,1}^i)|T_{m,1}^m\Big] f^4(T_{m,1}^m)\Big].
	\end{equation}
	Now, by using that $(Z_{m,1}^m,M_{m,1}^m)$ and $U_{m,1}^m$  are conditionally independent given $T_{m,1}^m$, we arrive at 
	\begin{equation}\label{eq:weak_conv1.2}
		\E\Big[\prod_{i = 1}^3f^i(P_{m,1}^i)|T_{m,1}^m\Big] = \E[f^1(Z_{m,1}^m) f^2(M_{m,1}^m)|T_{m,1}^m]\E[f^3(U_{m,1}^m)|T_{m,1}^m].
	\end{equation}
	Furthermore, note that $(Z_{m,1}^m, M_{m,1}^m)$ is independent of $T_{m,1}^m$ and thus
	\begin{equation}\label{eq:weak_conv1.3}
		\begin{split}
			&\E[f^1(Z_{m,1}^m) f^2(M_{m,1}^m)|T_{m,1}^m] = \E[f^1(Z_{m,1}^m) f^2(M_{m,1}^m)]\\
			&=\E[\E[f^1(Z_{m,1}^m)|M_{m,1}^m] f^2(M_{m,1}^m)] \\
			&= \E[f^1(\bar{Z}_m(1))]\E[\mathbbm{1}_{M_{m,1}^m = 1}f^2(1)] +\E[f^1(\bar{Z}_m(2))]\E[\mathbbm{1}_{M_{m,1}^m = 2}f^2(2)],
		\end{split}
	\end{equation}
	where $\bar{Z}_m(j)$ is a random variable distributed according to  
	\begin{equation*}
		\P(\bar{Z}_m(j) \geq b) = \frac{\nu^j_m([b,\infty))}{\nu^j_m((0,\infty))},
	\end{equation*}
	for $j \in \{1,2\}$ and $b>0$.
	Furthermore, we get 
	\begin{equation}\label{eq:weak_conv_1.4}
		\E[f^3(U_{m,1}^m)|T_{m,1}^m] = \frac{1}{\int_\R \mathfrak{V}^1(x_0, m ,T_{m,1}^m)(x)dx}\int_\R f^3(x) \mathfrak{V}^1 (x_0, m, T_{m,1}^m) (x)dx.
	\end{equation}
	Together, \eqref{eq:weak_conv1} , \eqref{eq:weak_conv1.2}, \eqref{eq:weak_conv1.3} and \eqref{eq:weak_conv_1.4} yield 
	\begin{equation}
		\begin{split}
			&\E\Big[\prod_{i = 1}^4f^i(P_{m,1}^i)\Big]\\
			&  = \sum_{j = 1}^2 \E[f^1(\bar{Z}_m(j))]  \E[\mathbbm{1}_{M_{m,1}^m = j}f^2(j)] \\
			& \qquad \times\E\Big[\frac{1}{\int_\R \mathfrak{V}^1 (x_0, m ,T_{m,1}^m)(x)dx}\int_\R f^3(x) \mathfrak{V} ^1(x_0, m ,T_{m,1}^m)(x)dxf^4(T_{m,1}^m)\Big].
		\end{split}
	\end{equation}
	Clearly $\bar{Z}_m(j)$ and $M_{m,1}^m$ converge in distribution to $\bar{Z}_\infty(j)$ and $M_{\infty,1}^\infty$, respectively. Hereby, $\bar{Z}_\infty(j)$ is distributed according to  
	\begin{equation*}
		\P(\bar{Z}_\infty(j) \in (0,b]) = \frac{\nu^{j}((0,b]) + d_j \mathbbm{1}_{b = \infty}}{\nu^{j}((0,\infty)) + d_j }, \quad j = 1,2 , b\in \bar{\R}^+.
	\end{equation*}
	Thus, it only remains to 
	show the convergence of 
	\begin{equation*}
		\frac{1}{\int_\R \mathfrak{V}^1(x_0, m ,T_{m,1}^m)dx}\int_\R f^3(x) \mathfrak{V}^1 (x_0, m, T_{m,1}^m) dx.
	\end{equation*}
	To this end, note that $S_{m,1}^m$ has the same distribution for all $m$ and that $Z_{m,0} = m$ converges 
	weakly to $Z_{\infty,0} = \infty$ with respect to the metric $d$ as $m \rightarrow \infty$. Along similar lines as in  the 
	proof of Lemma~\ref{lem:theta_cont} we get that the map
	\begin{equation*}
		(z, s) \mapsto \inf\{ r\geq 0 \colon \int_0^r \la V_t(\mu + z\delta_x), 1\ra dt >s\} 
	\end{equation*} 
	is continuous and bounded with respect to the canonical metric on $\overline{\R}_{+} \times \R_+$. Note that 
	\begin{equation*}
		T^m_{m,1} =  \inf\{ r\geq 0 \colon \int_0^r \la V_t(\mu + m\delta_x), 1\ra dt >S_{m,1}^m/\nu^m((0,\infty])\} . 
	\end{equation*}
	Thus, we get that $T_{m,1}^m$ converges weakly to $T_{\infty,1}^\infty$. 
	Furthermore, the map
	\begin{equation*}
		(t,z) \mapsto \frac{1}{\int_\R \mathfrak{V}^1(x_0, z ,t)(x)dx}\int_\R f^3(x) \mathfrak{V}^1 (x_0, z, t)(x) dx
	\end{equation*}
	is continuous and bounded for all $(t, z) \in (0,\infty)\times \overline{\R}_+$ . Since $\P(T_{m,1}^m = 0) =0$ the convergence 
	\begin{equation*}
	\begin{split}
		&\lim_{m\rightarrow \infty }\E\Big[\frac{1}{\int_\R \mathfrak{V}^1 (x_0, m ,T_{m,1}^m)(x)dx}\int_\R f^3(x) \mathfrak{V} ^1(x_0, m ,T_{m,1}^m)(x)dxf^4(T_{m,1}^m)\Big]\\
		& = \E\Big[\frac{1}{\int_\R \mathfrak{V}^1 (x_0, \infty ,T_{\infty,1}^\infty)(x)dx}\int_\R f^3(x) \mathfrak{V} ^1(x_0, \infty,T_{\infty,1}^\infty)(x)dxf^4(T_{\infty,1}^\infty)\Big],
	\end{split}
	\end{equation*}
	follows.  Along similar lines to the above we get that 
	\begin{equation*}
		\begin{split}
			&  \E\Big[\prod_{i = 1}^4f^i(P_{\infty,1}^i)\Big]\\
			&  = \sum_{j = 1}^2\E[f^1(\bar{Z}_{\infty}(1))]\E[\mathbbm{1}_{M_{\infty,1}^\infty =j}f^2(M_{\infty,1}^\infty)]\\
			& \qquad \times\E\Big[ \frac{1}{\int_\R \mathfrak{V}^1(x_0, \infty ,T_{\infty,1}^\infty)(x)dx}\int_\R f^3(x) \mathfrak{V}^1 (x_0, \infty, T_{\infty,1}^\infty)(x) dxf^4(T_{\infty,1}^\infty)\Big]
		\end{split}
	\end{equation*}
	and thus \eqref{eq:weak_conv1} holds.\\
	Now we assume that \eqref{eq:weak_conv_cond} holds for some $k\geq 1$ and prove that it also holds for $k +1$. 
	To this end let $f^i_\ell$ be bounded, continuous functions on $D^i$, for $i \in \{1,\ldots 4\}$ and $\ell \in \{1,\ldots k+1\}$. 
	Note that $(Z_{m,k+1}^m, M_{m,k+1}^m)$ is independent of $(\pmb{U}_{k+1}^{m,m}, \pmb{Z}_{k}^{m,m}, \pmb{T}_{k+1}^{m,m}, \pmb{M}_{k}^{m,m})$ . Thus it is enough to show 
	\begin{equation*}
		\begin{split}
			&\lim_{m \rightarrow \infty }\E\Big[\prod_{\ell = 1}^{k}\prod_{i = 1}^4f^i_\ell(P_{m,\ell}^i)f_{k+1}^3(U_{k+1,m}^m)f_{k+1}^4(T_{m,k+1}^m)\Big] \\
			&\qquad = \E\Big[\prod_{\ell = 1}^{k}\prod_{i = 1}^4f_\ell^i(P_{\infty,\ell}^i)f_{k+1}^3(U_{\infty,k+1}^\infty)f_{k+1}^4(T_{\infty,k+1}^\infty)\Big].
		\end{split}
	\end{equation*}
	First, consider 
	\begin{equation*}
		\E\Big[\prod_{\ell = 1}^{k}\prod_{i = 1}^3f^i_\ell(P_{m,\ell}^i)f_{k+1}^4(T_{m,k+1}^m)\Big],
	\end{equation*}
	and note that by \eqref{eq:phikexplain} we have
	\begin{equation*}
		f_{k+1}^4(T_{m,k+1}^m) = f_{k+1}^4(\Theta^k ( \pmb{U}_k^{m,m}, \pmb{T}_{k}^{m,m}, \pmb{Z}_{k}^{m,m}, S_{m,k+1}^{m}/\nu_m((0,\infty]))).
	\end{equation*}
	Thus, by Lemma~\ref{lem:theta_cont} and the continuous mapping theorem,
	 it follows that
	\begin{equation*}
		(\pmb{U}_{k}^{m,m}, \pmb{T}_{k+1}^{m,m}, \pmb{Z}_{k+1}^{m,m}, \pmb{M}_{k+1}^{m,m}) \Rightarrow (\pmb{U}_{k}^{ \infty,\infty}, \pmb{T}_{k+1}^{\infty,\infty}, \pmb{Z}_{k+1}^{\infty,\infty}, \pmb{M}_{k+1}^{\infty,\infty}),
	\end{equation*}
	as $m\rightarrow \infty$. Again, since  
	\begin{equation*}
		(\pmb{u}, \pmb{t}, \pmb{z}, t) \mapsto \frac{1}{\int_\R V_{t}( \mathfrak{V}^k(\pmb{u}, \pmb{z}, \pmb{t}))(x) dx} \int_\R f^3_{k+1} (x)V_{t}( \mathfrak{V}^k(\pmb{u}, \pmb{z}, \pmb{t}))(x) dx
	\end{equation*}
	is bounded and continuous on $\R^k \times (0,\infty)^k\times\bar{\R}_+^k \times(0,\infty)$, the convergence 
	\begin{equation*}
		(\pmb{U}_{k+1}^{m,m}, \pmb{T}_{k+1}^{m,m}, \pmb{Z}_{k+1}^{m,m}, \pmb{M}_{k+1}^{m,m}) \Rightarrow (\pmb{U}_{k+1}^{ \infty,\infty}, \pmb{T}_{k+1}^{\infty,\infty}, \pmb{Z}_{k+1}^{\infty,\infty}, \pmb{M}_{k+1}^{\infty,\infty}), \quad \text{ as } m\rightarrow \infty,
	\end{equation*}
	follows. This finishes the proof.
\end{proof}

\begin{lemma}\label{lem:conv_drift_exp}
	Let $a\in \R$ and let $X_0\in \mcC_{tem}^+$. Then there  exists $T_0>0$, independent of $X_0$, such that the following convergences hold:
	\begin{enumerate}[label=(\roman*)]
		\item \label{eq:conv1}
		\begin{equation*}
			\begin{split}
				&\lim_{m\rightarrow \infty}\E\Big[ (-1)^{J_t^{m,m}} \exp\Big(-\la Y^{m,m}_{t}, X_0\ra-a\int_0^t\la Y^{m,m}_s,1\ra ds\Big)\Big]\\
				& = \E\Big[ (-1)^{J_t^{\infty,\infty}} \exp\Big(-\la Y^{\infty,\infty}_{t}, X_0\ra-a\int_0^t\la Y^{\infty,\infty}_s,1\ra ds\Big)\Big],
			\end{split}
		\end{equation*}
		\item \label{eq:conv2}
		\begin{equation*}
			\begin{split}
				&	\lim_{n \rightarrow \infty}\E\Big[ (-1)^{J_t^{n,m}} \exp\Big(-\la Y^{n,m}_{t}, X_0\ra-a\int_0^t\la Y^{n,m}_s,1\ra ds\Big)\Big] \\
				&= \E\Big[ (-1)^{J_t^{\infty,m}} \exp\Big(-\la Y^{\infty,m}_{t}, X_0\ra-a\int_0^t\la Y^{\infty,m}_s,1\ra ds\Big)\Big],
			\end{split}
		\end{equation*}
		\item\label{eq:conv3}
		\begin{equation*}
			\begin{split}
				&	\lim_{m \rightarrow \infty}\E\Big[ (-1)^{J_t^{\infty,m}} \exp\Big(-\la Y^{\infty,m}_{t}, X_0\ra-a\int_0^t\la Y^{\infty,m}_s,1\ra ds\Big)\Big] \\
				&= \E\Big[ (-1)^{J_t^{\infty,\infty}} \exp\Big(-\la Y^{\infty,\infty}_{t}, X_0\ra-a\int_0^t\la Y^{\infty,\infty}_s,1\ra ds\Big)\Big],
			\end{split}
		\end{equation*}
	for all $t \in [0,T_0]$.
	\end{enumerate}
\end{lemma}
\begin{proof}
	We only show \ref{eq:conv1}, since \ref{eq:conv2} and \ref{eq:conv3} follow along similar lines. 
	Let $T_0>0 $ be such that \eqref{eq:est_lem65} is satisfied with $\gamma = 2a$. 
	Fix arbitrary $t\in [0,T_0]$ and let $\epsilon > 0$. 
By \eqref{eq:ineq_times}, there exists 
	a sequence of random variables $(\bar{R}_i)_{i \geq 1}$ such that 
	\begin{equation*}
		\P(\bar{T}_{m,i}^m \leq t) \leq \P(\bar{R}_i \leq t), \quad \forall m,i\geq 1, \quad t\geq 0.
	\end{equation*}
	Furthermore, by using Lemma~\ref{lem:number_jumps}, it is not hard to see that $\lim_{i \rightarrow \infty} \bar{R}_i = \infty$ almost surely. 
	Thus, there exists $k \in \N$ such that 
	\begin{equation}\label{eq:cond_T}
		\sup_{m\in \N}\P(\bar{T}_{m,k}^m\leq t) <\epsilon/2
	\end{equation}
	and 
	\begin{equation*}
		\P(\bar{T}_{\infty,k}^\infty\leq t) <\epsilon/2.
	\end{equation*}
	Thus we arrive at 
	\begin{equation*}
		\begin{split}
			&\Big| \E\Big[ (-1)^{J_t^{m,m}} \exp\Big(-\la Y^{m,m}_{t}, X_0\ra-a\int_0^t\la Y^{m,m}_s,1\ra ds\Big)\Big]\\
			& \quad -\E\Big[ (-1)^{J_t^{\infty,\infty}} \exp\Big(-\la Y^{\infty,\infty}_{t}, X_0\ra-a\int_0^t\la Y^{\infty,\infty}_s,1\ra ds\Big)\Big]\Big| \\
			& \leq C\epsilon + \sum_{i = 1}^k \Big| \E\Big[ (-1)^{J_t^{m,m}}\exp\Big(-\la Y^{m,m}_{t}, X_0\ra-a\int_0^t\la Y^{m,m}_s,1\ra ds\ra\Big)\mathbbm{1}_{t \in [\bar{T}_{m,i-1}^m,\bar{T}_{m,i}^m)}\Big] \\
			& \qquad 
			- \E\Big[ (-1)^{J_t^{\infty,\infty}}\exp\Big(-\la Y^{\infty,\infty}_{t}, X_0\ra-a\int_0^t\la Y^{\infty,\infty}_s,1\ra ds\Big) \mathbbm{1}_{t \in [\bar{T}_{\infty,i-1}^\infty,\bar{T}_{\infty,i}^\infty)}\Big] \Big|,
		\end{split}
	\end{equation*}
	where we used that 
	\begin{equation*}
		 \sup_{m}\E\Big[ \exp\Big(-2a \int_0^t \la Y_s^{m,m},1\ra ds\Big)\Big] < \infty. 
	\end{equation*}
	Note that by \eqref{eq:frakV_k_identity} we get for $i \in \{2,\ldots, k\}$   that 
	\begin{equation}\label{eq:ymm_identity_1}
		Y_t^{m,m}(x) = 	\mathfrak{V}^{i}(\pmb{U}_{i-1}^{m,m}, \pmb{Z}_{i-1}^{m,m} , (\pmb{T}_{i-1}^{m,m}, t-\bar{T}_{m,i-1}^{m}))(x), \quad x\in \R,
	\end{equation}
	and for $t  \in (\bar{T}_{m,i-1}^m,\bar{T}_{m,i}^m) $ and 
	\begin{equation}\label{eq:yii_indentity_1}
		Y_t^{\infty,\infty}(x) = 	\mathfrak{V}^{i}(\pmb{U}_{i-1}^{\infty,\infty}, \pmb{Z}_{i-1}^{\infty,\infty} , (\pmb{T}_{i-1}^{\infty,\infty}, t-\bar{T}_{\infty,i-1}^{\infty}))(x), \quad x\in \R,
	\end{equation}
	and for $ t  \in (\bar{T}_{\infty,i-1}^\infty,\bar{T}_{\infty,i}^\infty) $. 
    Furthermore, note that 
	\begin{equation}\label{eq:J_mm_identity}
		J_t^{m,m} = \sum_{\ell= 1}^{i-1} \mathbbm{1}_{M_{m,\ell}^m= 2}, \quad t \in [\bar{T}_{i-1,m}^m, \bar{T}_{i,m}^m).
	\end{equation}
	Now, it is sufficient to show that
	\begin{equation}\label{eq:weak_con_summand}
		\begin{split}
			&\Big| \E\Big[ (-1)^{J_t^{m,m}}\exp\Big(-\la Y^{m,m}_{t}, X_0\ra-a\int_0^t\la Y^{m,m}_s,1\ra ds\ra\Big)\mathbbm{1}_{t \in [\bar{T}_{m,i-1}^m,\bar{T}_{m,i}^m)}\Big] \\
			& \qquad 
			- \E\Big[ (-1)^{J_t^{\infty,\infty}}\exp\Big(-\la Y^{\infty,\infty}_{t}, X_0\ra-a\int_0^t\la Y^{\infty,\infty}_s,1\ra ds\Big) \mathbbm{1}_{t \in [\bar{T}_{\infty,i-1}^\infty,\bar{T}_{\infty,i}^\infty)}\Big] \Big|
		\end{split}
	\end{equation}
	converges to zero as $m \rightarrow \infty$. Let $\epsilon >0$ be arbitrary and choose $K >0$ sufficiently large. 
	We get that
	\begin{equation}\label{eq:expKbd1}
		\begin{split}
			&\sup_m\E\Big[\Big|\exp\Big(-a\int_0^t\la Y^{m,m}_s,1\ra ds\ra\Big) -\exp\Big(-a\int_0^t\la Y^{m,m}_s,1\ra ds\ra\Big)\wedge K \Big|\Big]\\
			& = \sup_m\E\Big[\exp\Big(-a\int_0^t\la Y^{m,m}_s,1\ra ds\ra\Big)\mathbbm{1}_{\exp\big(-a\int_0^t\la Y^{m,m}_s,1\ra ds\ra\big) >K}\Big]\\
			& \leq \frac{C}{K}<\epsilon/2,
		\end{split}
	\end{equation}
	for $K >0$ sufficiently large by using Lemma~\ref{lem:exponential_moment}. Similarly, we can ensure that 
	\begin{equation}\label{eq:expKbd2}
		\E\Big[\Big|\exp\Big(-a\int_0^t\la Y^{\infty,\infty}_s,1\ra ds\ra\Big) -\exp\Big(-a\int_0^t\la Y^{\infty,\infty}_s,1\ra ds\ra\Big)\wedge K \Big|\Big] < \epsilon/2.
	\end{equation}
	  Due to Lemma~\ref{lem:weak_tuzm} and \eqref{eq:ymm_identity_1}, \eqref{eq:yii_indentity_1} and \eqref{eq:J_mm_identity} we also get that
	\begin{equation*}
		\begin{split}
			&\Big| \E\Big[ (-1)^{J_t^{m,m}}\Big(\exp\Big(-\la Y^{m,m}_{t}, X_0\ra-a\int_0^t\la Y^{m,m}_s,1\ra ds\ra\Big)\wedge K\Big) \mathbbm{1}_{t \in [\bar{T}_{m,i-1}^m,\bar{T}_{m,i}^m)}\Big] \\
			& \quad 
			- \E\Big[ (-1)^{J_t^{\infty,\infty}}\Big(\exp\Big(-\la Y^{\infty,\infty}_{t}, X_0\ra-a\int_0^t\la Y^{\infty,\infty}_s,1\ra ds\Big)\wedge K\Big) \mathbbm{1}_{t \in [\bar{T}_{\infty,i-1}^\infty,\bar{T}_{\infty,i}^\infty)}\Big] \Big|,
		\end{split}
	\end{equation*}
	converges to zero as $m \rightarrow \infty$. 
	 This together with \eqref{eq:expKbd1} and \eqref{eq:expKbd2} demonstrates the convergence of \eqref{eq:weak_con_summand} to zero as $m \rightarrow \infty$ and thus finishes the proof.  
\end{proof}

\begin{lemma}\label{lem:conv_yh}
Let $T_0 > 0$ be such that \eqref{eq:est_lem65} in Lemma~\ref{lem:exponential_moment} holds for $\gamma = 2a$ and $K= \nu((0,\infty])$. Furthermore, let $(X_{t})_{t\geq 0}$ be a solution to \eqref{eq:SPDE1} with initial condition $X_0$ satisfying \eqref{eq:cond_x_0}. Furthermore, let $(g_n)_{n\geq 1}$ be a sequence of measurable bounded functions on $\R$, that converge boundedly pointwise to a bounded, measurable function $g$ on $\R$.
 Then it follows for all $T \leq T_0$ that
	\begin{equation*}
		\lim_{n \rightarrow \infty }\E\Big[\Big(\int_0^T \la Y_{s-}^{n,m}, |g_n(X_{T-s})- g(X_{T-s})| \ra ds\Big)^2\Big] = 0,
	\end{equation*}
for all $m\in \R_+$.
\end{lemma}
\begin{proof}
Fix $m\in \R_+$. To lighten notation, we suppress the superscript $m$ for the rest of the proof. Moreover, fix arbitrary $T\leq T_0$.
Note that by Lemma~\ref{lem:exponential_moment} with $K = \nu((0,\infty])$ we have
\begin{equation}\label{eq:est_lem610}
	\sup_{n}\E\Big[\Big(\int_0^T \la Y_{s-}^{n}, |g_n(X_{T-s})- g(X_{T-s})| \ra ds\Big)^4\Big] \leq 	c\sup_{n}\E\Big[\Big(\int_0^T \la Y_{s-}^{n}, 1\ra ds\Big)^4\Big] < \infty.
\end{equation}
Furthermore, we can use a similar argument to the proof of Lemma~\ref{lem:conv_drift_exp} to see that there exists $M\in \N$ such that
\begin{equation*}
	\P(\inf_n \bar{T}_{n,M} \leq T ) <\epsilon, 
\end{equation*}
and by  Lemma~\ref{lem:weak_tuzm} it follows that 
\begin{equation}\label{eq:weak_conv_3}
\begin{split}
	&(U_{n,1}, \ldots,U_{n,M}, T_{n,1}, \ldots, T_{n,M} , Z_{n,1},\ldots, Z_{n,M}) \\
&\qquad \Rightarrow (U_{\infty,1}, \ldots,U_{\infty,M}, T_{\infty,1},\ldots ,T_{\infty,M}, Z_{\infty,1}, \ldots, Z_{\infty,M})
\end{split}
\end{equation}
  and  thus by the Skorohod representation theorem there exists a probability space such that 
the convergence \eqref{eq:weak_conv_3} holds almost surely.
We now show that
\begin{equation*}
	\lim_{n \rightarrow \infty}\int_0^T \la Y_{s-}^n, |g_n(X_{T-s})- g(X_{T-s})| \ra ds = 0
\end{equation*}
for almost all $\omega$ in the event 
\begin{equation*}
	A_\epsilon = \{ \inf_n \bar{T}_{n,M} >T\}.
\end{equation*}
To this end, note that  there exists  $D(\omega) >0$  for $\omega\in A_\epsilon$ such that almost surely on $A_\epsilon$,
\begin{equation*}
	\sup_{n} |U_{n,k}(\omega)| < D(\omega), \quad D(\omega) <\infty,
\end{equation*}
for all $k \leq M$, since $(U_{n,k})_{n\geq 1}$ converges almost surely for each fixed $k$.
Let $\tilde{\epsilon} > 0$ be arbitrary. Then, by using  Lemma~\ref{lem:est_pde} and Proposition~\ref{prp:initial_trace_prop} we get the following.
There exists $\bar{D}(\omega)> D(\omega)$ such that for $\omega \in A_\epsilon$ we have
\begin{equation*}
	\int_0^T\int_{B(0,\bar{D})^c}  Y^n_{s-}(x)dxds \leq M\int_0^T \int_{B(0,\bar{D})^c}  \big(\mathcal{W}(D,0)\big)(x) dxds \leq \tilde{\epsilon},
\end{equation*}
with $\bar{D}(\omega) <\infty$ a.e. on $A_\epsilon$. We continue to work with $\omega \in A_\epsilon$. There exists $\delta_1 > 0$ such that 
\begin{equation*}
	M\int_0^{\delta_1}\la\mathcal{W}_s(0,0),1\ra ds < \tilde{\epsilon}.
\end{equation*}
Finally, once again by Proposition~\ref{prp:very_sing_asymp}, there exists $K(\omega)> 0$ such that 
\begin{equation*}
	M\sup_{s\in (\delta_1,T), x\in B(0,\bar{D}+D) }\mathcal{W}_{s}(0,0)(x) \leq K<\infty,
\end{equation*}
almost surely on $A_\epsilon$.
Thus, we get 
\begin{equation*}
\begin{split}
&	\int_0^T \la Y_{s-}^n, |g_n(X_{T-s})- g(X_{T-s})| \ra ds \\
&\leq C\tilde{\epsilon} + \sum_{k = 1}^M \int_{\bar{T}_{k-1,n}}^{\bar{T}_{k,n}\wedge T} \int_{B(0,\bar{D})}Y_{s-}^n(x) |g_n(X_{T-s}(x)) - g(X_{T-s}(x))|dxds\\
& \leq C\tilde{\epsilon} + \int_0^{T} \int_{B(0,\bar{D})}  \big(V_s(\mu+m\delta_{x_0})\big)(x) |g_n(X_{T-s}(x))-g(X_{T-s}(x))|dxds \\
& \quad+  \sum_{k = 2}^M \int_{\bar{T}_{k-1,n}}^{\bar{T}_{k,n}\wedge T} \int_{B(0,\bar{D})}Y_{s-}^n(x) |g_n(X_{T-s}(x)) - g(X_{T-s}(x))|dxds,
\end{split}
\end{equation*}
where $C = \sup_{x\in \R}|g(x)| +\sup_{n\geq 1}\sup_{x\in \R}|g_n(x)| <\infty$ is independent of $n$.
Furthermore, for $\delta_1>0$ sufficiently small we get 
\begin{equation*}
	\begin{split}
		&\sum_{k = 2}^M \int_{\bar{T}_{k-1,n}}^{\bar{T}_{k,n}\wedge T} \int_{B(0,\bar{D})}Y_{s-}^n(x) |g_n(X_{T-s}(x)) - g(X_{T-s}(x))|dxds\\
		&\leq CM\int_0^{\delta_1} \la \mathcal{W}_s(0,0),1\ra ds + \sum_{k = 2}^M \int_{\bar{T}_{k-1,n+\delta_1}}^{\bar{T}_{k,n}\wedge T} \int_{B(0,\bar{D})}Y_{s-}^n(x) |g_n(X_{T-s}(x)) - g(X_{T-s}(x))|dxds \\
		& \leq C\tilde{\epsilon} + M \big(\sup_{s\in (\delta_1,T), x\in B(0,\bar{D}+D) }\mathcal{W}_{s}(0,0)(x) \big) \int_0^T\int_{B(0,\bar{D})} |g_n(X_{T-s}(x)) - g(X_{T-s}(x))|dxds \\
		&\leq C\tilde{\epsilon} + K\int_0^T\int_{B(0,\bar{D})} |g_n(X_{T-s}(x)) - g(X_{T-s}(x))|dxds. 
	\end{split}
\end{equation*}
By the dominated convergence theorem, we get that 
\begin{equation*}
	\int_0^T\int_{B(0,\bar{D})} |g_n(X_{T-s}(x)) - g(X_{T-s}(x))|dxds \rightarrow 0,
\end{equation*}
as $n\rightarrow \infty$.  This demonstrates that
\begin{equation*}
	\lim_{n \rightarrow \infty} \int_0^T \la Y_{s-}^n, |g_n(X_{T-s})- g(X_{T-s})| \ra ds \leq 2C\tilde{\epsilon},
\end{equation*}
almost surely on $A_\epsilon$. Since $\tilde{\epsilon}>0$ was arbitrary, this proves that 
\begin{equation*}
		\lim_{n \rightarrow \infty} \int_0^T \la Y_{s-}^n, |g_n(X_{T-s})- g(X_{T-s})| \ra ds = 0,
\end{equation*}
for almost all $\omega \in A_{\epsilon}$.
Thus, by uniform integrability (see \eqref{eq:est_lem610}), we get
\begin{equation*}
\begin{split}
	\lim_{n \rightarrow \infty} \E\Big[\Big(\int_0^T \la Y_{s-}^n, |g_n(X_{T-s})- g(X_{T-s})| \ra ds\Big)^2\mathbbm{1}_{A_\epsilon}] = 0 .
\end{split}
\end{equation*}
Finally, by Hölder's inequality and \eqref{eq:est_lem610}, we get 
\begin{equation*}
	 \E\Big[\Big(\int_0^T \la Y_{s-}^n, |g_n(X_{T-s})- g(X_{T-s})| \ra ds\Big)^2\mathbbm{1}_{(A_\epsilon)^c}] \leq   c\epsilon\sqrt{  \E \Big[\Big(\int_0^T \la Y_s^n,1\ra ds\Big)^4\Big]} \leq C\epsilon,
\end{equation*}
and since $\epsilon > 0$ was arbitrary, this finishes the proof. 
\end{proof}

\subsection{Proof of Proposition~\ref{prop:duality}}
\label{sec:appr_dual}
In this section, we prove Proposition~\ref{prop:duality}. 
Throughout this section, we assume that there exists a weak solution $(X_t)_{t \geq 0}$ to \eqref{eq:SPDE1} with $X_0\sim \eta \in \mathcal{P}_p(\mathcal{C}_{tem}^+)$ such that \eqref{eq:cond_x_0} is satisfied.  
When $h$ is continuous, the existence of a solution was already resolved in \cite[Theorem 1.1]{Shiga94}. The existence of solutions when $h$ is discontinuous is proven in Section~\ref{sec:weak_exist} and relies on the duality relation applied to continuous approximations of the drift and convergence results established below.
Furthermore, we denote by $\E^X$ the conditional expectation 
given the $\sigma$-algebra $\mathcal{F}^Y$ generated by $(Y,J)$ and note that $X$ is independent of $\mathcal{F}^Y$. Similarly, we denote by $E^Y$ the conditional expectation 
given the $\sigma$-algebra $\mathcal{F}^X$ generated by $X$.\\

In order to establish the duality relation, we consider approximations of the process $Y$ defined in \eqref{eq:def_procc} and first derive an approximate duality. Then we use the results from Section~\ref{sec:convergence_dual} to pass to the limit.
To this end, we denote by $Y^n$ the process constructed in \eqref{eq:def_procc} with $\nu$ replaced by
\begin{equation*}
	\nu_n = \nu^1_{|[0,n]} + \nu^2_{|[0,n]} +d_1\delta_n +d_2\delta_n
\end{equation*}
and $Y_0 = \mu+m\delta_{x_0}$, where $\mu\in \mathcal{M}_F^+$, $x_0 \in \R$ and $m\in \N$. Recall that we define $d_1,d_2\geq 0$ and $a\in\R$ as in \eqref{eq:def_d}.
Accordingly, we define
\begin{equation}\label{eq:h_nduality}
\begin{split}
	h_n(x)&= \int_0^n (1-\ee^{-\lambda x}) \nu^1(d\lambda) +\int_0^n (1+\ee^{-\lambda x}) \nu^2(d\lambda) \\
    & \qquad+ d_1(1-\ee^{-nx}) +d_2(1+\ee^{-nx}) + a, \qquad  \qquad x \in \R_+,
\end{split}
\end{equation}
for $n\in \N$, which converges pointwise to $h$ given by
\begin{equation}\label{eq:h_duality}
\begin{split}
	h(x) &= \int_0^\infty (1-\ee^{-\lambda x}) \nu^1(d\lambda) +\int_0^\infty (1+\ee^{-\lambda x}) \nu^2(d\lambda)\\
    & \qquad + d_1(\mathbbm{1}_{x>0}) +d_2(1+\mathbbm{1}_{x=0}) + a,\qquad \qquad x \in \R_+.
\end{split}
\end{equation}
\begin{lemma}\label{lem:x_bound}
Let $X$ be a solution to \eqref{eq:SPDE1} with initial value $X_0 \sim \eta \in \mathcal{P}_p(\mathcal{C}^+_{tem})$, $p>2$,  such that \eqref{eq:cond_x_0} holds.
Then it follows that 
	\begin{equation*}
		\sup_{s\in [0,t]}\sup_{x \in \R}\E\big[X_s(x)\big] < \infty,
	\end{equation*}
for all $t >0$.
\end{lemma}
\begin{proof}
	Note that, since $h$ is bounded, 
	\begin{equation*}
	\begin{split}
		&\E[X_t(x)] = \E[\big(S_tX_0\big)(x)] + \int_0^t \int_\R p_{t-s}\E[h(X_{s}(y))]dyds \\
		&\leq \E[\big(S_tX_0\big)(x)] + C\int_0^t\int_\R p_{t-s}(x-y)dyds\\
		& = \E[\big(S_tX_0\big)(x)] + Ct,
	\end{split}
	\end{equation*}
for all $x\in \R$ and $t \geq 0$. This implies the result by standard estimates.
\end{proof}

\begin{lemma}\label{lem:ito}
	Let $X$ be a solution to \eqref{eq:SPDE1} with initial value $X_0 \sim \eta \in \mathcal{P}_p(\mathcal{C}^+_{tem})$, $p>2$, such that \eqref{eq:cond_x_0} holds.
	Let $(v_s)_{s\geq 0}$ be $\R$-valued, absolutely continuous, monotonically increasing or decreasing, locally bounded and such that $v_0 = 0$ and let $T >0$. Then it follows that 
	\begin{equation*}
		\begin{split}
			&\E^X\Big[\exp(-\la V_t(\mu), X_{T-t}\ra + v_t)\Big] = \E^X\Big[\exp(-\la V_0(\mu),X_T\ra)\Big]\\
			& \qquad  +\E^X\Big[ \int_0^t \exp(-\la V_s(\mu), X_{T-s}\ra+ v_s) \big(\la V_s(\mu), h(X_{T-s})\ra+ v_s'\big) ds\Big], \quad t \in [0,T],
		\end{split}
	\end{equation*}
for all finite measures $\mu$.
\end{lemma}
\begin{proof}
	The proof follows along similar lines as the proof of \cite[Lemma 2.6]{Mytnik98}, but we include it for the sake of completeness. Let $\mu$ be an arbitrary nonnegative finite measure,
	and let $(\mu_k)_{k \geq 1}$ be a sequence in $\mathcal{C}_c^\infty$ that converges weakly to $\mu$. Then it follows that $V_\cdot(\mu_k)\in \mathcal{C}^{1,2}(\R_+\times \R)\cap L^2_{loc}(\R_+\times \R)$, where  $L^2_{loc}(\R_+\times \R)$ denotes the set of measurable functions $\psi\colon [0,\infty)\times \R \rightarrow \R$ such that 
	\begin{equation}\label{eq:L2est_v}
		\int_0^t \int_\R \psi_s^2(x)dxds < \infty,  \quad \forall t >0. 
	\end{equation} 
	Now, let $T\geq 0$ and $t\in [0,T]$. 
	Thus, using \eqref{eq:L2est_v} we get by It\^o's formula 
	\begin{equation}\label{eq:ito_1}
		\begin{split}
			&\E^X\Big[\exp(-\la V_{0}(\mu_k), X_{T}\ra)\Big] = \E^X\Big[\exp(-\la V_{T-t}(\mu_k),X_t\ra +v_{T-t})\Big]\\
			& \qquad  -\E^X\Big[ \int_t^T \exp\big(-\la V_{T-s}(\mu_k), X_{s}\ra+v_{T-s}\big) \big(\la V_{T-s}(\mu_k), h(X_{s})\ra+v_{T-s}' \big)ds\Big],
		\end{split}
	\end{equation}
	since $V_t(\mu_k)$ solves \eqref{eq:PDE} and 
	\begin{equation}\label{eq:martingale_1}
		\E^X\Big[ \int_t^T \exp(-\la V_{T-s}(\mu_k), X_{s}\ra+ v_{T-s})dM_{s}(V_{T-\cdot}(\mu_k))\Big] = 0.
	\end{equation}
Here, $M$ is an orthogonal martingale measure such that for all $\psi\in \mathcal{C}(\R_+\times \R) \cap L^2_{loc}(\R_+\times \R)$, $M_t(\psi)$ is a square integrable martingale with quadratic variation given by 
	\begin{equation*}
		\la M(\psi) \ra_t = \int_0^t \la \psi_s^2,X_s\ra ds, \quad t\geq0.
	\end{equation*}
	Hereby, \eqref{eq:martingale_1} follows due to  the fact that
	\begin{equation*}
		M_{s}(V_{T-\cdot}(\mu_k)) 
	\end{equation*}
is a square-integrable martingale since $V(\mu_k)$ satisfies \eqref{eq:PDE} and since,
	by Lemma~\ref{lem:x_bound} and \eqref{eq:L2est_v}, we have for $t\in [0,T]$ that
	\begin{equation*}
		\begin{split}
		&\E[\la M_{\cdot}(V_{T-\cdot}(\mu_k) \ra_t] =\E\Big[\int_0^t \int_\R V_{T-s}(\mu_k)(x)^2 X_{s}(x)dxds\Big]\\
		& \leq \int_0^t \int_\R V_{T-s}(\mu_k)(x)^2dxds\sup_{s\in [0,T]}\sup_{x\in \R}\E[X_s(x)] <\infty.
		\end{split}
	\end{equation*}
	 This and the inequality 
	\begin{equation}
	\exp(-\la V_{T-s}(\mu_k), X_{s}\ra+ v_{T-s}) \leq \exp(v_{T-s}), \quad \forall s \leq t,
	\end{equation}
	 shows that $\int_0^t \exp(-\la V_{T-s}(\mu_k), X_s\ra+ v_{T-s})dM_s(V_{T-\cdot}(\mu_k))$ is a martingale, which implies \eqref{eq:martingale_1}.
	It remains to show the convergence as $k \rightarrow \infty$. Using \cite[Lemma 2.1]{Mytnik02}, the dominated convergence theorem  and setting $\tilde{t} = T-t$ in \eqref{eq:ito_1} yields
	\begin{equation*}
		\begin{split}
			&\E^X\Big[\exp(-\la V_{\tilde{t}}(\mu), X_{T-\tilde{t}}\ra + v_{\tilde{t}})\Big] = \E^X\Big[\exp(-\la V_0(\mu),X_{T}\ra)\Big]\\
			& \qquad  +\E^X\Big[ \int_0^{\tilde{t}} \exp(-\la V_s(\mu), X_{T-s}\ra+ v_s) \big(\la V_s(\mu), h(X_{T-s})\ra+ v_s'\big) ds\Big], \quad\forall  \tilde{t} \in [0,T].
		\end{split}
	\end{equation*}
\end{proof}
\noindent Fix arbitrary $m\in \R_+$, $x_0 \in \R$ and $\mu\in \mathcal{M}_F^+$. 
Until the end of this section let $Y^n$ denote the process $Y^{n,m}$ constructed in Section~\ref{sec:convergence_dual} with $Y_0 = \mu +m\delta_{x_0}$. We also suppress the upper index $m$ from $T_{n,k}^m$, $\bar{T}_{n,k}^m$, $U_{n,k}^m$, $Z_{n,k}^m$ and $M_{n,k}^m$. 
Furthermore, let the processes $Y^n$ be independent of $X$ for all $n \geq 1$. 
To this end define the point process
\begin{equation*}
	p^{(n)} \colon D_{p^{(n)}} \rightarrow \{1,2\} \times\R_+\times \R ,
\end{equation*}
with 
\begin{equation*}
	D_{p^{(n)}} = \{\bar{T}_{n,1}, \bar{T}_{n,2}, \ldots\},
\end{equation*}
and
\begin{equation*}
	p^{(n)}(\bar{T}_{n,k}) = (M_{n,k},Z_{n,k}, U_{n,k}), \forall k \geq 1.
\end{equation*}
The associated counting measure is given by
\begin{equation*}
	\mathfrak{P}^{(n)}(t,B) = \#\{ s\in D_{p^{(n)}} | s \leq t, p^{(n)}(s )\in B)\},
\end{equation*}
for measurable $B\subseteq \{1,2\} \times\R_+\times \R $.
\begin{lemma}\label{lem:compensator}
	The compensator of $\mathfrak{P}^{(n)}$ is given by
	\begin{equation}\label{eq:compensator}
		\hat{\mathfrak{P}}^{(n)} (t, B_1 \times B_2\times B_3) = \int_0^t \sum_{i = 1}^2 \mathbbm{1}_{i \in B_1}\int_{ B_2} \int_{B_3} Y^{n}_{s-}(x) dx \nu_n^i(d\lambda)ds,
	\end{equation}
for measurable sets $B_1 \subseteq \{1,2\}$, $B_2 \subseteq \R_+$ and $B_3 \subseteq \R$.
\end{lemma}
\begin{proof}
	We have that 
	\begin{equation*}
		\E[\mathfrak{P}^{(n)}(t,B)] \leq \E[\mathfrak{N}^{(n)}(t,\{1,2\} \times\R_+\times \R)] \leq \nu_n((0,\infty)) \E\Big[\int_0^t\int_\R Y_{s-}^n(x)dxds\Big] \leq C(T),
	\end{equation*}
	where the last inequality follows by Proposition~\ref{prop:main}. This means that $\mathfrak{P}^{(n)}$ is a monotone increasing, adapted, integrable process and thus, by the Doob-Meyer theorem the compensator exists. A calculation shows that it takes the form \eqref{eq:compensator}.
\end{proof}

\begin{lemma}\label{lem:approx_dual}
Let $(X_t)_{t\geq 0}$ be a solution to \eqref{eq:SPDE1} with $X_0 \sim \eta \in \mathcal{P}_p(\mathcal{C}_{tem}^+)$, $p>2$, such that \eqref{eq:cond_x_0} holds and let $a$ be defined as in \eqref{eq:def_d}. Let $Y^n$ be as above, independent of $X$. Then there exists a time $T_0>0$ depending only on $\la \mu,1\ra$ and $a$ such that for all $t \leq T \leq T_0$  we have
	\begin{equation*}
		\begin{split}
			&\E\Big[(-1)^{J_t^n}\exp\Big(-\la Y_t^n, X_{T-t}\ra-a \int_0^t\la Y_s^n, 1 \ra ds\Big)\Big] = \E\Big[\exp(-\la Y_0,X_T\ra)\Big]\\
			& + \E\Big[ \int_0^t(-1)^{J_{s-}^n} \exp\Big(-\la Y_{s-}^n, X_{T-s}\ra-a \int_0^s\la Y_r^n, 1 \ra dr\Big) \Big(\la Y_{s-}^n, h(X_{T-s})\ra - \la Y_{s-}^n, h_n(X_{T-s})\ra \Big) ds\Big].
		\end{split}
	\end{equation*}
\end{lemma}
\begin{proof}
	 First, choose $T_0 >0$ small enough so that  Lemma~\ref{lem:exponential_moment} is satisfied with $\gamma = 2a$ and $K = \nu((0,\infty])$, let $T \leq T_0$ and set $v_t = -a\int_0^t\la Y^n_s,1\ra ds$ for $t \geq 0$.
	By Lemma~\ref{lem:ito} we have 
	\begin{equation}\label{eq:duality_pf_1}
		\begin{split}
			&\E^X\Big[(-1)^{J^n_{\bar{T}_{n,k+1}-}}\exp\Big(-\la Y_{\bar{T}_{n,k+1}-}^n, X_{T-\bar{T}_{n,k+1}}\ra -a \int_0^{\bar{T}_{n,k+1-}}\la Y_s^n, 1 \ra ds\Big)\Big] \\
			&= \E^X\Big[(-1)^{J^n_{\bar{T}_{n,k}}}\exp\Big(-\la Y_{\bar{T}_{n,k}},X_{T-\bar{T}_{n,k}}\ra -a \int_0^{\bar{T}_{n,k}}\la Y_s^n, 1 \ra ds\Big)\Big]\\
			& \quad  + \E^X\Big[ \int_{\bar{T}_{n,k}}^{\bar{T}_{n,k+1}} (-1)^{J^n_{s-}}\exp\Big(-\la Y_{s-}^n, X_{T-s}\ra -a \int_0^{s}\la Y_r^n, 1 \ra dr\Big) \la Y_{s-}^n, h(X_{T-s})-a\ra ds\Big],
		\end{split}
	\end{equation} 
	since $t \mapsto X_t$ is continuous and $J_{\bar{T}_{n,k+1}-} = J_s$ for all $s \in (\bar{T}_{n,k}, \bar{T}_{n,k+1})$. Let $\pmb{\Delta}Y_s = Y_{s}-Y_{s-}$ for all $s\geq 0$, then we get 
	\begin{equation*}
		\begin{split}
			&	\E^X\Big[(-1)^{J^n_{\bar{T}_{n,k+1}}}\exp\Big(-\la Y_{\bar{T}_{n,k+1}}^n, X_{T-\bar{T}_{n,k+1}}\ra-a \int_0^{\bar{T}_{n,k+1-}}\la Y_s^n, 1 \ra ds\Big)\Big] \\
			& = \E^X\Big[(-1)^{J^n_{\bar{T}_{n,k+1}-}}\exp\Big(-\la Y_{\bar{T}_{n,k+1}-}^n, X_{T-\bar{T}_{n,k+1}}\ra+v_{\bar{T}_{n,k+1}}\Big)\Big] \\
			& +\exp\big(v_{\bar{T}_{n,k+1}}\big) \E^X\Big[(-1)^{J_{\bar{T}_{n,k+1}}}\exp\big(-\la Y_{\bar{T}_{n,k+1}}^n, X_{T-\bar{T}_{n,k+1}}\ra\big) \\
			&\hspace{7cm} -(-1)^{J^n_{\bar{T}_{n,k+1}-}}\exp\big(-\la Y_{\bar{T}_{n,k+1}-}^n, X_{T-\bar{T}_{n,k+1}}\ra\big)\Big] \\
			& = \E^X\Big[(-1)^{J^n_{\bar{T}_{n,k}}}\exp\Big(-\la Y_{\bar{T}_{n,k}},X_{T-\bar{T}_{n,k}}\ra +v_{\bar{T}_{n,k}}\Big)\Big]\\
			&+ \E^X\Big[ \int_{\bar{T}_{n,k}}^{\bar{T}_{n,k+1}} (-1)^{J^n_{s-}}\exp\Big(-\la Y_{s-}^n, X_{T-s}\ra +v_{s}\Big) \la Y_{s-}^n, h(X_{T-s})-a\ra ds\Big]\\
			& -\exp\big(v_{\bar{T}_{n,k+1}}\big)\E^X\Big[(-1)^{J^n_{\bar{T}_{n,k+1}-}}\exp\big(-\la Y_{\bar{T}_{n,k+1}-}^n, X_{T-\bar{T}_{n,k+1}}\ra\big) \\
			& \times \Big((1-\exp(-\la \pmb{\Delta} Y^n_{\bar{T}_{n,k+1}}, X_{T-\bar{T}_{n,k+1}}))\mathbbm{1}_{M_{n,k+1} = 1} +\big(1+\exp\big(-\la \pmb{\Delta} Y^n_{\bar{T}_{n,k+1}}, X_{T-\bar{T}_{n,k+1}}\ra\big)\big)\mathbbm{1}_{M_{n,k+1} = 2}\Big)\Big],
		\end{split}
	\end{equation*}
	where in the second equality we used \eqref{eq:duality_pf_1} and the definition of $J^n$.
	This leads to 
	\begin{equation}\label{eq:duality_pf_3}
		\begin{split}
			&\E^X\Big[(-1)^{J^n_t}\exp\big(-\la Y_t^n, X_{T-t}\ra\big)\exp\Big(-a \int_0^t\la Y^n_s,1\ra ds\Big)\Big] = \E\Big[\exp\big(-\la Y_0,X_T\ra\big)\Big]\\
			& +\E^X\Big[ \int_0^t (-1)^{J^n_{s-}}\exp\Big(-\la Y_{s-}^n, X_{T-s}\ra-a \int_0^{s}\la Y_r^n, 1 \ra dr\Big) \la Y_{s-}^n, h(X_{T-s})-a\ra ds\Big] \\
			& - \int_0^t \int_{\{1,2\}} \int_{\R_+}\int_{\R}(-1)^{J^n_{s-}}\E_X\Big[\exp\Big(-\la Y^n_{s-},X_{T-s}\ra - a \int_0^{s}\la Y_r^n, 1 \ra dr\Big) \\
			& \qquad \times \Big(\big(1- \ee^{-\lambda X_{T-s}(x)}\big)\mathbbm{1}_{r = 1} + (1+\ee^{-\lambda X_{T-s}(x)})\mathbbm{1}_{r = 2}\Big)\Big]\mathfrak{P}^{(n)}(dr,d\lambda,dx,ds).
		\end{split}
	\end{equation}
	By Lemma~\ref{lem:compensator} the compensator of $\mathfrak{P}^{(n)}$ is given by \eqref{eq:compensator}.
	Furthermore, by Lemma~\ref{lem:exponential_moment} it follows that 
	\begin{equation*}
		\E\Big[\exp\Big(- a \int_0^{T_0}\la Y_r^n, 1 \ra dr\Big)\Big] < \infty.
	\end{equation*} 
	Now take the expectation $\E^Y$ on both sides of \eqref{eq:duality_pf_3} and use the definition of $h_n$ and $h$, given by \eqref{eq:h_nduality}  and \eqref{eq:h_duality}, respectively, to get the desired result.
\end{proof}

\begin{corollary}\label{cor:approximate_duality}
	Let $(X_t)_{t\geq 0}$ be a solution to \eqref{eq:SPDE1} with $X_0 \sim \eta \in \mathcal{P}_p(\mathcal{C}^+_{tem})$, $p>2$,  and let $a$ be defined as in \eqref{eq:def_d}. Let $Y^n$ be as above, independent of $X$.
	Then there exists $T_0 >0$ , depending only on $a$ and $\la \mu,1\ra$, such that for each $X_0 \sim \eta \in \mathcal{P}_p(\mathcal{C}^+_{tem})$ satisfying \eqref{eq:cond_x_0}, we get
	\begin{equation*}
		\lim_{n \rightarrow \infty}\E[(-1)^{J^n_T}\exp\Big(-\la Y_T^n, X_{0}\ra-a \int_0^T\la Y_s^n, 1 \ra ds\Big)\Big] = \E[\exp(-\la Y_0,X_T\ra)],
	\end{equation*}
for all $T \leq T_0$. 
\end{corollary}
\begin{proof}
We only prove the claim for the case $a \leq 0$, since the proof follows along similar lines when $a>0$ with only minor modifications.
	By Lemma~\ref{lem:approx_dual} there exists $T_0>0$ such that for all $T\leq T_0$ we have
	\begin{equation*}
		\begin{split}
			&\lim_{n \rightarrow \infty}\E\Big[(-1)^{J^n_T}\exp\Big(-\la Y_T^n, X_{0}\ra-a \int_0^T\la Y_s^n, 1 \ra ds\Big)\Big] \\
			&  = \E[\exp\big(-\la Y_0,X_T\ra\big)]\\
			&+\lim_{n \rightarrow \infty} \E\Big[ \int_0^T(-1)^{J^n_{s-}}\exp\Big(-\la Y_{s-}^n, X_{T-s}\ra - a \int_0^s \la Y_{r-}^n, 1\ra dr\Big) \Big(\la Y_{s-}^n, h(X_{T-s})\ra - \la Y_{s-}^n, h_n(X_{T-s})\ra \Big) ds\Big].
		\end{split}
	\end{equation*}
	Furthermore, by  Lemma~\ref{lem:conv_yh} and Lemma~\ref{lem:exponential_moment} we have 
	\begin{equation*}
		\begin{split}
			&\lim_{n \rightarrow \infty}\Big| \E\Big[ \int_0^T (-1)^{J^n_{s-}}\exp\Big(-\la Y_{s-}^n, X_{T-s}\ra -a \int_0^s \la Y_{r-}^n, 1\ra dr\Big) \Big(\la Y_{s-}^n, h(X_{T-s})\ra - \la Y_{s-}^n, h_n(X_{T-s})\ra \Big) ds\Big]\Big| \\
			& \leq \lim_{n \rightarrow \infty} \E\Big[\exp\Big(-a \int_0^T \la Y_s^n,1\ra ds\Big) \int_0^T \la Y_{s-}^n , |h_n(X_{T-s})-h(X_{T-s})|\ra ds\Big] \\
			& \leq \lim_{n \rightarrow \infty} \E\Big[\exp\Big(-2a \int_0^T \la Y_s^n,1\ra ds\Big)\Big]^{1/2} \E\Big[\Big(\int_0^T \la Y_{s-}^n , |h_n(X_{T-s})-h(X_{T-s})|\ra ds\Big)^2\Big]^{1/2} =0.
		\end{split}
	\end{equation*}
	This concludes the proof.
\end{proof}
Finally, we can prove Proposition~\ref{prop:duality}.
\begin{proof}[Proof of Proposition~\ref{prop:duality}]
	The assertion follows immediately from Corollary~\ref{cor:approximate_duality} and Lemma~\ref{lem:conv_drift_exp}.
\end{proof}
\section{Proof of Weak Existence}\label{sec:weak_exist}
In this section, we prove the weak existence part of Theorem~\ref{th:uniqueness}. Again, our arguments were inspired by \cite{Barnes24}.\\
 Note that we only need to prove weak existence in the case where $h$ is discontinuous due to \cite[Theorem 1.1]{Shiga94}.  That means we need to consider the case where $b_0\neq b_1$. For simplicity, we thus only consider the case where $\nu^1 = \nu^2\equiv 0$. 
The proof relies on considering a sequence $(X^n)_{n\geq 1}$ of solutions to \eqref{eq:SPDE1}, where $h$ is approximated by continuous functions $h_n$  given by
\begin{equation*}
 h_n(x) = d_1(1-\ee^{- n x}) + d_2(1+\ee^{-n x})+a, \quad x\geq 0,
\end{equation*} 
where $d_1,d_2$ and $a$ are given according to \eqref{eq:def_d}. The sequence of functions $(h_n)_{n\geq1}$ converges boundedly pointwise to $h$.
 As mentioned above, we get the existence of a weak solution $X^n$ by applying \cite[Theorem 1.1]{Shiga94} for each $n\in N$. 
When $h$ is replaced by $h_n$, we also get the weak uniqueness of solutions by applying the results from Section~\ref{sec:weak_unique}. Furthermore, we can apply the duality relation from Corollary~\ref{cor:approximate_duality} for $X^n$, which is conditional on its existence. \\ 
In a first step, we prove tightness of the sequence $(X^n)_{n\geq 1}$, which can be achieved in a similar manner to the proofs in \cite{Shiga94}. 
Then, after passing to a weakly convergent subsequence and applying the Skorohod
representation theorem, we realize copies of this subsequence
on a common probability space such that they converge almost surely
to a limit $X$. We then show that $X$ is a solution to
\eqref{eq:SPDE1}. The main difficulty is to prove, on this common probability space, that the drift terms $h_n(X^n)$ converge to
$h(X)$  in a suitable $L^2$ sense, despite the discontinuity of $h$. Similar to the approach from \cite{Barnes24}, we use the duality
relation of Proposition~\ref{prop:duality}, together with the
convergence results for the approximating dual processes established
in Section~\ref{sec:duality}, to prove this claim.

\subsection{Tightness}
Recall that $X_0 \in \mathcal{P}_p(\mathcal{C}_{tem})$ for $p>2$. We can apply standard methods to derive the following lemma.
\begin{lemma}\label{lem:tight}
\label{lem:tightness-approximating-solutions}
Let \(T>0\), let \(p>2\), and let \(\eta \in \mathcal P_p(\mcC^+_{tem})\).
Let \((X^n)_{n\ge 1}\) be \(\mcC([0,T],\mcC^+_{tem})\)-valued weak solutions of
\[
  dX^n_t(x)
  =
  \frac12 \Delta X^n_t(x)\,dt
  + h_n(X^n_t(x))\,dt
  + \sqrt{X^n_t(x)}\, \dot{W}^n(t,x),
  \qquad X^n_0 \sim \eta,
\]
where $ \dot{W}^n$ is space-time white noise.
Then the sequence \((X^n)_{n\ge1}\) is tight in
\[
  \mcC([0,T],\mcC^+_{tem}).
\]
Consequently, since \(T>0\) is arbitrary, \((X^n)_{n\ge1}\) is tight in
\(\mcC(\mathbb R_+,\mcC^+_{tem})\), endowed with the topology of locally uniform
convergence.
\end{lemma}

\begin{proof}
It is enough to prove tightness in \(\mcC([0,T],\mcC_{tem})\), since
\(\mcC^+_{tem}\) is closed in \(\mcC_{tem}\).  We use the seminorms
\[
  \|f\|_{-\lambda}
  :=
  \sup_{x\in\mathbb R} e^{-\lambda |x|}|f(x)|,
  \qquad \lambda>0,
\]
and the metric \(d_{\mathrm{tem}}\) defined as
\[
  d_{\mathrm{tem}}(f,g)
  :=
  \sum_{k=1}^{\infty} 2^{-k}
  \bigl(\|f-g\|_{-1/k}\wedge 1\bigr)
\]
and endow $ \mcC([0,T],\mcC^+_{tem})$ with the uniform metric induced by $d_{\mathrm{tem}}$.

By the mild formulation,
\[
  X^n_t = S_tX^n_0 + A^n_t + Z^n_t,\qquad 0\le t\le T,
\]
where \(S_t\) denotes the heat semigroup with generator \(\frac12\Delta\),
\[
  A^n_t(x)
  :=
  \int_0^t S_{t-s}\bigl(h_n(X^n_s)\bigr)(x)\,ds,
\]
and
\[
  Z^n_t(x)
  :=
  \int_0^t\int_{\mathbb R}
    p_{t-s}(x-y)\sqrt{X^n_s(y)}\,W^n(ds,dy).
\]

We first consider the initial semigroup term.  Define
\[
  \Phi:\mcC_{tem}\to \mcC([0,T],\mcC_{tem}),
  \qquad
  \Phi(f)(t):=S_tf .
\]
Standard estimates show that $\Phi$ is continuous and thus the sequence \((S_\cdot X^n_0)_{n\ge1}\) is tight in \(\mcC([0,T],\mcC_{tem})\).
The tightness of  \((A^n)_{n\ge1}\) follows by standard heat kernel estimates and the fact that \(\sup_n\sup_{r\ge0}|h_n(r)|<\infty\). 

It remains to treat the stochastic convolutions \((Z^n)_{n\ge1}\).  To this end, let $\epsilon >0$ be fixed but arbitrary.  Fix a compact set
\(G_\epsilon\subset \mcC^+_{tem}\)
such that $\eta(G_\epsilon^c)<\epsilon$
and put
\[
 E^n_G := \{X^n_0\in G_\epsilon\}.
\]
Once again, using standard estimates as in \cite[Section 6]{Shiga94} one can prove that there exists a compact set $K_{G_\epsilon} \subset \mcC([0,T],\mcC_{tem})$ so that
\begin{equation*}
    \sup_{n\geq 1} \mathbb{P}(Z^n \in K_{G_\epsilon}^c, E^n_{G_\epsilon}) < \epsilon. 
\end{equation*}
Since $\epsilon >0$ was arbitrary this implies that $(Z^n)_{n\geq 1}$ is tight in $\mcC([0,T],\mcC_{tem})$. 

We have proved tightness of \((S_\cdot X^n_0)_{n\ge1}\), \((A^n)_{n\ge1}\), and
\((Z^n)_{n\ge1}\).    
Therefore
\[
  X^n = S_\cdot X^n_0 + A^n + Z^n
\]
is tight in \(\mcC([0,T],\mcC_{tem})\).  Since each \(X^n\) is
\(\mcC^+_{tem}\)-valued, the sequence is tight in
\(\mcC([0,T],C^+_{tem})\).

Finally, tightness in \(\mcC(\mathbb R_+,\mcC^+_{tem})\), endowed with locally uniform
convergence, follows from tightness of the restrictions to \(\mcC([0,m],\mcC^+_{tem})\)
for every \(m\in\mathbb N\).
\end{proof}

\noindent Lemma~\ref{lem:tight} demonstrates by \cite[Lemma 6.3]{Shiga94} that 
the sequence $(X^n)_{n \geq 1}$ of $\mathcal{C}(\R_+, \mathcal{C}_{tem}(\R))$ valued random variables is tight. Thus, by Prohorov's theorem and Skorohod's representation theorem, there exists a subsequence $(X^{n_k})_{k\geq 1}$ and another sequence $(\tilde{X}^{n_k})_{k\geq 1}$ of  $\mathcal{C}(\R_+, \mathcal{C}_{tem}(\R))$-valued random variables such that the following holds.
\begin{enumerate}
	\item The sequence $(\tilde{X}^{n_k})_{k\geq 1}$ is defined on a common probability space;
	\item For all $k \geq 1$ we have that $X^{n_k} \overset{law}{=} \tilde{X}^{n_k}$;
	\item The sequence $(\tilde{X}^{n_k})_{k\geq 1}$ converges almost surely to a limit $X$ with respect to the topology of $\mathcal{C}(\R_+, \mathcal{C}_{tem}(\R))$.
\end{enumerate}
In the following we shall denote $\tilde{X}^{n_k}$ by $X^{n}$, with a slight abuse of notation, for the sake of readability. 
\subsection{Convergence of $h_n(X^n)$}

\begin{lemma}\label{lem:lpconv}
	There exists $T_0>0$ such that 
	 for all $x\in \R$ and $0\leq t \leq T_0$ it holds that
	\begin{equation*}
		 \exp(-nX^n_t(x)) \rightarrow \mathbbm{1}_{X_t(x) = 0}
	\end{equation*}
in $L^q$ for $q\in[1,\infty)$ as $n \rightarrow \infty$.  
\end{lemma}
\begin{proof}
	The proof is similar to the proof of \cite[Lemma 5.3]{Barnes24}. 
	First, we prove the $L^1$ convergence. Let $T_0>0$ be such that Proposition~\ref{prop:duality} holds for $X^{n}$, for all $n \geq 0$. Note that
	it is possible to choose $T_0>0$ uniformly in $n$, since it only depends on $\sup_n\nu_n^1(0,\infty)+\nu_n^2(0,\infty)$, where $\nu_n^1, \nu_n^2$ are given by
	\begin{equation*}
		\nu_n^i = \nu^i_{|[0,n]} + d_i\delta_n, \quad i \in \{1,2\}, n\geq 1.
	\end{equation*}
	 For simplicity denote $X_t(x)$ by $Z$ and $X^n_t(x)$ by $Z_n$. 
	Below we use the fact that almost surely $Z_n, Z \geq 0$.
	Thus we get 
	\begin{equation}\label{eq:hnXn1}
		\begin{split}
			&\E[|\exp(-nZ_n)- \mathbbm{1}_{Z =0}|] \\
			&= \E[|\exp(-nZ_n)- \mathbbm{1}_{Z =0}| \mathbbm{1}_{Z = 0}] +  \E[|\exp(-nZ_n)- \mathbbm{1}_{Z =0}| \mathbbm{1}_{Z >0}] \\
			& =\E[|\exp(-nZ_n)- 1| \mathbbm{1}_{Z = 0}] +  \E[|\exp(-nZ_n)| \mathbbm{1}_{Z >0}] \\
			& =\E[\mathbbm{1}_{Z =0}] - \E[\exp(-nZ_n)] + 2 \E[|\exp(-nZ_n)| \mathbbm{1}_{Z >0}].
		\end{split}
	\end{equation}
Since $Z_n \rightarrow Z$ almost surely, we get that $ 2 \E[|\exp(-nZ_n)| \mathbbm{1}_{Z >0}] \rightarrow 0$. By Lemma~\ref{lem:conv_drift_exp} and by applying Corollary~\ref{cor:approximate_duality} with $X^n$, we furthermore get that 
\begin{equation*}
\begin{split}
&	\E[\mathbbm{1}_{Z =0}]  = \lim_{m\rightarrow \infty}\E[\exp(-mX_t(x))] \\
&=  \lim_{m\rightarrow \infty}\lim_{n\rightarrow \infty}\E[\exp(-mX^{n}_t(x))] \\
&= 	\lim_{m\rightarrow \infty}\lim_{n \rightarrow \infty}\E\Big[ (-1)^{J_t^{n,m}}\exp\Big(-\la Y^{n,m}_{t}(x), X_0\ra-a\int_0^t \la Y_s^{n,m},1\ra ds\Big)\Big]\\
&= \E\Big[ (-1)^{J_t^{\infty,\infty}}\exp\Big(-\la Y^{\infty,\infty}_{t}(x), X_0\ra-a\int_0^t \la Y_s^{\infty,\infty},1\ra ds\Big)\Big],
\end{split}
\end{equation*}
and 
\begin{equation*}
\begin{split}
&\lim_{n \rightarrow \infty}	\E[\exp(-nZ_n)] = \lim_{n\rightarrow \infty}\E\Big[ (-1)^{J_t^{n,n}}\exp\Big(-\la Y^{n,n}_{t}(x), X_0\ra-a\int_0^t \la Y_s^{n,n},1\ra ds\Big)\Big]\\
& \quad= \E\Big[ (-1)^{J_t^{\infty,\infty}}\exp\Big(-\la Y^{\infty,\infty}_{t}(x), X_0\ra-a\int_0^t \la Y_s^{\infty,\infty},1\ra ds\Big)\Big].
\end{split}
\end{equation*}
Thus it follows that
\begin{equation*}
	 \E[\mathbbm{1}_{Z =0}] - \E[\exp(-nZ_n)] \rightarrow 0,
\end{equation*}
as $n \rightarrow \infty$.
By the above and \eqref{eq:hnXn1} the $L^1$ convergence of $\exp(-nX_t^n(x))$ to $\mathbbm{1}_{X_t(x) = 0}$ follows. From this it is standard to derive the convergence in $L^q$ for $q \in [1,\infty)$. 
\end{proof}

\begin{proof}[Proof of Weak Existence]
	First, we prove weak existence of solutions up to time $T_0 >0$ such that Proposition~\ref{prop:duality} holds. Note that $T_0$
	does not depend on the initial condition. Thus, we can extend our solution provided that 
	$\P\circ X_{T_0}^{-1} \in \mathcal{P}_p(\mathcal{C}_{tem}^+)$ and satisfies \eqref{eq:cond_x_0}. This, however, readily follows by Lemma~\ref{lem:x_bound} and by the estimate 
    \begin{equation}\label{eq:est_Pp}
  \sup_{n\ge1}\sup_{0\le s\le t}
  \int_{\mathbb R}
    e^{-\lambda |x|}
    \mathbb E\!\left[(X^n_s(x))^p\right]\,dx
  <\infty,
    \end{equation} 
    which follows by standard estimates as in \cite[Section 6]{Shiga94}.
	Fix arbitrary $\phi \in \mathcal{C}_c^\infty(\R)$. Recall that according to our conventions $X^n \rightarrow X$ in $\mathcal{C}(\R_+,\mathcal{C}_{tem}^+)$ almost surely as $n \rightarrow \infty$. Since \eqref{eq:SPDE1} with drift $h_{n}$ for $n \geq 1$ is unique in law, we have that almost surely, 
	\begin{equation}\label{eq:martingale_sbm}
	\begin{split}
	&	\int_\R \phi (x)X^{n}_t(x)dx - \int_\R \phi (x)X_0^n(x) dx \\
	&\qquad= \frac{1}{2}\int_\R\int_0^t \phi''(x) X_s^{n}(x) dsdx + \int_\R \int_0^t \phi(x) h_{n}(X_s^{n}(x))dsdx + M_t^{n}(\phi),
	\end{split}
	\end{equation}
where $(M_t^{n}(\phi))_{t\geq 0}$ is an $L^2$-martingale with quadratic variation
\begin{equation*}
	\la M_{\cdot}^{n}(\phi) \ra_t = \int_\R \int_0^t \phi(x)^2 X_s^{n}(x)dsdx.
\end{equation*}
From Lemma~\ref{lem:lpconv} we know that the second term on the right-hand side converges to 
\begin{equation*}
	 \int_\R \int_0^t \phi (x)h(X_s(x))dsdx 
\end{equation*}
in $L^2$. Furthermore,  due to \eqref{eq:est_Pp} and since $p >2$ we have that $(X^{n}_\cdot-X_\cdot)^2$ on the interval $[0,t]$ and restricted to the support of $\phi$ is uniformly 
integrable. For this reason, 
\begin{equation*}
	\frac{1}{2}\int_\R\int_0^t \phi''(x) X_s^{n}(x) dsdx \rightarrow \frac{1}{2}\int_\R\int_0^t \phi''(x) X_s(x) dsdx,
\end{equation*}
in $L^2$ as $n\rightarrow \infty$ and 
\begin{equation*}
		\int_\R \phi (x)X^{n}_t(x)dx \rightarrow \int_\R \phi (x)X_t(x)dx
\end{equation*}
in $L^2$ as $n\rightarrow \infty$.
This implies that  $M_t^{n}(\phi)$ converges in $L^2$ to a random variable $M_t(\phi)$ for all $t \in [0,T_0]$. Now, it is standard to deduce that 
$(M_t(\phi))_{t\in [0,T_0]}$ is a $L^2$-martingale with quadratic variation 
\begin{equation*}
		\la M_{\cdot}(\phi) \ra_t = \int_\R \int_0^t \phi(x)^2 X_s(x)dsdx, \quad t \in [0,T_0].
\end{equation*}
Thus, all the terms in \eqref{eq:martingale_sbm} converge and we get
\begin{equation*}
	\begin{split}
		&	\int_\R \phi (x)X_t(x)dx - \int_\R \phi (x)X_0(x) dx \\
		&\qquad= \frac{1}{2}\int_\R\int_0^t \phi''(x) X_s(x) dsdx + \int_\R \int_0^t \phi(x) h(X_s(x))dsdx + M_t(\phi),
	\end{split}
\end{equation*}
for $t \in [0,T_0]$. Finally, using standard arguments as in \cite[Proof of Lemma 2.4]{Konno88} and after extending the filtered probability space so that it supports an additional, independent space-time white noise, it is possible to show that $X$ satisfies \eqref{eq:SPDE1}.
\end{proof}
\section{Proof of Weak Uniqueness}\label{sec:weak_unique}
\begin{lemma}\label{lem:one_dim_unique}
	Let $X^1, X^2$ be two solutions to \eqref{eq:SPDE1} with $X_0^1,X_0^2\sim \mu \in \mathcal{P}_p(\mathcal{C}_{tem}^+) $ such that \eqref{eq:cond_x_0} is satisfied. Then, it follows that 
	\begin{equation*}
	X_t^1 \overset{d}{=} X_t^2,
	\end{equation*}
for all $t \in [0,T_0]$, where $T_0$ only depends on $a$.
\end{lemma}

\begin{proof}
First, we choose $T_0$ as in Proposition~\ref{prop:duality} under the assumption that $\la Y_0,1\ra =1$.
Now, fix $t \in [0,T_0]$ and let $f$ be a fixed but arbitrary nonnegative simple function. Then it follows for $g = \frac{1}{\|f\|_{L^1}}f$ that 
\begin{equation}\label{eq:equality}
	\E[\exp(- \lambda \la g , X_t^1 \ra)] = \E[\exp(- \lambda \la g , X_t^2 \ra)] \quad \forall \lambda \in [0,1],
\end{equation}
by Proposition~\ref{prop:duality}.
Since the map $\lambda \mapsto	\E[\exp(- \lambda \la g , X_t^i \ra)] $ is analytic on $(0, \infty)$ for $i \in \{1,2\}$,
\eqref{eq:equality} extends to all $\lambda \geq 0$. In particular, we have that 
\begin{equation*}
	\E[\exp(- \lambda \la f, X_t^1 \ra)] = \E[\exp(- \lambda \la f , X_t^2 \ra)] \quad \forall \lambda \geq 0.
\end{equation*}
Since $f$ was arbitrary it follows by \cite[Corollary 2.3]{Kallenberg17} that $X^1_t \overset{d}{=} X_t^2$.
\end{proof}

\begin{thm}
	The solution to \eqref{eq:SPDE1} is weakly unique in $\mathcal{C}( \R_+,\mathcal{C}_{tem}(\R))$.
\end{thm}
\begin{proof}
	It is enough to prove weak uniqueness on the interval $[0,T_0]$ and that $\P \circ (X_{T_0})^{-1} \in \mathcal{P}_p(\mathcal{C}_{tem}^+) $ and that $X_{T_0}$ satisfies \eqref{eq:cond_x_0},  since $T_0$ does not depend on the initial value $X_0$.
    Condition \eqref{eq:cond_x_0} follows from Lemma~\ref{lem:x_bound}. Since we are only considering $\mcC_{tem}^+$ valued solutions, we only need to verify the integrability condition to verify that $\P \circ (X_{T_0})^{-1} \in \mathcal{P}_p(\mathcal{C}_{tem}^+)$, which follows by standard estimates and the fact that $h$ is bounded.
	Thus, it remains to prove uniqueness on $[0,T_0]$. However, by using Lemma~\ref{lem:one_dim_unique}, this follows along 
	the same lines as the proof of \cite[Theorem 4.2]{Ethier05}. It is only necessary to show that $\P \circ (X_{t})^{-1} \in \mathcal{P}_p(\mathcal{C}_{tem}^+) $ and that $X_{t}$ satisfies \eqref{eq:cond_x_0} for all $t \in [0,T_0]$. 
    However, this once again follows from Lemma~\ref{lem:x_bound} and standard estimates.
\end{proof}

\section{Proof of Theorem~\ref{thm:scp}}\label{sec:scp}
Before we prove Theorem~\ref{thm:scp}, we derive the following helpful lemma.
\begin{lemma}\label{lem:integral_ineq}
	Let $s>0$, $a>0$ and $x,y \in (-\infty,a)$ such that $x \leq y$. Then 
	\begin{equation*}
		\int_{a}^\infty \mathcal{W}_s(x,0)(z)dz \leq \int_{a}^\infty \mathcal{W}_s(y,0)(z)dz. 
	\end{equation*}
\end{lemma}
\begin{proof}
	By Proposition~\ref{prp:very_sing_asymp} we have 
	\begin{equation*}
		\begin{split}
			&\int_{a}^\infty \mathcal{W}_s(x,0)(z)dz = s^{-1} \int_a^\infty f(s^{-1/2}|z-x|)dz\\
			& =  s^{-1} \int_{a-x+y}^\infty f(s^{-1/2}|z'-y|)dz' \leq \int_a^\infty \mathcal{W}_s(y,0)(z')dz',
		\end{split}
	\end{equation*}
	which proves the assertion.
\end{proof}

\begin{proof}[Proof of Theorem~\ref{thm:scp}]
We only prove the statement for the case  $b_1\geq 0$. When $b_1<0$ the statement follows by a comparison  with the solution when $b_1 = 0$, by using weak uniqueness, \cite[Corollary 2.4]{Shiga94} and standard approximation and tightness arguments. Thus, we assume without loss of generality that $b_1 \geq 0$ from now on. \\
The proof relies on the duality relation. To this end, let $Y$ be the process constructed in the beginning of Section~\ref{sec:duality} with $Y_0 = \infty\delta_0$ and let $I_t(Y_0)$ for $t\geq 0$ be defined as in \eqref{eq:defI_t}.
By Lemma~\ref{lem:exponential1}, there exists $T_0>0$ such that 
	\begin{equation}\label{eq:expon_support}
		\E[ \ee^{|I_{T_0}(Y_0)||}] < \infty.
	\end{equation}
  For $t= 0$ the assertion \eqref{eq:th_supp} follows by the assumptions of the theorem. Thus,  let $s\in (0,T_0]$. 
Using Proposition~\ref{prop:duality} and Remark~\ref{rem:T0} it follows that
\begin{equation}\label{eq:compact_initial}
	\begin{split}
		&\E \Big[ \int_\R \mathbbm{1}_{X_s(x)>0}dx\Big] = \lim_{n \rightarrow \infty} \E \Big[\int_\R \big(1-\ee^{-nX_s(x)}\big)dx\Big]\\
		& = \lim_{n \rightarrow \infty}  \int_\R \big(1-\E[\ee^{-nX_s(x)}]\big)dx\\
		& =  \lim_{n \rightarrow \infty} \int_\R \big(1-\E[\ee^{-\la X_0, Y_s(n\delta_x)\ra}]\big)dx.
	\end{split}
\end{equation}
Note, that it is enough to prove
\begin{equation*}
	\E \Big[\int_\R \mathbbm{1}_{X_s(x)>0}dx\Big]  < \infty.
\end{equation*}
We denote 
\begin{equation*}
	R_x = (|x|-R)/2,
\end{equation*}
where $R$ is such that $\supp(X_0) \subseteq B(0,R)$ and we also assume $R>2$.
Furthermore for $|x| > 2R$, denote
\begin{equation*}
	\tau_{R_x,n} = \inf\{ k \geq 0 \colon |U_{n,k}^x| \leq R_x+R \},
\end{equation*}
where $U_{n,k}^x$ is the location of the $k$-th jump of the process $Y(n \delta_x)$ and $U_{n,0}^x = x$.
Moreover, define
\begin{equation*}
	K_t^n = \max\{ k \in \N \colon \bar{T}_{n,k}^x \leq t\}, \quad t\geq 0,
\end{equation*}
where $\bar{T}_{n,k}^x $ is time of the $k$-th jump of the process $Y(n\delta_x)$. That is, $K^n_t$ is the number of jumps of $Y(n\delta_x)$ on the time interval $(0,t]$. 
We furthermore denote $T_{n,k}^x =\bar{T}_{n,k}^x - \bar{T}_{n,k-1}^x$, for $k \geq 1$. 
We obtain the following decomposition
\begin{equation} \label{eq:decomp_comapct}
	\begin{split}
		&\E[(1-\ee^{-\la X_0, Y_s(n\delta_x)})]\\
		& = \E[(1-\ee^{-\la X_0, Y_s(n\delta_x)}) \mathbbm{1}_{K_s^n >|x|^{1/4}}] \\
		& \quad +  \E[(1-\ee^{-\la X_0, Y_s(n\delta_x)}) \mathbbm{1}_{K_s^n \leq |x|^{1/4}}\mathbbm{1}_{\tau_{R_x,n} \leq K_s^n}]\\
		& \quad + \E[(1-\ee^{-\la X_0, Y_s(n\delta_x)}) \mathbbm{1}_{K_s^n \leq |x|^{1/4}}\mathbbm{1}_{\tau_{R_x,n} > K_s^n}]\\
		& = I_1^n(s,x) +I_2^n(s,x) +I_3^n(s,x)
	\end{split}
\end{equation}
Without loss of generality, we assume that $x >2R$.
First, we use \eqref{eq:expon_support} and the fact that $|K_s^n|$ is stochastically  dominated by $|I_s(Y_0)|$ in order to estimate $I_1^n(s,x)$. Hereby, the stochastic domination follows from \eqref{eq:ineq_times}. We arrive at
\begin{equation}\label{eq:I1}
	I_1^n(s,x)\leq \P(|I_s(Y_0)|^8 > x^2) \leq \frac{\E[|I_s(Y_0)|^8]}{x^2}\wedge1.
\end{equation}
Note, that \eqref{eq:expon_support} implies that $\E[|I_s(Y_0)|^8] <\infty$.
Next, we estimate $I_3^n(s,x)$ for $x>2R$. By using Lemma~\ref{lem:pde_est} and Proposition~\ref{prp:initial_trace_prop} we arrive at 
\begin{equation}\label{eq:I3}
	I_3^n(s,x) \leq \E[ \la X_0, Y_s(n\delta_x)\ra \mathbbm{1}_{K_s^n \leq |x|^{1/4}}\mathbbm{1}_{\tau_{R_x,n} > K^n_s}] \leq C\sum_{k=0}^{\lceil|x|^{1/4}\rceil} \la \mathbbm{1}_{(-\infty,R]},\mathcal{W}_{s-\bar{T}_{n,k}^x}(R_x+R,0) \ra,
\end{equation}
where we used Lemma~\ref{lem:integral_ineq} and the fact  that $U_{n,k}^x > R+R_x$ for all $0 \leq k \leq |x|^{1/4}+1$ and that $\supp(X_0) \subseteq B(0,R)$.
Note, that by Proposition~\ref{prp:very_sing_asymp} we have that $\mathcal{W}_t(0,0)(y) = t^{-1}f(t^{-1/2}|y|)$, $t >0, y\in \R$ such that
\begin{equation*}
	\sup_{y \in [s^{-1/2},\infty)}\frac{f(y)}{\ee^{-1/2 y^2}y} \leq c(s).
\end{equation*}
Therefore, for all $r>0$ we have
\begin{equation*}
	\begin{split}
	&\la \mathbbm{1}_{(-\infty,R]},\mathcal{W}_r(R_x+R,0)\ra \leq C\int_{-\infty}^R \mathcal{W}_r(R_x+R,0)(y)dy\\
	& = \int_{-\infty}^{-R_x} \mathcal{W}_r(0,0)(y)dy\leq C \int^{\infty}_{R_x} r^{-3/2}|y|\ee^{-\frac{y^2}{2r}}dy = Cr^{-1/2} \ee^{-\frac{R_x^2}{2r}},
	\end{split}
\end{equation*}
where we used that $R_x >1$ for $x>2R$. This leads to
\begin{equation}\label{eq:int_est_W_1}
	\la \mathbbm{1}_{(-\infty,R]},\mathcal{W}_r(R_x+R,0)\ra \leq  C r^{-1/2} (\ee^{-\frac{(R_x)^2}{2r}}).
 \end{equation}
 For fixed $y\in \R$ the map
 \begin{equation*}
 	r \mapsto r^{-1/2}\ee^{-\frac{y^2}{2r}}, \quad r >0,
 \end{equation*}
 is maximized for $r = y^2$ and is monotonously increasing for $r \leq y^2$. Recall that $r = s-\bar{T}_{n,k}^x \leq s$ and thus we get the estimate
 \begin{equation}\label{eq:int_est_W_2}
 	I_3^n(s,x) \leq C(|x|^{1/4}+1) \min(s, R_x^2)^{-1/2} \ee^{-\frac{(R_x)^2}{2\min(R_x^2,s)}}.
 \end{equation}
Finally, we estimate $I_2^n(s,x)$.
To this end, we break up the event into several parts that we can handle more easily.
\begin{equation}\label{eq:I_2_1}
	\begin{split}
		&	I_2^n(s,x) \leq \P( \tau_{R_x,n} \leq K^n_s ,K^n_s \leq |x|^{1/4}) \\
		& \leq \P(  U_{k,x}^n \leq R_x+R \text{ for some } k \in\{1, \ldots, K^n_s\}, K^n_s \leq |x|^{1/4} )\\
		& \leq \sum_{m = 1}^{\lceil |x|^{1/4}\rceil}\sum_{k = 1}^{m}\P( U_{n,j}^x > R_x+R, j\in \{1, \ldots k-1\}, U_{n,k}^x \leq R_x+R , K_s^n = m ).
	\end{split}
\end{equation}
Now, we denote $\eta_x = R_x/(|x|^{1/4}+1)  $. Then,
for $1\leq k\leq m$ we get
\begin{equation}\label{eq:I2_intermidate1}
	\begin{split}
		& \P( U_{n,j}^x > R_x+R, j\in \{1, \ldots k-1\}, U_{n,k}^x \leq R_x+R, K_s^n = m )\\
		& \leq \P( U_{n,k}^x \leq R_x+R,U_{n,j}^x > x-j\eta_x \text{ for all } j \in \{1, \ldots k-1\} , K_s^n = m) \\
		& \qquad  + \P(U_{n,j}^x \leq x-j\eta_x \text{ for some } j\in \{1, \ldots k-1\}, K_s^n = m).
	\end{split}
\end{equation}
We first estimate the first term on the right hand side of \eqref{eq:I2_intermidate1}. For $1\leq k\leq m \leq |x|^{1/4}+1$ we get
\begin{equation*}
	\begin{split}
		&\P( U_{n,k}^x \leq R_x+R,U_{n,j}^x > x-j\eta_x \text{ for all } j \in \{1, \ldots ,k-1\} , K_s^n = m) \\
		&\leq \P( U_{n,k}^x \leq R_x+R,U_{n,j}^x > x-j\eta_x \text{ for all } j \in \{1, \ldots, k-1\} , \bar{T}_{n,k}^x \leq s) \\
		& = \E[\mathbbm{1}_{\bar{T}_{k,x}^n \leq s}\mathbbm{1}_{U_{n,j}^x > x-j\eta_x , j\in \{1, \ldots, k-1\}} \mathbbm{1}_{ U_{n,\ell}^x  \leq R+R_x}]\\
		& =  \E[\mathbbm{1}_{\bar{T}_{n,k}^x \leq s}\mathbbm{1}_{U_{n,j}^x > x-j\eta_x , j\in \{1, \ldots, k-1\}} \E[\mathbbm{1}_{ U_{n,\ell}^x  \leq R+R_x}| \mathcal{G}_{n,x,k}]], 
	\end{split}
\end{equation*}
where $\mathcal{G}_{n,x,k}$ is the $\sigma$-algebra generated by $(U_{n,1}^x, \ldots U_{n,k-1}^x, T_{n,1}^x, \ldots T_{n,k}^x)$ . 
For $y\in \R$, using Lemma~\ref{lem:pde_est}, Proposition~\ref{prp:initial_trace_prop} and the construction of the process $Y(n\delta_x)$ we get that
\begin{equation}\label{eq:U_dist_thing}
	\begin{split}
		&\mathbbm{1}_{\bar{T}_{n,k}^x \leq s} \E[\mathbbm{1}_{ U_{n,k}^x  \leq y}| \mathcal{G}_{n,x,k}]\\
		& =\mathbbm{1}_{\bar{T}_{n,k}^x \leq s} \frac{\int_{-\infty}^{y} Y_{\bar{T}_{n,k}^x-}(n\delta_x)(y')dy'}{\int_\R Y_{\bar{T}_{n,k}^x-}(n\delta_x)(y')dy'}\leq C \sum_{j=0}^{k-1} \int_{-\infty}^{y} \mathcal{W}_{\bar{T}_{n,k}^x-\bar{T}_{n,j}^x}(U_{n,j}^x,0)(y')dy' ),
	\end{split}
\end{equation}
where we also used that
\begin{equation}\label{eq:lower_bound ineq}
	\inf_{n \geq 1}\int_\R Y_{\bar{T}_{n,k}^x-}(n\delta_x)(y')dy' \geq \int_\R V_s(1\delta_x)(y')dy' >c>0,
\end{equation}
for $\bar{T}_{n,k}^x \leq s$. Note that we again used Lemma~\ref{lem:pde_est}, Proposition~\ref{prp:initial_trace_prop} and the construction of $Y(n\delta_x)$ as well as \eqref{eq:monotonicity2} in order to derive \eqref{eq:lower_bound ineq}.
Using \eqref{eq:U_dist_thing} and Lemma~\ref{lem:integral_ineq}, we thus get that 
\begin{equation*}
	\begin{split}
		&\E[\mathbbm{1}_{ U_{n,k}^x  \leq R+R_x}| \mathcal{G}_{n,x,k}] \leq C \sum_{j=0}^{k-1} \int_{-\infty}^{R+R_x} \mathcal{W}_{\bar{T}_{n,k}^x-\bar{T}_{n,j}^x}(U_{n,j}^x,0)(y')dy' \\
		& \leq  C\sum_{j=0}^{k-1} \int_{-\infty}^{R+R_x} \mathcal{W}_{\bar{T}_{n,k}^x-\bar{T}_{j,x}^n}(x-(k-1)\eta_x,0)(y')dy' \\
		&\leq Ck\min(s,\eta_x^2)^{-1/2} \ee^{-\frac{(\eta_x)^2}{2\min(s,\eta_x^2)}}, 
	\end{split}
\end{equation*}
on the event $\{ \bar{T}_{n,k}^x \leq s, U_{n,j}^x > x-j\eta_x  j\in \{1, \ldots k-1\}\}$. Hereby, we again used Proposition~\ref{prp:very_sing_asymp} and the fact that $x-(k-1)\eta_x-R-R_x \geq  \eta_x >c>0$ for all $k \leq \lceil |x|^{1/4} \rceil$
to estimate the integral in the last inequality in the same way as in \eqref{eq:int_est_W_1} and \eqref{eq:int_est_W_2}. Thus, we arrive at 
\begin{equation}\label{eq:I2_intermediate3}
	\begin{split}
	&\P( U_{n,k}^x \leq R_x+R,U_{n,j}^x > x-j\eta_x \text{ for all } j \in \{1, \ldots, k-1\} , K_s^n = m) \\
	&\leq Ck\min(s,\eta_x^2)^{-1/2} \ee^{-\frac{(\eta_x)^2}{2\min(s,\eta_x^2)}}.
	\end{split}
\end{equation}
Now, we estimate the second term on the right hand side of \eqref{eq:I2_intermidate1}. First we can estimate
\begin{equation}\label{eq:I2_intermidate2}
	\begin{split}
		&\P(U_{n,j}^x \leq x-j\eta_x \text{ for some } j\in \{1, \ldots, k-1\}, K_s^n = m) \\
		& \leq \sum_{\ell = 1}^{k-1}  \P(U_{n,j}^x > x-j\eta_x  ,j\in \{1, \ldots, \ell-1\}, U_{n,\ell}^x  \leq x-\ell\eta_x, \bar{T}_{n,\ell}^x \leq s) ,
	\end{split}
\end{equation}
for $1\leq k \leq m \leq |x|^{1/4}+1$.
Similarly to above we get for $1\leq \ell \leq |x|^{1/4}+1$ that
\begin{equation*}
	\begin{split}
		&  \P(U_{n,j}^x > x-j\eta_x,  j\in \{1, \ldots, \ell-1\}, U_{n,\ell}^x  \leq x-\ell\eta_x, \bar{T}_{n,\ell}^x \leq s) \\
		& = \E[\mathbbm{1}_{\bar{T}_{n,\ell}^x \leq s}\mathbbm{1}_{U_{n,j}^x > x-j\eta_x  ,j\in \{1, \ldots, \ell-1\}} \mathbbm{1}_{ U_{n,\ell}^x  \leq x-\ell\eta_x}]\\
		& =  \E[\mathbbm{1}_{\bar{T}_{n,\ell}^x \leq s}\mathbbm{1}_{U_{n,j}^x > x-j\eta_x, j\in \{1, \ldots, \ell-1\}} \E[\mathbbm{1}_{ U_{n,\ell}^x  \leq x-\ell\eta_x}| \mathcal{G}_{n,x,\ell}]]\\
		& \leq  C\ell\min(s,\eta_x^2)^{-1/2} \ee^{-\frac{\eta_x^2}{\min(s,\eta_x^2)}},
	\end{split}
\end{equation*}
and thus 
\begin{equation}\label{eq:I_2_intermediate4}
	\begin{split}
	 &\P(U_{n,j}^x \leq x-j\eta_x \text{ for some } j\in \{1, \ldots, k-1\}, K_s^n = m)\\
	 & \leq Ck^2\min(s,\eta_x^2)^{-1/2} \ee^{-\frac{\eta_x^2}{\min(s,\eta_x^2)}}, \quad 1\leq k \leq m \leq |x|^{1/4}+1.
	 \end{split}
\end{equation}
Together, \eqref{eq:I_2_1}, \eqref{eq:I2_intermidate1}, \eqref{eq:I2_intermediate3} and \eqref{eq:I_2_intermediate4} yield
\begin{equation}\label{eq:I2}
	I_2(s,x) \leq C (|x|^{1/4}+1)^4 \min(s,\eta_x^2)^{-1/2}\ee^{-\frac{(\eta_x)^2}{2\min(s,\eta_x^2)}}.
\end{equation}
Plugging \eqref{eq:I1}, \eqref{eq:I3},\eqref{eq:I2} and \eqref{eq:decomp_comapct} into \eqref{eq:compact_initial} yields
\begin{equation*}
	\begin{split}
		&\lim_{n \rightarrow \infty} \int_\R \big(1-\E[\ee^{-\la X_0, Y_s(n\delta_x)\ra}]\big)dx\\
		& \leq \lim_{n \rightarrow \infty}\int_{B(0,2R)^c} I^n_1(s,x) + I^n_2(s,x) + I^n_3(s,x)dx + 4R\\
		& \leq C\int_{B(0,2R)^c} x^{-2}\wedge1 dx+ C\int_{B(0,2R)^c} \min\Big(1, (|x|^{1/4}+1) \min(s,R_x^2)^{-1/2} \ee^{-\frac{(R_x)^2}{2\min(s,R_x^2)}} \Big)dx\\
		& \quad + C \int_{B(0,2R)^c} \min\Big(1,(|x|^{1/4}+1)^4\min(s,\eta_x^2)^{-1/2} \ee^{-\frac{(\eta_x)^2}{2\min(s,\eta_x^2)}}\Big)dx + 4R.
	\end{split}
\end{equation*}
Note that there exists $R^\ast \geq 2R$ such that
\begin{equation*}
	R_x^2 \geq s, \eta_x^2 \geq s, \quad x \geq R^\ast.
\end{equation*}
Thus, we get 
\begin{equation*}
	\begin{split}
		&\lim_{n \rightarrow \infty} \int_\R \big(1-\E[\ee^{-\la X_0, Y_s(n\delta_x)\ra}]\big)dx\\
		& \leq C\Big(R^\ast + s^{-1/2}\int_{R^\ast}^\infty x\ee^{-\frac{(x-R)^2}{8s|x|^{1/2}}}dx \Big) ,< \infty. 
	\end{split}
\end{equation*}
where we used that $\eta_x \leq R_x$ for $x\geq R^\ast$.
This proves \eqref{eq:th_supp}.
\end{proof}

\section{Proof of Theorem~\ref{thm:survival}}\label{sec:survival}
The proof of Theorem~\ref{thm:survival} relies on a rescaling of \eqref{eq:SPDE1} as well as a comparison with solutions to
\begin{equation}\label{eq:FKKP}
\begin{split}
	d_t\bar{X}^\theta_t(x) &= \frac{1}{6}\Delta \bar{X}^\theta_t(x) + \theta \bar{X}^\theta_t(x) - \big(\bar{X}^\theta_t(x)\big)^2 + \sqrt{\bar{X}^\theta_t(x)}\dot{W}(t,x), \quad t > 0, \quad x\in \R,\\
	\bar{X}^\theta_0(x) &= f(x), \quad x\in \R,
\end{split}
\end{equation}
where $\theta>0$ and $f\in \mathcal{C}_c^+$. 
\begin{lemma}\label{lem:rescaling}
	Let $a,b,c>0$ and let $X$ be a solution to \eqref{eq:SPDE1}. Then 
	\begin{equation*}
		\tilde{X}_t(x) = cX_{at}(bx), \quad t \geq 0, \quad x\in \R, 
	\end{equation*}
	is a solution to 
	\begin{equation*}
	\begin{split}
		d_t \tilde{X}_t(x) &= \frac{a}{2b^2}\Delta\tilde{X}_t(x) + ach_c(\tilde{X}_t(x)) + \frac{c^{1/2}a^{1/2}}{b^{1/2}}\sqrt{\tilde{X}_t(x)} \dot{W}(t,x),\quad x\in \R ,t\geq 0,\\
	\tilde{X}_0(x) &= cX_0(bx), \quad x\in \R ,
		\end{split}
	\end{equation*}
	where $h_c(x) = h(c^{-1}x)$ for $x\in \R$.
\end{lemma}
	\begin{proof}
		This follows along the same lines as the proof of \cite[Lemma 2.1.2]{Muller94}.
	\end{proof}

\begin{lemma}\label{lem:firststep_survival}
	Let 
	\begin{equation*}
		h(x) =\theta( 1-e^{-\lambda x})
	\end{equation*}
	for some $\lambda, \theta >0$ and let $X$ be a solution to \eqref{eq:SPDE1} with $X_0 =  f\in \mathcal{C}_c^+$ such that $f\neq 0$. 
	Then it follows that 
	\begin{equation*}
		\P(\la X_t,1 \ra>0, \forall t > 0)>0. 
	\end{equation*}
\end{lemma}
\begin{proof}
	By \cite[Theorem 1]{Muller94} there exists $\theta^\ast >0$ such that 
	\begin{equation*}
		\P(\la \bar{X}^{\theta^\ast}_t,1\ra >0, \forall t >0) >0.
	\end{equation*}
	Now, we set $a = 1/3b^2$, and $c = 3b^{-1}$ for $b>0$. Then, by Lemma~\ref{lem:rescaling}, $\tilde{X}_t(x) = cX_{at}(bx)$ solves
	\begin{equation}\label{eq:rescaled}
		d_t \tilde{X}_t = \frac{1}{6} \Delta \tilde{X}_t(x) + b\theta(1-\ee^{-\frac{b}{3} \lambda \tilde{X}_t(x)})+ \sqrt{\tilde{X}_t(x)} \dot{W}(t,x), \quad x\in \R, t\geq 0.
	\end{equation}
	Suppose, that there exists $b>0$ such that 
	\begin{equation}\label{eq:domination_h_1}
		b\theta(1-\ee^{-\frac{b}{3} \lambda x}) \geq \theta^\ast x -x^2, \quad \forall x\geq 0.  
	\end{equation}
	By \cite[Theorem 3.3]{Kotelenez93} and standard approximation and tightness arguments it follows by weak uniqueness that $\tilde{X}$ and $\bar{X}$ can be defined on a common probability space such that  for all $t\geq 0$ we have
	\begin{equation*}
		\tilde{X}_{t} (x) \geq \bar{X}^{\theta^\ast}_{t}(x), \quad x\in \R,
	\end{equation*}
	almost surely. Again by weak uniqueness and by continuity of $\tilde{X}$ and $\bar{X}$ we  get 
	\begin{equation*}
		\P(\tilde{X}_t(x) \geq \bar{X}^{\theta^\ast}_t(x), \forall t\geq 0, x\in \R)=1.
	\end{equation*}
	This in particular implies that $\P\la (\tilde{X}_t,1\ra >0 ,\forall t >0) >0$ and thus $\P(\la X_t,1\ra >0, \forall t >0) >0$.
	Thus, it only remains to establish \eqref{eq:domination_h_1}, which is equivalent to proving 
	\begin{equation*}
		\frac{b^3 \lambda ^2}{9}\theta(1-\ee^{-x}) \geq \theta^\ast\frac{b \lambda }{3} x-  x^2, \quad \forall x \geq 0,
	\end{equation*}
	for some $b>0$.
	For $x\in [0,1]$ we may use that $(1-\ee^{-x}) \geq \frac{1}{2} x$ and thus in this regime it is sufficient to choose $b>0$ sufficiently large so that 
	\begin{equation*}
		\frac{b^3 \lambda ^2}{9}\theta \frac{1}{2} \geq \theta^\ast\frac{b \lambda }{3}.
	\end{equation*}
	Furthermore, note that 
	\begin{equation*}
		\Big(\theta^\ast\frac{b \lambda }{3} \Big)^2 \frac{1}{4} \geq   \theta^\ast\frac{b \lambda }{3} x-  x^2,
	\end{equation*}
	for $x\in\R$. Thus, for $x>1$,  we need to choose $b>0$  sufficiently large such that also 
	\begin{equation*}
			\frac{b^3 \lambda ^2}{9}\theta(1-\ee^{-1}) \geq 	\Big(\theta^\ast\frac{b \lambda }{3} \Big)^2 \frac{1}{4} .
	\end{equation*}
	That means, \eqref{eq:domination_h_1} is satisfied for 
	\begin{equation*}
		b> \max\Big( \sqrt{\frac{6\theta^\ast}{\lambda \theta}}, \frac{(\theta^\ast)^2}{4\theta (1-e^{-1})} \Big), 
	\end{equation*}
	which finishes the proof.
\end{proof}
\begin{proof}[Proof of Theorem~\ref{thm:survival}]
	Suppose that there exists $\lambda, \theta >0$ such that 
	\begin{equation}\label{eq:domination_h}
		h(x) \geq \theta (1-\ee^{-\lambda x}) ,\quad \forall x\geq 0. 
	\end{equation}
	Then, by weak uniqueness, \cite[Corollary 2.4]{Shiga94} and standard approximation and tightness arguments we have that 
	\begin{equation*}
		\P(\la X_t,1 \ra > 0, \forall t>0) \geq	\P(\la \hat{X}_t, 1\ra  > 0, \forall t>0)
	\end{equation*}
	where $\hat{X}$ is the solution to \eqref{eq:SPDE1} with $h (x) = \theta (1-\ee^{-\lambda x}) $ for all $x\in \R$. Thus, by Lemma~\ref{lem:firststep_survival} it is sufficient to verify \eqref{eq:domination_h}.  If $b_1>0$, then \eqref{eq:domination_h} is clearly satisfied for all $\lambda >0$ and $\theta \leq b_1$. When $b_1 = 0$ we have the estimate
	\begin{equation*}
		\int_0^\infty (1-\ee^{-\tilde{\lambda} x})\nu^1(d\tilde{\lambda}) \geq   (1-\ee^{-\lambda x}) \nu^1([\lambda, \infty)), 
	\end{equation*}
	where $\lambda >0$ is chosen in such a way that $\nu^1([\lambda, \infty)) >0$. Then, \eqref{eq:domination_h} is clearly satisfied for $\theta \leq \nu^1([\lambda, \infty))$.
\end{proof}

\bibliography{bibliography_m.bib}

\end{document}